\documentclass[11pt, a4paper]{article}

\usepackage[utf8]{inputenc}
\usepackage[english]{babel}

\usepackage{amsmath}
\usepackage{amssymb}
\usepackage{amsthm}
\usepackage{mathabx}
\usepackage{tikz}
\usepackage{tikz-cd}
\usepackage{todonotes}
\usepackage{hyperref}
\usepackage{cleveref}
\usepackage{footmisc}
\usepackage{bbm}
\usepackage{afterpage}
\usepackage{mathrsfs}
\usepackage{sseq}
\usepackage{extarrows}
\usepackage{enumitem}
\usepackage{stmaryrd}
\usepackage{pdflscape}

\newcommand{\alg}{\mathrm{alg}}
\newcommand{\aff}{\mathrm{aff}}

\newcommand{\co}{{\mathrm{co}}}
\newcommand{\cont}{{\mathrm{cont}}}
\newcommand{\Cpl}{{\mathrm{Cpl}}}
\DeclareMathOperator{\Def}{Def}

\newcommand{\dR}{\mathrm{dR}}
\DeclareMathOperator{\Endo}{End}

\DeclareMathOperator{\fib}{fib}

\newcommand{\ori}{\mathrm{ori}}
\newcommand{\lN}{{\mathrm{loc.N}}}
\DeclareMathOperator{\Lie}{Lie}

\newcommand{\loc}{\mathrm{loc}}
\newcommand{\m}{{\mathfrak{m}}}

\newcommand{\red}{{\mathrm{red}}}
\newcommand{\sHen}{\mathrm{sHen}}
\DeclareMathOperator{\Spec}{Spec}
\DeclareMathOperator{\Spf}{Spf}

\newcommand{\et}{\mathrm{\acute{e}t}}

\newcommand{\Sph}{\mathbf{S}}

\newcommand{\KU}{{\mathrm{KU}}}
\DeclareMathOperator{\KO}{KO}

\DeclareMathOperator{\MUP}{MUP}
\DeclareMathOperator{\TMF}{TMF}
\DeclareMathOperator{\TAF}{TAF}
\DeclareMathOperator{\Tmf}{Tmf}
\DeclareMathOperator{\tmf}{tmf}

\newcommand{\ad}{{\mathrm{ad}}}
\newcommand{\cn}{{\mathrm{cn}}}
\newcommand{\acn}{{\mathrm{acn}}}
\newcommand{\nc}{{\mathrm{nc}}}
\renewcommand{\o}{{\mathrm{or}}}
\newcommand{\un}{{\mathrm{un}}}

\newcommand{\1}{\mathbf{1}}

\DeclareMathOperator{\Ab}{Ab}

\DeclareMathOperator{\Aff}{Aff}

\newcommand{\BT}{{\mathrm{BT}^p}}
\newcommand{\BTn}{n}
\newcommand{\BTone}{1}
\newcommand{\BTtwo}{2}
\newcommand{\BTtwog}{{2g}}
\DeclareMathOperator{\Cat}{\mathscr{C}at}

\newcommand{\CAlg}{{\mathrm{CAlg}}}
\DeclareMathOperator{\CohomTh}{CohomTh}

\newcommand{\Ell}{{\mathrm{Ell}}}
\DeclareMathOperator{\Fun}{Fun}

\newcommand{\OrFGroup}{{\mathrm{OrFG}}}
\newcommand{\FGroup}{{\mathrm{FG}}}

\newcommand{\LT}{{\mathrm{LT}}}
\newcommand{\Mod}{\mathrm{Mod}}

\newcommand{\Nil}{{\mathrm{Nil}}}
\DeclareMathOperator{\OrDat}{OrDat}
\newcommand{\OrBT}{{\mathrm{OrBT}^p}}
\newcommand{\OrBTn}{{\mathrm{OrBT}^p_n}}

\DeclareMathOperator{\QCoh}{QCoh}

\DeclareMathOperator{\SpDM}{SpDM}
\DeclareMathOperator{\fSpDM}{fSpDM}
\DeclareMathOperator{\fDM}{fDM}
\DeclareMathOperator{\DM}{DM}
\newcommand{\Shv}{\mathcal{S}\mathrm{hv}}
\newcommand{\Sp}{\mathrm{Sp}}
\newcommand{\Spc}{\mathcal{S}}
\newcommand{\Set}{\mathrm{Set}}
\renewcommand{\top}{\mathrm{top}}
\newcommand{\Top}{\mathcal{T}\mathrm{op}}

\DeclareMathOperator{\AVar}{AVar}

\renewcommand{\lim}{{\mathrm{lim}\,}}
\DeclareMathOperator{\colim}{colim}

\newcommand{\Hom}{\mathrm{Hom}}

\DeclareMathOperator{\Map}{Map}

\newcommand{\id}{\mathrm{id}}

\newcommand{\op}{\mathrm{op}}

\newcommand{\SL}{\mathbf{SL}}
\newcommand{\Z}{\mathbf{Z}}

\DeclareMathOperator{\Ext}{Ext}

\newcommand{\C}{\mathcal{C}}

\newcommand{\complexproj}{\mathbf{CP}}
\newcommand{\CP}{\mathbf{CP}}

\newcommand{\D}{\mathscr{D}}

\newcommand{\E}{\mathbf{E}}
\newcommand{\EE}{\mathcal{E}}
\newcommand{\F}{\mathbf{F}}
\newcommand{\FF}{\mathcal{F}}
\newcommand{\G}{\mathbf{G}}
\newcommand{\GG}{\mathcal{G}}

\newcommand{\h}{\mathrm{h}}

\renewcommand{\L}{\mathcal{L}}

\newcommand{\M}{\mathcal{M}}

\renewcommand{\O}{\mathscr{O}}
\newcommand{\OO}{\mathfrak{O}}
\renewcommand{\P}{\mathcal{P}}

\newcommand{\Q}{\mathbf{Q}}
\newcommand{\QQ}{\mathcal{Q}}

\newcommand{\T}{\mathcal{T}}

\renewcommand{\u}{\mathfrak{U}}
\newcommand{\U}{\mathcal{U}}

\newcommand{\X}{\mathcal{X}}

\newcommand{\x}{\mathfrak{X}}
\newcommand{\xx}{\mathsf{X}}
\newcommand{\Y}{\mathcal{Y}}
\newcommand{\y}{\mathfrak{Y}}

\newcommand{\yy}{\mathsf{Y}}

\newcommand{\al}{\alpha}
\newcommand{\be}{\beta}

\newcommand{\Ga}{\Gamma}

\renewcommand{\k}{{\kappa}}

\theoremstyle{theorem}\numberwithin{equation}{section}
\newtheorem{theorem}[equation]{Theorem}
\crefname{theorem}{{th}.\!\!}{{ths}.\!\!}
\Crefname{theorem}{{Th}.\!\!}{{Ths}.\!\!}

\theoremstyle{luriestheorem}

\newtheorem*{luriestheorem*}{Lurie's theorem}

\crefname{theoremalph}{{th}.\!\!}{{ths}.\!\!}
\Crefname{theoremalph}{{Th}.\!\!}{{Ths}.\!\!}

\Crefname{problem}{{Prb}.\!\!}{{Prbs}.\!\!}
\newtheorem{prop}[equation]{Proposition}
\Crefname{prop}{{Pr}.\!\!}{{Prs}.\!\!}

\Crefname{lemma}{{Lm}.\!\!}{{Lms}.\!\!}
\newtheorem{cor}[equation]{Corollary}
\Crefname{cor}{{Cor}.\!\!}{{Cors}.\!\!}

\theoremstyle{definition}\numberwithin{equation}{section}
\newtheorem{mydef}[equation]{Definition}
\Crefname{mydef}{{Df}.\!\!}{{Dfs}.\!\!}

\Crefname{recall}{{Rcl}.\!\!}{{Rcls}.\!\!}
\newtheorem{construction}[equation]{Construction}
\Crefname{construction}{{Con}.\!\!}{{Cons}.\!\!}

\Crefname{ass}{{As}.\!\!}{{As}.\!\!}
\newtheorem{notation}[equation]{Notation}
\Crefname{notation}{{Nt}.\!\!}{{Nts}.\!\!}

\Crefname{situation}{{St}.\!\!}{{Sts}.\!\!}

\theoremstyle{remark}\numberwithin{equation}{section}
\newtheorem{example}[equation]{Example}
\Crefname{example}{{Ex}.\!\!}{{Exs}.\!\!}

\Crefname{nonexample}{{NonEx}.\!\!}{{NonEx}.\!\!}

\Crefname{claim}{{Clm}.\!\!}{{Clms}.\!\!}
\newtheorem{remark}[equation]{Remark}
\Crefname{remark}{{Rmk}.\!\!}{{Rmks}.\!\!}

\Crefname{idea}{{Id}.\!\!}{{Ids}.\!\!}
\newtheorem{warn}[equation]{Warning}
\Crefname{warn}{{Warn}.\!\!}{{Warns}.\!\!}

\Crefname{section}{{\textsection}\!\!}{{\textsection}\!\!}
\Crefname{subsection}{{\textsection}\!\!}{{\textsection}\!\!}
\Crefname{appendix}{{\textsection}\!\!}{{\textsection}\!\!}

\begin{document}
\title{On Lurie's theorem and applications}
\author{Jack Morgan Davies\footnote{Mathematisches Institut, Universität Bonn, Germany, \url{davies@math.uni-bonn.de}}}
\maketitle

\begin{abstract}
Lurie's theorem states that there exists a sheaf of ring spectra on the site of formally étale Deligne--Mumford stacks over the moduli stack of $p$-divisible groups of height $n$, which agrees with the classical Landweber exact functor theorem (LEFT) on affines. In other words, this theorem is a global, higher categorical refinement of the LEFT. In recent work, Lurie has introduced many of the ingredients one needs to prove this theorem, and in this article, we gather these ingredients together and prove Lurie's theorem. Applications of this theorem to Lubin--Tate theories, topological modular and automorphism forms, and Adams operations are also discussed.

\end{abstract}


\setcounter{tocdepth}{3}
{\tableofcontents}


\addcontentsline{toc}{section}{Introduction}
\section*{Introduction} 

One of Lurie's original motivations to develop spectral algebraic geometry is to give a moduli interpretation of the theorems of Goerss--Hopkins--Miller producing Lubin--Tate theories and topological modular forms. This ambition is outlined in Lurie's survey on the subject \cite{lurieecsurveyname}. Various authors, see \cite[\textsection8.1]{taf} for example, stated that Lurie had a proof of a statement along the following lines---everything is implicitly $p$-completed.

\begin{luriestheorem*}[Version 1]
	For each formally étale morphism $\x\to \M_\BTn$ from a Deligne--Mumford stack into the moduli stack of $p$-divisible groups of height $n$, the structure sheaf $\O_\x$ has a canonical lift to a sheaf of locally Landweber ring spectra $\O^\top_\x$.
\end{luriestheorem*}

The goal of this article is to provide an accurate statement and a complete proof of this theorem, as well as describe some applications, including the theorems of Goerss--Hopkins--Miller above. To contextualise the statement of Lurie's theorem and give a more precise statement, let us recall Quillen's theorem and the birth of elliptic cohomology.\\

A homotopy commutative ring spectrum $E$ is said to be \emph{even $2$-periodic} if $\pi_\ast E$ contains a unit of degree 2 and $\pi_{2n+1}E=0$ for all $n\in \Z$. Given such an $E$, the formal spectrum $\Spf (E^0 \complexproj^\infty)$ has the structure of a (smooth $1$-parametre) \emph{formal group} coming from the multiplicative structure of $\complexproj^\infty$. Quillen's theorem (\cite{quillenstheorem} is the remarkable fact that the homotopy groups of the even $2$-periodic complex cobordism spectrum $\MUP$ classify formal groups, or more precisely, that the stack associated with the Hopf algebroid $(\MUP_0, \MUP_0\MUP)$ can be identified with the \emph{moduli stack of formal groups} $\M_{\FGroup}$. This passage, from algebraic topology to arithmetic geometry has proven useful in homotopy theory, especially in providing global structure to the stable homotopy category; see \cite{bluebook,agnestobi,orangebook,greenbook} for an overview and \cite{MRW77,moravacobordism} for early applications from this perspective. One can also use this connection to construct cohomology theories from certain formal groups.\\

Given an elliptic curve $\xx$ over a ring $R$, the completion of $\xx$ at its identity section $\widehat{\xx}$ (\cite[\textsection IV]{silvermaneasy}) is a formal group. In this situation, write $f_\xx\colon \Spec R \to \M_{\FGroup}$ for the map defined by $\widehat{\xx}$. The \emph{Landweber exact functor theorem} (LEFT) states that if $f_\xx$ is flat, then the assignment sending a spectrum $Y$ to the pullback of the quasi-coherent sheaf $\MUP_\ast Y$ along $f_\xx$ defines a homology theory. Applying Brown representability, this homology theory produces a spectrum $E_\xx$. Furthermore, a theorem of Hovey--Strickland \cite[Cor.2.15]{hoveysticklandMoravalocal} states that the passage from homology theories to spectra is faithful when applied after the LEFT. This means that $E_\xx$ comes equipped with a homotopy commutative multiplication giving us an identification of rings $\pi_0 E_\xx = R$. The LEFT also shows that $E_\xx$ is even $2$-periodic and identifies the formal groups $\Spf (E_\xx^0 \complexproj^\infty) = \widehat{\xx}$. In other words, $E_\xx$ is an \emph{elliptic cohomology theory} for $\xx$; see \cite{landweberconferencebook} for the classical statement of the LEFT, or \cite[\textsection5]{naumann07} and \cite[Th.0.0.1]{ec2} for this modern reformulation. A ring spectrum $E$ that can be constructed using the LEFT is said to be \emph{Landweber exact}.\\

The above produces elliptic cohomology theories, but only over affine schemes $\Spec R$. It is also hopeless to try and glue together elliptic cohomology theories using the LEFT, as the desired diagrams of spectra lie in the \emph{homotopy category of spectra} $\h\Sp$, which famously has few limits and colimits. We want to therefore replace $\h\Sp$ with the $\infty$-category $\CAlg$ of \emph{$\E_\infty$-ring spectra}, or simply \emph{$\E_\infty$-rings}, the homotopy coherent version of homotopy commutative ring spectra. The advantage of this $\infty$-category $\CAlg$ is that it has all limits and colimits, now meant in an $\infty$-categorical sense.\footnote{Classically this was phrased with model categories. Another facet of Lurie's wider project was to develop higher algebra and spectral algebraic geometry purely in the language of $\infty$-categories.} The \emph{Goerss--Hopkins--Miller theorem} (GHMT, see \cite[Th.1.2]{bourbakigoerss}) is then a partial lift of the LEFT to the $\infty$-category $\CAlg$. More precisely, the GHMT states that the LEFT can be refined to a functor $\O^\top$ from the category étale maps $\Spec R \to \M_\Ell$ to the \emph{moduli stack of elliptic curves} $\M_\Ell$ to the $\infty$-category $\CAlg$. This functor was originally defined on the compactification $\overline{\M}_\Ell$, but this article will not immediately apply to this scenario; see \Cref{nocompactifcation}. The GHMT also constructs $\O^\top$ as an \emph{étale sheaf}, meaning that one can recover the value of $\O^\top$ on a Deligne--Mumford stack $S$ which is étale over $\M_\Ell$ by taking a limit of the Landweber exact values of $\O^\top$ on an affine cover of $S$. For example, one can apply this to the moduli stack of elliptic curves itself $\M_\Ell$, from which one obtains the universal theory elliptic cohomology theory $\TMF=\O^\top(\M_\Ell)$ of \emph{topological modular forms}. This spectrum, and its connective form $\tmf$, had been a long time coming. In \cite{hopkinsfirsttmficm}, Hopkins discusses his discovery of $\tmf$, there written as $eo_2$, together with Mahowald and Miller, as well as connections to differential geometry through the Witten genus. This was originally done using $\E_1$-rings, there call $\mathbf{A}_\infty$-ring spectra, and the further refinement to $\E_\infty$-rings found in the GHMT was only made possible due to the obstruction theory of Goerss--Hopkins \cite{gh04}. The first published construction of $\tmf$ and $\TMF$ appears in \cite[\textsection12]{tmfbook}; see \cite[\textsection6]{handbooktmf} for a further discussion.\\

Lurie asked if the functor $\O^\top$ in the GHMT, whose original construction uses technical obstruction theory, could be recovered more naturally as the structure sheaf of a moduli problem in spectral algebraic geometry. This ambition is outlined throughout \cite{lurieecsurveyname} and carried out in detail in \cite[\textsection7]{ec2}. The answer is: $\O^\top$ is the structure sheaf for the \emph{moduli stack of oriented elliptic curves}, an object naturally defined using spectral algebraic geometry. Aside from being an alternative construction, the methods of proof suggest a vast generalisation of the GHMT from elliptic curves to \emph{$p$-divisible groups}. The following is often referred to as \emph{Lurie's theorem}, which first appeared without proof in \cite[Th.8.1.4]{taf}. We have still omitted some technical details here; see \Cref{maintheorem} for the precise statement.

\begin{luriestheorem*}[Version 2]
Fix a prime $p$ and an integer $n\geq 1$. There is an étale sheaf of $\E_\infty$-rings $\O^\top_\BTn$ on the site of Deligne--Mumford stacks admitting a formally étale morphism to the moduli stack $\M_\BTn$ of $p$-divisible groups of height $n$ such that its value on a $p$-divisible group $\G$ over a ring $R$ is an $\E_\infty$-ring $\EE$ with the following properties:
\begin{enumerate}
\item $\EE$ has a canonical complex orientation and is Landweber exact.
\item There is a canonical isomorphism of rings $\pi_0 \EE\simeq R$.
\item The homotopy groups $\pi_\ast\EE$ vanish for all odd integers and otherwise $\pi_{2k} \EE$ canonically isomorphic to the $k$-fold tensor product of the line bundle of invariant differentials of the formal group $\widehat{\G}$ associated with $\G$.
\item There is a canonical isomorphism of formal groups between $\widehat{\G}$ and $\Spf (\EE^0 \complexproj^\infty)$.
\end{enumerate}
\end{luriestheorem*}

In some sense, this sheaf $\O^\top_\BTn$ gives a kind of structure sheaf to a moduli stack of oriented $p$-divisible groups, although we will not make this intuition precise.\\

This relates to Version 1 as for each Deligne--Mumford stack $S$ equipped with a formally étale morphism to $\M_\BTn$, the above theorem associates with each étale morphism $\Spec R \to S$ an $\E_\infty$-ring $\EE$ with $\pi_0 \EE \simeq R$. In other words, this assignment refines the structure sheaf of $S$ to a sheaf $\O^\top_S$ of locally Landweber $\E_\infty$-rings. The pair $(S,\O^\top_S)$ is a \emph{spectral Deligne--Mumford stack}. There is also an implicit functoriality in Lurie's theorem: given a morphism $f\colon T \to S$ between two such Deligne--Mumford stacks over $\M_\BTn$, then we obtain a morphism between the associated spectral Deligne--Mumford stacks $(T, \O^\top_T) \to (S, \O^\top_S)$, and this assignment is functorial.\\

Interest in this theorem stems from its applications, all originally due to Lurie \cite{ec2} and Behrens--Lawson \cite{taf}, which we discuss in \Cref{applicationsandactions}:

\begin{itemize}
\item All of the \emph{Lubin--Tate theories} arising from a perfect field $\kappa$ and a formal group $\widehat{\G}$ of exact height $n$ over $\kappa$ can be recovered from $\O^\top_\BTn$; see \Cref{ltbttheoryeiessection}. The functoriality of this construction with respect to automorphisms of formal groups recovers the action of the \emph{extended Morava stabiliser group} on such $\E_\infty$-rings, as studied in \cite{gh04}.
\item One can apply Lurie's theorem to the moduli stack of elliptic curves to recover the ($p$-completion of the) sheaf $\O^\top$ and also the $p$-completion of the universal elliptic cohomology theory $\TMF_p$; see \Cref{tmfsection}. 
\item A more in-depth study of dimension $g$ abelian varieties with PEL structure yields cohomology theories called \emph{topological automorphic forms} due to Behrens--Lawson; see \Cref{tafsection}. The only known construction of these theories requires Lurie's theorem.
\end{itemize}

Lurie's theorem also begins to make up for an interesting gap left between the GHMT and the LEFT. Indeed, the LEFT constructs cohomology theories which are functorial in formal groups. In particular, automorphisms of a formal group act upon the associated spectrum in the stable homotopy category. The statement of the GHMT only lifts automorphisms of an elliptic curve to automorphisms on the associated $\E_\infty$-ring in $\CAlg$. Lurie's theorem shows that automorphisms of the underlying $p$-divisible groups of elliptic curves also act on these $\E_\infty$-rings. The simplest nontrivial case of this yields stable Adams operations discussed in (\Cref{adamsoperationsection}). Further applications to topological modular forms appear in \cite{adamsontmf,mythesis,heckeontmf,realspectra,heighttwojatthree}.

\subsection*{Outline}
In the language of \cite{ec2}, for each formally étale map $\G\colon \Spf B \to \M_n$ of classical stacks, we want to define $\O^\top_n(\G)$ as the \emph{orientation classifier} of the \emph{universal spectral deformation} of $\G_0$. In a little more detail, after we define many of the necessary objects and precisely state Lurie's theorem (\Cref{lurietheoremsection}), we discuss the notion of \emph{formally étale} morphisms in spectral algebraic geometry (\Cref{fetsubsection}). Applying this theory to the classical moduli stack of $p$-divisible group of height $n$, we obtain a canonical lift of formally étale morphisms $\G\colon \x \to \M_n$ to formally étale morphisms over the spectral moduli stack of $p$-divisible groups (\Cref{defonmoduliofbdivsection}). In other words, a canonical passage for such objects from classical to spectral algebraic geometry. This produces the \emph{universal spectral deformation} mentioned above. Next, we define a spectral moduli stack $\M^\ori_\BTn$ of \emph{oriented} $p$-divisible groups of height $n$ (\Cref{orientationparaphgra}) à la Lurie's \cite[\textsection4]{ec2}. By definition, base change along the forgetful map from $\M^\ori_\BTn$ to the spectral moduli stack of $p$-divisible groups defines the \emph{orientation classifiers} from above. Finally, we prove that the process given by first taking universal spectral deformations followed by orientation classifiers defines an étale sheaf $\O^\top_\BTn$ which satisfies the conditions of Lurie's theorem (\Cref{finalsectioninproof}). We then apply Lurie's theorem to construct a variety of well-known $\E_\infty$-rings as well as some operations and actions thereupon (\Cref{applicationsandactions}). In an appendix we summarise some technical facts about formal spectral Deligne--Mumford stacks used elsewhere in this article (\Cref{appendix}).

\subsection*{Acknowledgements}
As I have and will continue to mention, the ideas present here are due to Jacob Lurie and can be found spread out among \cite{lurieecsurveyname,ec1,ec2,sag,ec3}---my intention here is to collect those specific ideas needed to prove Lurie's theorem. I would like to thank my PhD supervisor Lennart Meier for suggesting this project and for his immeasurable help. I also had some useful conversations with Alice Hedenlund, Gijs Heuts, Achim Krause, Tommy Lundemo, Magdalena K\c{e}dziorek, Luca Pol, Adrien Sauvaget, and Yuqing Shi, and thank you to Lennart, Tommy, and Theresa Rauch for reading some drafts. This article would also not be what it is without the support of Paul Goerss. Many sincere thanks to the anonymous referees, both of them, for their patience and clarity of suggestions, all of which greatly improved this work. Much of the revising and editing was done while I was in residence at Institut Mittag-Leffler in Djursholm, Sweden during the first quarter of 2022, so this work is supported by the Swedish Research Council under grant no. 2016-06596. At Bonn now, I am an associate member of the Hausdorff Center for Mathematics at the University of Bonn (\texttt{DFG GZ 2047/1}, project ID \texttt{390685813}).


\addcontentsline{toc}{section}{Conventions}
\section*{Conventions}

Now, and forever, fix a prime $p$.

\subsection*{Higher categories and higher algebra}
We will make free and extensive use of the language of $\infty$-categories, higher algebra, and spectral algebraic geometry, following \cite{htt}, \cite{ha}, \cite{sag}, and especially the conventions listed in \cite{ec2}. In particular:
\begin{itemize}
\item For an $\infty$-category $\C$ and two objects $X$ and $Y$ of $\C$, we will write $\Map_\C(X,Y)$ for the mapping space and $\Hom_\C(X,Y)$ if $\C$ happens to be the nerve of a $1$-category.
\item Commutative rings and abelian groups will be treated as discrete $\E_\infty$-rings and spectra. The $\infty$-category of discrete $\E_\infty$-rings is denoted by $\CAlg^\heartsuit$. Moreover, the smash product of spectra will be written as $\otimes$ and the $\infty$-categorical completion functor will be written as $(-)^\wedge_I$ following \cite[\textsection7]{sag}.
\item All module categories $\Mod_R$ refer to the stable $\infty$-category of $R$-modules, where $R$ is an $\E_\infty$-ring. In particular, if $R$ is a discrete commutative ring, then $\Mod_R$ will be the stable $\infty$-category of $R$-module spectra, and \textbf{not} the abelian 1-category of $R$-modules. The same holds for $\infty$-categories of quasi-coherent sheaves.
\item Following \cite{ec2} (and contrary to \cite{ec1,sag}), we will write $\Spec R$ for the nonconnective spectral Deligne--Mumford stack associated with an $\E_\infty$-ring $R$.
\end{itemize}
Moreover, all $n$-categories are $(n,1)$-categories, for $n=1,2,\infty$.

\subsection*{Sites and sheaves}
Lurie's theorem concerns higher categorical sheaves. The $\infty$-categories that we consider as sites are not necessarily (essentially) small, so we a priori do need to be careful about potential size issues. However, we are interested in constructing particular functors and proving they are sheaves, so we only really need to step into a large universal to quantify our definition of a sheaf; see \cite[\textsection6.2.2]{htt} for a wider discussion.

\begin{mydef}[Sheaves]
Let $\T$ be an $\infty$-category $\T$ with finite coproducts and finite limits, $\tau$ be a Grothendieck topology on $\T$, and $\C$ an $\infty$-category $\C$ with small limits. Suppose further that associated with each covering family $\{Y_\al \to X\}$ of an object $X$ in $\T$, there is a finite subfamily $Y_\be$ such that $\coprod Y_\be \to X$ is a $\tau$-cover. Given a $\tau$-cover $Y\to X$ in $\T$, we define its \emph{\v{C}ech nerve} as $Y_\bullet \to X$, so $Y_n = Y^{\times_X (n+1)}$. A functor $F\colon \T^\op\to \C$ is a \emph{$\C$-valued $\tau$-sheaf on $\T$} if for each object $X$ in $\T$ and each $\tau$-cover $Y\to X$ in $\T$, the natural map $F(X) \to \lim_{\Delta} F(Y_\bullet)$ is an equivalence. Let us write $\P_\C(\T) = \Fun(\T^\op, \C)$ for the $\infty$-category of $\C$-valued presheaves on $\T$ and $\Shv_\C(\T)$ for the $\infty$-subcategory of $\tau$-sheaves. If $\C=\Spc$ is the $\infty$-category of spaces, we will drop the $\C$ from this notation.
\end{mydef}

The restriction that each $\tau$-covering family has a finite subcover will not bother us, as in practice our stacks will all be quasi-compact and quasi-separated.\\

A \emph{hypercover} is a generalisation of a cover in a Grothendieck site. In general, our sheaves, including the sheaf occurring in the statement of Lurie's theorem will be hypersheaves. Following \cite[\textsection A]{sag}, this variation on a sheaf comes with a more concrete description.

\begin{mydef}[Hypersheaves {(\cite[Df.A.5.7.1]{sag})}]\label{hypercoversandthelike}
Let $\Delta_{s,+}$ denote the 1-category whose objects are linearly ordered sets of the form $[n]=\{0<1<\cdots<n\}$ for $n\geq -1$, and whose morphisms are strictly increasing functions. We will omit the $+$ when considering the full $\infty$-subcategory with $n\geq 0$. If $\T$ is an $\infty$-category, we will refer to a functor $X_\bullet\colon \Delta_{s,+}^\op\to\T$ as an \emph{augmented semisimplicial object of $\T$}. When $\T$ admits finite limits, then for each $n\geq 0$, we can associate to an augmented semisimplicial object $X_\bullet$ the \emph{$n$th matching object} $M_n(X_\bullet)$ and its associated \emph{matching map}
\[X_n\to \underset{{[i]\hookrightarrow [n]}}{\lim} X_i = M_n(X_\bullet)\]
where the limit above is taken over all injective maps $[i]\hookrightarrow [n]$ such that $i<n$. Given a collection of morphisms $S$ inside $\T$, we call an augmented semisimplicial object $X_\bullet$ an \emph{$S$-hypercover} (for $X_{-1}=X$) if all the natural matching maps belongs to $S$, for every $n\geq 0$. Given a Grothendieck topology $\tau$ on $\T$, then a presheaf of spectra $F$ on $\T$ is called a \emph{$\tau$-hypersheaf} if for all $\tau$-hypercovers $X_\bullet\to X$, the natural map $F(X)\to {\lim}_{\Delta_s} F(X_\bullet)$ is an equivalence of spectra. Some useful general references for the prefix \emph{hyper} in the homotopy theory of sheaves are \cite{clausenmathew1} and \cite[\textsection A]{sag}.
\end{mydef}

Our favourite examples will be when $\T$ is the $\infty$-category $\Aff^\cn$, $\Aff$, $\C_{A_0}$, or $\C_A$, and $S$ is either fpqc or \'{e}tale covers.\\

Given $\T$ and $\tau$ from \Cref{hypercoversandthelike}, then for each $\tau$-cover $Y\to X$ in $\T$, the associated \v{C}ech nerve $Y_\bullet$ is a $\tau$-hypercover. This and the fact that $\Delta_s$ is a final subdiagram of $\Delta$ show that $\tau$-hypersheaves are $\tau$-sheaves. From the definitions, we immediately see that given a containment of collections of morphisms $S\subseteq S'$, then $S'$-hypersheaves are $S$-hypersheaves. This all yields the following diagram of implications, which will often be used implicitly:
\[\begin{tikzcd}
{\mbox{fpqc hypersheaf}}\ar[r, Rightarrow]\ar[d, Rightarrow]	&	{\mbox{fpqc sheaf}}\ar[d, Rightarrow]	\\
{\mbox{\'{e}tale hypersheaf}}\ar[r, Rightarrow]				&	{\mbox{\'{e}tale sheaf}.}
\end{tikzcd}\]

\subsection*{Topological rings and formal stacks}
The study of deformation theory comes hand-in-hand with the study of rings with a topology and the associated algebraic geometry. We will follow the definition of an adic $\E_\infty$-ring from \cite[Df.0.0.11]{ec2}, except we will only consider the \emph{connective} case; recall an $\E_\infty$-ring $R$ is \emph{connective} if $\pi_d R=0$ for $d< 0$.

\begin{mydef}\label{adicringsdefiniton}
An \emph{adic ring} $A$ is a classical ring with a topology defined by an $I$-adic topology for some finitely generated ideal of definition $I\subseteq A$. Morphisms between adic rings are continuous ring homomorphisms, defining a $1$-category $\CAlg^\heartsuit_\ad$ of classical adic rings. This comes equipped with a functor to $\CAlg^\heartsuit$, the category of classical rings, which forgets the ideal of definition. An \emph{adic $\E_\infty$-ring} is a connective $\E_\infty$-ring $A$ with an adic structure on $\pi_0 A$. We define the \emph{$\infty$-category of adic $\E_\infty$-rings} as the fibre product
\[\CAlg^\cn_\ad= \CAlg^\cn\underset{\CAlg^\heartsuit}{\times} \CAlg^\heartsuit_\ad\]
so a morphism of adic $\E_\infty$-rings is a map of $\E_\infty$-rings which is continuous on $\pi_0$ for the given adic topologies. An adic $\E_\infty$-ring $A$ is said to be \emph{complete} if it is complete with respect to an ideal of definition $I$; see \cite[Df.7.3.1.1 \& Th.7.3.4.1]{sag}. An $\E_\infty$-ring $R$ is \emph{local} if $\pi_0 R$ is a local ring, and we call an adic $\E_\infty$-ring $R$ local if the topology on $\pi_0 R$ is defined by the maximal ideal of $\pi_0 R$. We give $\CAlg^\heartsuit_\ad$ and $\CAlg^\cn_\ad$ the usual Grothendieck topologies (fpqc, \'{e}tale, etc.) via the forgetful functors to $\CAlg^\heartsuit$ and $\CAlg^\cn$, respectively. There is also a functor $\CAlg^\cn \to \CAlg^\cn_\ad$ which equips a connective $\E_\infty$-ring with the adic topology associated with the ideal $(0)$.
\end{mydef}

The geometric definition of a formal spectral Deligne--Mumford stack follows.

\begin{mydef}\label{formalspdm}
Following \cite[Pr.8.1.2.1]{sag}, let $\Spf\colon\CAlg^\cn_\ad\to \infty\Top_{\CAlg}^\loc$ be the functor, from the $\infty$-category of adic $\E_\infty$-rings to the $\infty$-category of locally spectrally ringed $\infty$-topoi (\cite[Df.1.4.2.1]{sag}), defined as the left Kan extension of $\Spec\colon \CAlg^\cn \to \infty\Top_{\CAlg}^\loc$ along the inclusion $\CAlg^\cn \to \CAlg^\cn_\ad$. Roughly speaking, $\Spf A$ is the colimit of all $\Spec B$ where $B$ is an $A$-algebra such that $\Spec B\to \Spec A$ factors through the closed subscheme $K$ of $\Spec A$ associated with an ideal of definition $I$. A \emph{formal spectral Deligne--Mumford stack} is a strictly Henselian (\cite[Df.1.4.2.1]{sag}) spectrally ringed $\infty$-topos with a cover by affine formal spectral Deligne--Mumford stacks; see \cite[Df.8.1.3.1]{sag}. Let $\fSpDM$ denote the full $\infty$-subcategory of $\infty\Top_\CAlg^\sHen$ spanned by formal spectral Deligne--Mumford stacks. Similarly, one can define a $2$-category $\fDM$ of classical formal Deligne--Mumford stacks (\Cref{classicalformalspectra}) where we further assume all such objects are \emph{locally Noetherian}; see \Cref{locallynoetheridefinito} below.
\end{mydef}

\begin{remark}[Invariance of $\Spf$ under completion]\label{completioninimplicitlol}
The formal spectrum $\Spf A$ is equivalent to $\Spf A^\wedge_I$ where $I$ is any finitely generated ideal of definition for $A$; see \cite[Rmk.8.1.2.4]{sag}. In particular, when we write $\Spf A$, we will implicitly complete $A$ with respect to a choice of finitely generated ideal of definition $I$.
\end{remark}

\begin{mydef}\label{locallynoetheridefinito}
Let $\x=(\X, \O_\x)$ be a formal spectral Deligne-Mumford stack, so $\X$ is an $\infty$-topos and $\O_\x$ is a sheaf of $\E_\infty$-rings on $\X$. We call an object $U$ inside $\X$ \emph{affine} if the locally spectrally ringed $\infty$-topos $(\X_{/U}, \O_\x|_U)$ is equivalent to $\Spf A$ for some adic $\E_\infty$-ring $A$. We will also say that $\x$ is \emph{locally Noetherian} if for every affine object $U$ of $\X$, the $\E_\infty$-ring $\O_\x(U)$ is Noetherian in the sense of \cite[Df.7.2.4.30]{ha}.
\end{mydef}

Note that $\Spf B$ is locally Noetherian if and only if $B$ itself is a Noetherian $\E_\infty$-ring; see \cite[Pr.8.4.2.2]{sag}.

\subsection*{Functor of points}
The classical moduli stack of $p$-divisible groups is neither a Deligne--Mumford nor an Artin stack. This necessitates our use of a functorial point of view, for both classical and spectral (formal) algebraic geometry.

\begin{notation}\label{presheavesandaffinestuff}
Denote the opposite category $\CAlg^\op$ as $\Aff$ along with all of the super/subscripts such as $(-)^\cn$, $(-)_\ad$, and $(-)^\heartsuit$, as they apply to $\CAlg$.
\end{notation}

Our justification for working in presheaf $\infty$-categories such as $\Fun(\CAlg^\cn,\Spc)=\P(\Aff^\cn)$ is justified by the following commutative diagram of fully faithful functors:
\begin{equation}\label{diagramoffullyfaithfulstuff}\begin{tikzcd}
{\Aff^\heartsuit_\lN}\ar[dd, "{\mathrm{(a)}}"]\ar[rr]\ar[rd]	&&	{\Aff^\cn}\ar[rd]\ar[dd, "{\mathrm{(c)}}" {yshift=5pt}]	&&	\\
	&	{\Aff^\heartsuit_{\ad, \lN}}\ar[rr, crossing over]	&&	{\Aff^\cn_\ad}\ar[dd]	&	\\
{\DM_\lN}\ar[rr]\ar[rd, "{\mathrm{(b)}}"]	&&	{\SpDM}\ar[rd, "{\mathrm{(d)}}"]	&&						\\
	&	{\fDM}\ar[rr]\ar[from=uu, crossing over]	&&	{\fSpDM}\ar[r]	&	{\P(\Aff^\cn).}
\end{tikzcd}\end{equation}
The $\lN$ subscript denotes those full $\infty$-subcategories spanned by locally Noetherian objects, see \Cref{locallynoetheridefinito}, which is required as $\fDM$ is defined with an implicit locally Noetherian hypothesis. The upper face is defined using \Cref{presheavesandaffinestuff}. The lower face consists of the (nerve of the) $2$-category $\DM_{\lN}$ of \emph{locally Noetherian Deligne--Mumford stacks} (\cite[Df.1.2.4.1 \& Cor.1.2.4.8]{sag}), the $\infty$-category $\SpDM$ of \emph{spectral Deligne--Mumford stacks} (\cite[Df.1.4.4.2]{sag}), the $2$-category $\fDM$ of \emph{formal Deligne--Mumford stacks} (\Cref{classicalformalspectra}), and finally the $\infty$-category of \emph{formal spectral Deligne--Mumford stacks} $\fSpDM$ of \Cref{formalspdm}. The definitions and fully faithfulness of the functors above are explained in \Cref{fffffstacks}, except the functors (a)-(d), which can be justified as follows:
\begin{itemize}
\item[(a)]	is fully faithful as this holds without the locally Noetherian hypotheses; see \cite[Rmk.1.2.3.6]{sag} and restrict to the underlying $2$-category.
\item[(b)]	is fully faithful by using part (d) below and \Cref{spectralanddiscrete}. Indeed, if $G\circ F$ and $G$ are fully faithful, then so is $F$.
\item[(c)]	is fully faithful by making a connective version of \cite[Rmk.1.4.7.1]{sag}; this is justified by \cite[Cor.1.4.5.3]{sag}.
\item[(d)]	is fully faithful as both $\SpDM$ and $\fSpDM$ being defined as full $\infty$-subcategories of $\infty\Top_\CAlg^\loc$ and the fact that spectral Deligne--Mumford stacks are examples of formal spectral Deligne--Mumford stacks by \cite[p. 628]{sag}.
\end{itemize}

Similarly, we will consider most of classical algebraic geometry as living in the $2$-category $\Fun(\CAlg^\heartsuit, \Spc_{\leq 1})$ which we then embed inside the $\infty$-category $\P(\Aff^\heartsuit)$ using the limit preserving inclusion $\Spc_{\leq 1}\to \Spc$. Quasi-coherent sheaves of functors in $\P(\Aff^\cn)$ are defined as in \cite[Df.6.2.2.1]{sag}.

\begin{warn}[Quasi-coherent sheaves on formal spectral Deligne--Mumford stacks]\label{warningabovealmostconnectivestuff}
When we consider quasi-coherent sheaves on a formal spectral Deligne--Mumford stack $\x$, then what we write as $\QCoh(\x)$ is what Lurie would write as $\QCoh(h_\x)$. In other words, we consider the $\infty$-categories of quasi-coherent sheaves of formal spectral Deligne--Mumford stacks through their functors of points. By \cite[Cor.8.3.4.6]{sag}, we see that these two notations are equivalent as long as one restricts to \emph{almost connective} quasi-coherent sheaves on both sides. Recall that a quasi-coherent sheaf is \emph{almost connective} if it is bounded below, so if some finite suspension is connective; see \cite[Var.8.2.5.7 \& Rmk.8.2.5.9]{sag}. All of our quasi-coherent sheaves of interest will be cotangent complexes, which are almost connective by definition (\cite[Df.17.2.4.2]{sag}).
\end{warn}

\subsection*{Cotangent complexes}
Given a natural transformation $X\to Y$ between functors in $\P(\Aff^\cn)$ which admits a cotangent complex (\cite[Df.17.2.4.2]{sag}), we write this cotangent complex as $L_{X/Y}$ and consider it as an object of $\QCoh(X)$. A few specific cases can be made more explicit---thank you to an anonymous referee for vastly simplifying example 3 for us.

\paragraph{(1)}	If $\xx\to \yy$ is a morphism of spectral Deligne--Mumford stacks and $X\to Y$ is the associated transformation of functors in $\P(\Aff^\cn)$, then $L_{X/Y}$ is equivalent to $L_{\xx/\yy}$ under the equivalence of categories $\QCoh(X)\simeq \QCoh(\xx)$ by \cite[Cor.17.2.5.4]{sag}. If $\xx=\Spec B$ and $\yy=\Spec A$, then we have further identifications of $L_{X/Y}$ with $L_{B/A}$ under the equivalence of $\infty$-categories $\QCoh(\Spec A)\simeq \Mod_A$; see \cite[Lm.17.1.2.5]{sag}.

\paragraph{(2)}	 If $\x$ is a formal spectral Deligne--Mumford stack, and $X$ is the associated functor in $\P(\Aff^\cn)$, then $L_{X}$ is equivalent to $L^\wedge_{\x}$, the \emph{completed cotangent complex} of \cite[Df.17.1.2.8]{sag}, under the equivalence of categories $\Theta_\x\colon\QCoh(\x)^\acn\xrightarrow{\simeq} \QCoh(X)^\acn$ of \cite[Cor.8.3.4.6]{sag}, where the superscript $\acn$ indicates full $\infty$-subcategories of almost connective objects, meaning bounded below objects. If $\x=\Spf A$ for an adic $\E_\infty$-ring $A$, then $L_{\Spf A}\simeq (L_A)^\wedge_I$ (under the equivalence of $\infty$-categories $\QCoh(\Spf A)\simeq \Mod_A^\Cpl$, where $I$ is a finitely generated ideal of definition for the topology on $\pi_0 A$; see \cite[Ex.17.1.2.9]{sag}.

\paragraph{(3)}	 If $f\colon \x\to \y$ is a morphism of formal Deligne--Mumford stacks and $F\colon X\to Y$ is the associated morphism of functors in $\P(\Aff^\cn)$, then the cofibre $L_{\x/\y}^\wedge$ of the natural map $f^\star L_\y\to L_\x$ is naturally equivalent to $L_{X/Y}$ under the equivalence of categories $\Theta_\x\colon\QCoh(\x)^\acn \xrightarrow{\simeq} \QCoh(X)^\acn$; see \cite[Df.17.1.2.8]{sag} for a definition of $L_{\x/\y}$. Indeed, the naturality of $\Theta_\x$ in $\x$ (\cite[Con.8.3.4.1]{sag}) yields an equivalence $\Theta_\x\circ f^\star\simeq F^\ast\circ \Theta_\y$ of functors. The desired identification then follows from the existence of the fibre sequences
\[f^\star L_\y\to L_\x\to L_{\x/\y},\qquad\qquad F^\ast L_Y\to L_X\to L_{X/Y},\]
the absolute case above (2), and the fact that $\QCoh(\x)^\acn$ and $\QCoh(X)^\acn$ are closed under fibre sequences; see \cite[Cor.8.2.4.13 \& Pr.6.2.3.4]{sag}, respectively.\\

Due to the equivalences above, we will drop the completion symbol from our notation for the cotangent complex between formal spectral Deligne--Mumford stacks. The following standard properties of the cotangent complex of functors will be used without explicit reference:

\begin{itemize}
\item By \cite[Pr.7.4.3.9]{ha}, given $A\to B$ in $\CAlg^\cn$, we have a natural equivalence in $\Mod_{\pi_0 B}$
\[\pi_0 L_{B/A}\simeq \Omega^1_{\pi_0 B/\pi_0 A}.\]
\item By \cite[Pr.17.2.5.2]{sag}, for composable transformations $X\to Y\to Z$ in $\P(\Aff^\cn)$, where each functor (or each transformation) has a cotangent complex, we obtain a canonical fibre sequence in $\QCoh(X)$
\[\left.L_{Y/Z}\right|_X\to L_{X/Z}\to L_{X/Y}.\]
\item Given transformations $X\to Y\gets Y'$ where $L_{X/Y}$ exists, then $L_{X\times_Y Y'/Y'}$ exists and is naturally equivalent to $\pi_1^\ast L_{X/Y}$ where $\pi_1\colon X\times_Y Y'\to X$; see \cite[Rmk.17.2.4.6]{sag}.
\end{itemize}

\begin{warn}[Topological vs algebraic cotangent complexes]
The cotangent complexes in this article are \textbf{not} the same as those developed by Andr\'{e} and Quillen; see \cite{andrehomology}, \cite{quillencotangentcomplex}, and \cite[\href{https://stacks.math.columbia.edu/tag/08P5}{08P5}]{stacks}. In particular, for an ordinary commutative ring $R$ considered as a discrete $\E_\infty$-ring, then $L_R$ is what some call the \emph{topological cotangent complex} of Basterra; see \cite{taqofbasterra} for the original definitions and \cite[\textsection25.3]{sag} for a quantification of this discrepancy.
\end{warn}

\subsection*{Deformation theory}
We use ideas from classical deformation theory as well as Lurie's spectral deformation theory, so we take a moment here to clarify our definitions.

\begin{mydef}\label{deformationtheory}
Let $\G_0$ be a $p$-divisible group over a commutative ring $R_0$ and write $\CAlg^\Cpl_\ad$ for the $\infty$-subcategory of $\CAlg^\cn_\ad$ spanned by complete connective adic $\E_\infty$-rings. Define a functor $\Def_{\G_0}\colon \CAlg^\Cpl_\ad\to \Spc$ by the formula
\[\Def_{\G_0}(A)=\underset{I}{\colim} \left( \BT(A)\underset{\BT(\pi_0 A/I)}{\times}\Hom_{\CAlg^\heartsuit}(R_0,\pi_0 A/I)\right)\]
where the colimit is indexed over all finitely generated ideals of definition $I$ for $\pi_0 A$. A priori this construction yields an $\infty$-category, but \cite[Lm.3.1.10]{ec2} states this is an $\infty$-groupoid. Let $(A,\G, I, \al\colon R_0 \to \pi_0 A/I, \be\colon \G|_{\pi_0 A/I} \simeq \al^\ast \G_0)$ be a deformation of $\G_0$, so an object of $\Def_{\G_0}(A)$. We say $\G$ is the \emph{universal spectral deformation} of $\G_0$ with \emph{spectral deformation ring} $A$ if for every $B$ in $\CAlg^\Cpl_\ad$, the natural map
\begin{equation}\label{universaldefdef}\Map_{\CAlg^\Cpl_\ad}(A,B)\xrightarrow{\simeq} \Def_{\G_0}(B)\end{equation}
is an equivalence. If $A$ is discrete, we say $\G$ is the \emph{universal classical deformation} of $\G_0$ with \emph{classical deformation ring} $A$ if for every discrete $B$ in $\CAlg^\Cpl_\ad$, the natural map (\ref{universaldefdef}) is an equivalence. If such universal spectral (or classical) deformations $(R,\G)$ exist, they are uniquely determined by the pair $(R_0,\G_0)$.
\end{mydef}

The above definition agrees with that in \cite[Df.3.1.11]{ec2} in the cases that the $A$ above is connective. Indeed, in this case, if $B$ is a nonconnective complete adic $\E_\infty$-ring, as defined by \cite[Df.0.0.11]{ec2}, the fact connective cover is a right adjoint and $\BT(B)=\BT(\tau_{\geq 0}B)$ by definition, we obtain the following:
\[\Map_{\CAlg^\ad}(R,B)\simeq \Map_{\CAlg^\ad}(R,\tau_{\geq 0} B)\simeq \Def_{\G_0}(\tau_{\geq 0} B)\simeq \Def_{\G_0}(B).\]

\begin{remark}\label{classicalcomesfromnonclassical}
If a spectral deformation ring $R$ exists for a pair $(R_0,\G_0)$, then a classical deformation ring also does, and it can be taken to be $\pi_0 R$. Indeed, if $B$ is a discrete object of $\CAlg^\Cpl_\ad$ as in \Cref{deformationtheory}, then the fact the truncation functor is a left adjoint on connective objects yields the equivalences
\[\Def_{\G_0}(B)\simeq\Map_{\CAlg^\Cpl_\ad}(R,B)\simeq \Map_{\CAlg^\Cpl_\ad}(\pi_0 R,B)\]
showing that $\pi_0 R$ is the classical deformation ring of $(R_0,\G_0)$.
\end{remark}

\subsection*{Complex-oriented cohomology theories}

Finally, we will need some notions from algebraic topology. 

\begin{mydef}\label{complexperiodic}
Fix an $\E_\infty$-ring $A$. We say that $A$ is:
\begin{enumerate}
	\item	\emph{complex-orientable} if the unit map $e\colon\Sph\to A$ admits a factorisation $\overline{e}$:
\[\Sph\simeq \Sigma^{-2}S^2\simeq \Sigma^{-2}\complexproj^1\to \Sigma^{-2}\complexproj^\infty \xrightarrow{\overline{e}} A;\]
see \cite[\textsection II]{bluebook} or \cite[\textsection4.1.1]{ec2}.
	\item	\emph{weakly 2-periodic} if $\Sigma^2 A$ is a locally free $A$-module of rank $1$, or equivalently, that $\pi_2 A$ is a locally free $\pi_0 A$-module of rank $1$ and the natural map $\pi_2 A\otimes_{\pi_0 A}\pi_{-2} A\to \pi_0 A$ is an equivalence.
	\item	\emph{complex periodic} if $A$ is \emph{complex orientable} and \emph{weakly 2-periodic}.
\end{enumerate} 
\end{mydef}

Suppose $A$ is an even $2$-periodic $\E_\infty$-ring, meaning that $\pi_{2n+1} A=0$ for all $n\in \Z$ and there is a unit of $\pi_\ast A$ of degree $2$. Then $A$ a classical exercise shows that $A$ is complex-oriented and the unit in $\pi_2 A$ acts as a trivialisation of the $\pi_0 A$-module $\pi_2 A$. In other words, even $2$-periodic $\E_\infty$-rings are complex periodic.\\

Associated with each complex-orientable $\E_\infty$-ring $A$ is a formal group $\Spf A^0(\CP^\infty)$; see \cite[\textsection II.2]{bluebook}. In fact, we can say a little more in spectral algebraic geometry.

\begin{mydef}\label{quillenfg}
By \cite[Con.4.1.13]{ec2}, a complex periodic $\E_\infty$-ring $A$ comes with an associated \emph{Quillen formal group} $\widehat{\G}^\QQ_A$ over $A$. The \emph{classical Quillen formal group} $\widehat{\G}^{\QQ_0}_A $ is the image of $\widehat{\G}^\QQ_A$ under the functor $\FGroup(A)\to \FGroup(\pi_0 A)$, or equivalently as the formal spectrum $\Spf A^0 \complexproj^\infty$. 
\end{mydef}

The definition of the classical Quillen formal group above does not depend on a particular complex orientation of $A$. If we were to fix a choice of complex orientation as well as a choice of unit in $\pi_2 A$, then we would obtain a coordinate for our classical Quillen formal group $\widehat{\G}^{\QQ_0}_A $ and hence produce a \emph{formal group law}; see \cite[\textsection2]{goerssquasicoherent}. This is summarised by the canonical bijection of sets
\[\{\text{Coordinates on } \widehat{\G}^{\QQ_0}_A\} \simeq \{\text{Complex orientations of }A\} \times \{\text{Units in }\pi_2 A\}\]
from \cite[Ex.5.3.7]{ec2}. 

\begin{mydef}
	A formal group $\widehat{\G}$ over a ring $R$ is \emph{Landweber exact} if the defining map from $\Spec R$ to the moduli stack of formal groups is flat. A complex periodic $\E_\infty$-ring is \emph{Landweber exact} if its associated classical Quillen formal group is.
\end{mydef}

There is also a fundamental line bundle associated with a formal group.

\begin{mydef}
	Recall from \cite[\textsection4.2.5]{ec2}, the \emph{dualising line} of a formal group $\widehat{\G}$ over a commutative ring $R$ is the $R$-linear dual of its Lie algebra $\Lie(\widehat{\G})$. This \emph{Lie algebra} of a formal group is the tangent space of $\widehat{\G}$ over $R$ at the unit section $\O_{\widehat{\G}}\to R$; see \cite{zinkcartiertheorie} for a discussion about Lie algebras associated with formal groups.
\end{mydef}



\section{The statement of Lurie's theorem}\label{lurietheoremsection}
First, let us fix some universal notation for this article. 

\begin{notation}[Fixed adic $\E_\infty$-ring $A$]\label{fixeda}
Let $A$ denote some fixed complete local Noetherian adic $\E_\infty$-ring with perfect residue field of characteristic $p$. Write $A_0$ for $\pi_0 A$, $\m_A$ for the maximal ideal of $A_0$, and $\k_A$ for the residue field. Let us further assume that $A$ is flat over the sphere $\Sph$, simply to avoid the ambiguous notation if $A$ is an Eilenberg--Mac\ Lane spectrum. Everything in this article, except \Cref{fetsubsection}, is now implicitly over $\Spf A$, so we have applied $-\times \Spf A$ implicitly throughout. The reader might keep in mind the initial case of the $p$-complete sphere $A=\Sph_p$ with associated $A_0$ the $p$-adic integers $\Z_p$. This will be the usual choice. Others include the spherical Witt vectors $\Sph W(\kappa)$ of a perfect field $\kappa$ of characteristic $p$; see \cite[\textsection5.1]{ec2}.
\end{notation}


First, let us recall the definition of a $p$-divisible group over an $\E_\infty$-ring; see \cite[Df.2.0.2]{ec2} for this definition, and \cite[\textsection6]{ec1} or \cite[\textsection2]{ec3} for wider discussions. 

\begin{mydef}[{\cite[Df.6.5.1]{ec1} \& \cite[Df.2.0.2]{ec2}}]\label{definitionofourcategories}
Let $R$ be a connective $\E_\infty$-ring. A \emph{$p$-divisible (Barsotti--Tate) group over a connective $\E_\infty$-ring $R$} is a functor $\G\colon \CAlg^\cn_R\to \Mod_\Z^\cn$ with the following properties:
\begin{enumerate}
\item For every connective $\E_\infty$-$R$-algebra $B$, the $\Z$-module $\G(B)[1/p]$ vanishes.
\item For every finite abelian $p$-group $M$, the functor
\[\CAlg_R^\cn\to \Spc\qquad B\mapsto \Map_{\Mod_\Z}(M, \G(B))\]
is corepresented by a finite flat $\E_\infty$-$R$-algebra.
\item The map $p\colon \G\to \G$ is locally surjective with respect to the finite flat topology.
\end{enumerate}
A $p$-divisible group over a general $\E_\infty$-ring $R$ is a $p$-divisible group over its connective cover. The $\infty$-category $\BT(R)$ of $p$-divisible groups over an $\E_\infty$-ring $R$ is the full $\infty$-subcategory of $\Fun(\CAlg_{\tau_{\geq 0}R}^\cn, \Mod_\Z^\cn)$ spanned by $p$-divisible groups.
\end{mydef}

This leads to the definition of the main object of interest in this article: the spectral moduli stack $\M_\BT$ and its substacks $\M_\BTn$.

\begin{mydef}[{\cite[Df.3.2.1]{ec2}}]\label{moduliofpdivdef}
 Let $\M_\BT$ be the \emph{moduli stack of $p$-divisible groups}, which is the functor inside $\P(\Aff^\cn)$ sending $R$ to the $\infty$-groupoid core $\BT(R)^\simeq$; see \cite[Df.3.2.1]{ec2}. We say a $p$-divisible group $\G$ has \emph{height $n$} if the $\E_\infty$-$R$-algebra corepresenting the functor
\[\CAlg_R^\cn\to \Spc\qquad B\mapsto \Map_{\Mod_\Z}(\Z/p\Z, \G(B))\]
is finite of rank $p^n$. This notion of height is invariant under base change, as if $A$ is a finite flat $R$-algebra rank $r$, then for any map of $\E_\infty$-rings $R\to R'$, a degenerate Tor-spectral sequence shows that $A\otimes_R R'$ is a finite flat $R'$-algebra of rank $r$. From this, we further define the open subfunctor $\M_{\BTn}$ for all $n\geq 1$ consisting of all $p$-divisible groups of \emph{height $n$}.
\end{mydef}

\begin{mydef}
Let $\x$ be a formal spectral Deligne--Mumford stack. A \emph{$p$-divisible group over $\x$} is a natural transformation $\G\colon\x\to \M_\BT$ in $\P(\Aff^\cn)$. By \cite[Pr.3.2.2(4)]{ec2}, this is equivalent to a coherent family of $p$-divisible groups $\G_{B_i}$ on $\Spec (B_i)^\wedge_{J_i}$, where the collection $\{\Spf B_i\to \x\}_i$ form an affine \'{e}tale cover of $\x$ and $J_i$ is an ideal of definition for $B_i$. We say $\G$ has \emph{height $n$} if this map factors through $\M_\BTn$. 
\end{mydef}

Similarly, one can define formal groups over $\E_\infty$-rings and study their associated moduli stacks.

\begin{mydef}[{\cite[Df.1.6.1]{ec2}}]\label{formalgroupsdefintion}
	Recall that an \emph{abelian group object} in an $\infty$-category $\C$ with finite limits is a functor from the opposite category of free abelian groups of finite rank to $\C$ which commutes with finite products. Write the $\infty$-category of such objects as $\Ab(\C)$. A \emph{formal group} over an $\E_\infty$-ring $A$ is an object of $\Ab(\fSpDM_{/\Spec \tau_{\geq 0}A})$ which is isomorphic in $\fSpDM_{/\Spec \tau_{\geq 0}A}$ to $\Spf B$, where $B$ is an $\tau_{\geq 0}A$-algebra in $\CAlg^\cn_\ad$ such that there is an isomorphism of graded $\pi_\ast A$-modules
	\[\pi_\ast B \simeq \prod_{n\geq 0} (\mathrm{Sym}^n_{\pi_0 A}(M)\otimes_{\pi_0 A} \pi_\ast A)\]
	for some finite finite rank projective $\pi_0 A$-module $M$ which generates an ideal of definition in $\pi_0 A$; this is a combination of \cite[Pr.1.4.11 \& Df.1.6.1]{ec2}.
\end{mydef}

The following construction connects $p$-divisible groups and formal groups.

\begin{remark}[Identity component]\label{connectedcomponentremark}
	Recall from \cite[Th.2.0.8]{ec2}, for each $p$-divisible group $\G$ over a $p$-complete $\E_\infty$-ring $R$ there is a unique formal group $\G^\circ$ over $R$ such that on connective $\E_\infty$-$\tau_{\geq 0}R$-algebras $A$ which are truncated and $p$-nilpotent we can describe $\G^\circ(A)$ as the fibre of $\G(A)\to \G(A^\red)$ induced by the quotient by the nilradical; see \cite[(2.2)]{tatepdiv} for a classical reference.
\end{remark}

Our main object of study is the spectral moduli stack $\M_\BTn$, although we are also interested in its relationship to the underlying classical moduli stack.

\begin{notation}\label{heartsuitnotation}
For a functor $\M\colon \CAlg^\cn\to \Spc$, write $\M^\heartsuit$ for its restriction along the inclusion $\CAlg^\heartsuit\to \CAlg^\cn$. This commutes with finite products
\[(X\times Y)^\heartsuit\xrightarrow{\simeq} X^\heartsuit\times Y^\heartsuit.\]
\end{notation}

Let us also recall what almost perfection means in spectral algebraic geometry.

\begin{mydef}\label{almostperfectdefin}
	Recall that a quasi-coherent sheaf $\FF$ on a functor $X\colon \CAlg^\cn\to \Spc$ is \emph{almost perfect} if for all connective $\E_\infty$-rings $R$ and all morphisms of presheaves $\eta\colon \Spec R\to X$, the $R$-module $\eta^\ast\FF$ is almost perfect; see \cite[Df.7.2.4.10]{ha}. If $R$ is a Noetherian $\E_\infty$-ring (\cite[Df.7.2.4.30]{ha}), then by \cite[Pr.7.2.4.17]{ha}, an $R$-module $M$ is almost perfect if $\pi_k M$ is a finitely presented $\pi_0 R$-modules for every integer $k$ and $\pi_d M=0$ for $d \ll 0$.
\end{mydef}

Finally, we need some classical algebro-geometric adjectives in the spectral world.

\begin{mydef}
	We say that a formal spectral Deligne--Mumford stack $\x$ is \emph{quasi-compact} (qc) if every étale cover $\{\x_\al \to \x\}$ of $\x$ has a finite subcover. We say $\x$ is \emph{quasi-separated} (qs) if the pullback of any étale map $\y\to \x$, where $\y$ is qc, along the diagonal $\x \to \x\times \x$ is itself qc; see \Cref{qcqsdefinition} for the relative formulation.
\end{mydef}

We can now define the sites upon which we will soon define our sheaves of $\E_\infty$-rings. The only property of a morphism of stacks we have not defined yet is \emph{formally étale}, although this is the subject of \Cref{fetsubsection}, and is largely analogous to the classical definition found in \cite[\href{https://stacks.math.columbia.edu/tag/02HFN}{02HFN}]{stacks}.

\begin{mydef}\label{definitionofsites}
Define $\C_{A_0}$ as the full $\infty$-subcategory of $\P(\Aff^\heartsuit)_{/{\M}^\heartsuit_{\BTn}}$ spanned by those objects $\G_0\colon\x_0\to{\M}_{\BTn}^\heartsuit$ where $\x_0$ is a locally Noetherian qcqs formal Deligne--Mumford stack with perfect residue fields at all closed points, the cotangent complex $L_{\x_0/{\M}_{\BTn}}$ is almost perfect inside $\QCoh(\x_0)$, and $\G_0$ is formally \'{e}tale in $\P(\Aff^\heartsuit)$. Note that the cotangent complex $L_{\x_0/{\M}_\BTn}$ exists as the absolute cotangent complexes for $\x_0$ and ${\M}_{\BTn}$ both exist; a consequence of \cite[Pr.17.2.5.1]{sag} and \cite[Pr.3.2.2]{ec2}, respectively. Similarly, define $\C_A$ as the full $\infty$-subcategory of $\P(\Aff^\cn)_{/{\M}_{\BTn}}$ spanned by those objects $\G\colon\x\to{\M}_{\BTn}$ where $\x$ is a locally Noetherian qcqs formal spectral Deligne--Mumford stack with perfect residue fields at all closed points, and $\G$ is formally \'{e}tale in $\P(\Aff^\cn)$. We will endow $\C_{A_0}$ and $\C_{A}$ with both the fpqc and \'{e}tale topologies through the forgetful map to $\P(\Aff^\heartsuit)$ and $\P(\Aff^\cn)$, respectively.
\end{mydef}

The precise version of Lurie's theorem can now be stated.

\begin{theorem}[Lurie's theorem]\label{maintheorem}
Using the notation of \Cref{definitionofsites}, there is an \'{e}tale hypersheaf of $\E_\infty$-rings $\O^\top_\BTn$ on $\C_{A_0}$ such that for a formal affine $\G_0\colon\Spf B_0\to {\M}_\BTn^\heartsuit$ in $\C_{A_0}$ the $\E_\infty$-ring $\O^\top_\BTn(\G_0)=\EE$ has the following properties:
\begin{enumerate}
\item $\EE$ is {complex periodic} and {Landweber exact}.
\item There is a natural equivalence of rings $\pi_0 \EE\simeq B_0$ and $\EE$ is complete with respect to an ideal of definition for $B_0$. In particular, $\EE$ is $\m_A$-complete, hence also $p$-complete.
\item The groups $\pi_k \EE$ vanish for all odd integers $k$. Otherwise, there are natural equivalences of $B_0$-modules $\pi_{2k}\EE\simeq \omega_{\G_0}^{\otimes k}$.
\item There is a natural equivalence of formal groups $\G_0^\circ\simeq \widehat{\G}^{\QQ_0}_\EE$ over $B_0$.
\end{enumerate}
\end{theorem}

We have included a few more details than in the original statement (\cite[Th.8.1.4]{taf}) by incorporating some work of Behrens--Lawson involving Landweber exactness.\\

Now that we have stated Lurie's theorem, let us explain how some of the adjectives defining $\C_A$ and $\C_{A_0}$ come into play.

\paragraph{(Locally Noetherian)}	Our formal Deligne--Mumford stacks are locally Noetherian (\Cref{locallynoetheridefinito}) as completions in the classical world and spectral world do not agree otherwise; see \cite[Warn.8.1.0.4]{sag}. Moreover, in the world of spectral algebraic geometry such objects are better behaved, as detailed by Lurie in \cite[\textsection8.4]{sag}. For example, it will only be obvious that locally Noetherian formal spectral Deligne--Mumford stacks have truncations; see \Cref{truncatationsoffspdmstacksaregood}.
\paragraph{(Qcqs)}	The adjective qcqs is a mild finiteness condition. When a Deligne--Mumford stack $X$ is qcqs, then it has an étale cover $\Spec A\to X$ (qc) and the fibre product $P=\Spec A\times_X \Spec A$ is also an étale cover $\Spec B\to P$ (qs). Moreover, this can be iterated. In other words, if a scheme $X$ is qcqs, it has an étale hypercover by affine schemes; see \Cref{etalehypercoversandsuch}.  This article could be written again, with the word \emph{separated} replacing the word \emph{quasi-separated} and deleting all occurrences of the prefix \emph{hyper}, although the extra generality of hypersheaves can be useful in practice; see \cite{clausenmathew1}.
\paragraph{(Formal geometry)}	Our study of $\M_\BT$ and $\M_n$ is essentially through deformation theory, hence the appearance of formal algebraic geometry. As stated in \cite[Rmk.3.2.7]{ec2}:
\begin{center}
\emph{``The central idea in the proof of Theorem 3.1.15 (of \cite{ec2}) is (\ldots) to guarantee the representability of $\M_\BT$ in a \emph{formal} neighborhood of any sufficiently nice $R$-valued point.''}
\end{center}
\paragraph{(Closed points have perfect residue fields)}	A crucial step in showing our definition of $\O^\top_\BTn$ satisfies the conditions of \Cref{maintheorem} is to reduce ourselves to closed points of the affine objects of $\C_{A_0}$. This allows us to use the Lubin--Tate theories of \cite[\textsection5]{ec2}; see \Cref{wehavesomeuniversaldeformations}.
\paragraph{(Formally \'{e}tale over ${\M}_\BTn$)}	One inspiration for \Cref{maintheorem} is the Serre--Tate theorem, which posits that ${\M}_{\Ell}^\heartsuit$ is formally \'{e}tale over ${\M}_{\BTtwo}^\heartsuit$ (after base change over $\Spf \Z_p$). Phrased differently, in this situation, the deformation theories of ${\M}_{\Ell}^\heartsuit$ and ${\M}_{\BTtwo}^\heartsuit$ agree. The phrase formally \'{e}tale is simply a convenient way to package our deformation theory; see \Cref{fetsubsection}.
\paragraph{(Cotangent complex conditions in $\C_{A_0}$)}	These conditions are essentially finiteness hypotheses, however, they are necessary to apply a deep existence criterion of Lurie (\Cref{impofrepofderham}).\\

Let us now discuss a simple criterion for checking if an object lies in $\C_{A_0}$.

\begin{mydef}\label{finitepresentationdef}
A morphism $f\colon\x_0\to \Spf A_0$ of classical formal Deligne--Mumford stacks is \emph{locally of finite presentation} if for all \'{e}tale morphisms $\Spf B_0\to \x_0$, the induced morphisms of rings $A_0\to B_0$ are of finite presentation. It suffices to check this on a fixed collection of \'{e}tale morphisms $\Spf B_0\to \x_0$ covering $\x_0$. We say $f$ is of \emph{finite presentation} if $f$ is locally of finite presentation and quasi-compact (\Cref{qcqsdefinition}).
\end{mydef}

\begin{prop}\label{simplercrierion}
Let $\G_0\colon \x_0\to {\M}_{\BTn}^\heartsuit$ be a $p$-divisible group defined on a classical formal Deligne--Mumford stack $\x_0$ of finite presentation over $\Spf A_0$ of \Cref{fixeda} such that the associated map into ${\M}_{\BTn}^\heartsuit$ is formally \'{e}tale. Then $\G_0$ lies in $\C_{A_0}$.
\end{prop}

\begin{proof}
Notice that $\x_0$ is locally Noetherian, qcqs, and has all residue fields corresponding to closed points perfect of characteristic $p$ as the morphism $\x_0\to \Spf A_0$ is of finite presentation. Indeed, for locally Noetherian one can use \cite[\href{https://stacks.math.columbia.edu/tag/00FN}{00FN}]{stacks}, for qcqs one can use \cite[\textsection D]{hortzwedhorn}, and the residue fields are perfect as finite field extensions of perfect fields are perfect by \cite[\href{https://stacks.math.columbia.edu/tag/05DU}{05DU}]{stacks}. It remains to show that the cotangent complex in question $L=L_{\x_0/{\M}_{\BTn}}$ is almost perfect. To see this, we consider the composition in $\P(\Aff^\cn)$
\[\x_0\xrightarrow{\G_0} {\M}_{\BTn}\xrightarrow{\pi_2} \Spf A,\]
recall from \Cref{fixeda} that $\M_\BTn$ is short for $\M_\BTn\times \Spf A$. This induces the following fibre sequence in $\QCoh(\x_0)$:
\[\G_0^\ast L_{{\M}_{\BTn}/\Spf A}\to L_{\x_0/\Spf A}\to L\]
Abbreviating the above to $\G_0^\ast L_1\to L_2\to L$, we first focus on $\G_0^\ast L_1$. As a quasi-coherent sheaf on a formal spectral Deligne--Mumford stack $\x_0$, to see $\G_0^\ast L_1$ is almost perfect, it suffices to see that $\eta^\ast \G_0^\ast L_1$ is almost perfect inside $\QCoh(\xx)$ for every morphism $\eta\colon\xx\to\x_0$ where $\xx$ is a spectral Deligne--Mumford stack; see \cite[Th.8.3.5.2]{sag}. Using the base change equivalence
\[L_1=L_{{\M}_{\BTn}/\Spf A}\simeq \pi_1^\ast L_{\M_\BTn}\]
it suffices to show $L_1'=\eta^\ast \G_0^\ast\pi_1^\ast L_{\M_\BTn}$ is almost perfect. By \cite[Cor.8.3.5.3]{sag}, it suffices to check the affine case of $\xx=\Spec R$, where $R$ is a connective $\E_\infty$-ring. Note that $p$ is nilpotent in $\pi_0 R$ as $\Spec R$ maps into $\Spf A$, and $p\in \m_A$ by assumption; see \Cref{fixeda}. By \cite[Pr.3.2.5]{ec2}, we then see that $L_1'$ is almost perfect in $\Mod_R$. Therefore, $\G^\ast_0 L_1$ is almost perfect. \\

Focusing on $L_2$ now, we consider the composition $\x_0\to \Spf A_0\to \Spf A$ and the induced fibre sequence of quasi-coherent sheaves over $\x_0$:
\begin{equation}\label{nowthecofibresequenceonfdmstacks}\left. L_{\Spf A_0/\Spf A}\right|_{\x_0}\to L_{\x_0/\Spf A}=L_2\to L_{\x_0/\Spf A_0}.\end{equation}
By \Cref{truncationsofnoetherianfspdmaregood}, we see $L_{\Spf A_0/\Spf A}$ is almost perfect in $\QCoh(\Spf A_0)$, and pullbacks preserve almost perfectness (\cite[Cor.8.4.1.6]{sag}), hence the first term of (\ref{nowthecofibresequenceonfdmstacks}) is almost perfect. To see the third term of (\ref{nowthecofibresequenceonfdmstacks}) is almost perfect, we may work locally and replace $\x_0$ with $\Spf B_0$ where $B_0$ is a complete discrete adic ring. In this case, we use the assumption that $A_0\to B_0$ is of finite presentation, which implies $L_{B_0/A_0}$ is almost perfect in $\Mod_{B_0}$; see \cite[Th.7.4.3.18]{ha}. By \cite[Pr.7.3.5.7]{sag}, $L_{B_0/A_0}$ is complete with respect to an ideal of definition $J$ for $B_0$, and it follows the $B_0$-module
\[L_{B_0/A_0}\simeq \left(L_{B_0/A_0}\right)^\wedge_J\simeq L_{\Spf B_0/\Spf A_0}\]
is almost perfect. Therefore $L_2$ is almost perfect, so $L$ itself is almost perfect.
\end{proof}


\section{Formally \'{e}tale natural transformations}\label{fetsubsection}
At the heart of spectral algebraic geometry is deformation theory---Lurie (\cite[p.1385]{sag}) even goes as far as to state the heuristic principle:
\[\{\text{spectral algebraic geometry}\}=\{\text{classical algebraic geometry}\}+\{\text{deformation theory}\}.\]
The adjective \emph{formally \'{e}tale} will help us navigate between the two worlds of classical and spectral algebraic geometry using Lurie's spectral deformation theory. More concretely, given a (nice enough) formally \'{e}tale morphism $\x_0\to {\M}$, where $\x_0$ is a classical formal stack, there is a universal spectral lift of $\x_0$. This process allows us to lift objects in classical algebraic geometry to spectral algebraic geometry without changing the underlying classical object; see \Cref{maintheorembeforetopoi}. Lurie often uses an affine and étale version of such a statement, for example, that the functor
\[\pi_0 \colon \CAlg_A^\et \xrightarrow{\simeq} \CAlg_{\pi_0 A}^\et \]
from étale $\E_\infty$-rings over $A$ to étale commutative rings over $\pi_0 A$, is an equivalence; see \cite[Th.7.5.4.2]{ha}. \Cref{maintheorembeforetopoi} is a generalisation of this to the small formally étale site $\C_A$ over $\M_n$.

\subsection{On presheaves of discrete rings}
Let us first consider formally \'{e}tale maps between presheaves of discrete rings. Recall that a map of rings $\widetilde{R} \to R$ is a \emph{square-zero extension} of $R$ if it is surjective with square-zero kernel.

\begin{mydef}\label{formallyetalemap}
A natural transformation $f\colon X\to Y$ of functors in $\P(\Aff^\heartsuit)$ is called \emph{formally \'{e}tale} if for all square-zero extensions $\widetilde{R} \to R$ the natural commutative diagram of spaces
\[\begin{tikzcd}
{X(\widetilde{R})}\ar[r]\ar[d]	&	{X({R})}\ar[d]	\\
{Y(\widetilde{R})}\ar[r]		&	{Y({R})}
\end{tikzcd}\]
is Cartesian. Moreover, we say that $f$ is \emph{formally unramified} if the fibres of the map
\[X(\widetilde{R})\to X(R)\underset{Y(R)}{\times} Y(\widetilde{R})\]
are either empty or contractible.
\end{mydef}

From this definition, we immediately obtain the following.

\begin{prop}\label{superbaseicproeprtoes}
Formally \'{e}tale morphisms in $\P(\Aff^\heartsuit)$ are closed under composition. If $X\xrightarrow{f} Y\xrightarrow{g} Z$ are composable morphisms in $\P(\Aff^\heartsuit)$ such that $g$ is formally unramified and $gf$ is formally \'{e}tale, then $g$ is formally \'{e}tale. Formally \'{e}tale (resp. unramified) morphisms are closed under base change.
\end{prop}

Let us now relate \Cref{formallyetalemap} to the definitions found in classical algebraic geometry.



\begin{prop}\label{formallyetmeanswhat}
Let $f\colon X\to Y$ be a natural transformation of functors in $\P(\Aff^\heartsuit)$. Then $f$ is formally \'{e}tale if and only if for every ring $R$, every square-zero extension of rings $\widetilde{R}\to R$, and every commutative diagram of the form
\begin{equation}\label{squarezerolift}\begin{tikzcd}
{\Spec R}\ar[r]\ar[d]		&	{X}\ar[d, "f"]	\\
{\Spec \widetilde{R}}\ar[r]	&	{Y}
\end{tikzcd}\end{equation}
the mapping space
\[\Map_{\P(\Aff^\heartsuit)_{\Spec R//Y}}(\Spec \widetilde{R}, X)\]
is contractible; here $\P(\Aff^\heartsuit)_{\Spec R//Y}$ denotes the slice $\infty$-category under $\Spec R$ and over $Y$. 
\end{prop}

In other words, if $X\to Y$ is formally étale, there exists a unique lift $\Spec \widetilde{R}\to X$ for (\ref{squarezerolift}).

\begin{proof}
Given a ring $R$, a square-zero extension $\widetilde{R}\to R$, and a commutative diagram (\ref{squarezerolift}), consider the following diagram of spaces:
\[\begin{tikzcd}
{\Map_{R//Y}(\widetilde{R}, X)}\ar[r]\ar[d]	&	{\Map_{/Y}(\widetilde{R}, X)}\ar[r]\ar[d]	&	{\Map_{/Y}(R, X)}\ar[d]	\\
{\Map_{R/}(\widetilde{R}, X)}\ar[r]\ar[d]	&	{\Map(\widetilde{R}, X)}\ar[r]\ar[d]		&	{\Map(R, X)}\ar[d]	\\
{\Map_{R/}(\widetilde{R}, Y)}\ar[r]		&	{\Map(\widetilde{R}, Y)}\ar[r]			&	{\Map(R, Y)}
\end{tikzcd}\]
By definition, the rows and columns are fibre sequences, the fibres in this diagram have been taken with respect to the maps from (\ref{squarezerolift}), we have suppressed the functor $\Spec$, and we have abbreviated the categories above to express only the over/under categories. In other words, $\Map_{R//Y}(-,-)$ is the mapping space in the slice category $\P(\Aff^\cn)_{\Spec R//Y}$. By the Yoneda lemma, the lower-right square is naturally equivalent to (\ref{squarezerolift}). Hence $f$ is formally étale if and only if the space in the upper-left corner is contractible. This finishes the proof.
\end{proof}

Let us list some instances of formally \'{e}tale morphisms found in algebraic geometry.

\begin{example}[Formally \'{e}tale morphisms of schemes]
In the setting of classical algebraic geometry, we usually take the existence of a unique map $\Spec \widetilde{R}\to X$ (under $\Spec R$ and over $Y$) as the definition of a formally \'{e}tale maps of rings (or schemes); see \Cref{formallyetmeanswhat}. An object in $\P(\Aff^\heartsuit)$ represented by a scheme factors through $\Fun(\CAlg^\heartsuit, \Set)$, as mapping spaces between classical schemes are discrete, and we see \Cref{formallyetmeanswhat} precisely matches \cite[\href{https://stacks.math.columbia.edu/tag/02HG}{02HG}]{stacks}.
\end{example}

\begin{example}[Classical Serre--Tate theorem]\label{classicalserretatetheorem}
For an $\E_\infty$-ring $R$, we write $\AVar_g(R)$ for the $\infty$-category of (strict) spectral abelian varieties over $R$; see \cite[Df.1.5.1]{ec1}. In other words, an object is an abelian group object (see \Cref{formalgroupsdefintion}) in the $\infty$-subcategory of (nonconnective) Deligne--Mumford stacks $\xx$ over $\Spec R$ such that the map $\xx \to \Spec R$ is flat and the map $\tau_{\geq 0} \xx \to \Spec \tau_{\geq 0} R$ is proper, locally almost of finite presentation, geometrically reduced, and geometrically connected. The classical Serre--Tate theorem (see \cite[p.854]{serretatevolII} for the original source, or \cite[Th.7.0.1]{ec1} for statement of the spectral version) states that if $\widetilde{R}\to R$ is a square-zero extension of commutative rings and $p$ is nilpotent within them, then the diagram of 1-groupoids
\begin{equation}\label{serretatediagram}\begin{tikzcd}
{\AVar_g(\widetilde{R})^\simeq}\ar[r]\ar[d, "{[p^\infty]}"]	&	{\AVar_g({R})^\simeq}\ar[d, "{[p^\infty]}"]	\\
{\mathrm{BT}^p_{2g}(\widetilde{R})^\simeq}\ar[r]					&	{\mathrm{BT}^p_{2g}({R})^\simeq}
\end{tikzcd}\end{equation}
is Cartesian, where $\mathrm{BT}^p_{2g}(A)$ is the $\infty$-subcategory of $\BT(A)$ spanned by $p$-divisible groups of height $2g$. This implies the morphism of classical moduli stacks $[p^\infty]\colon \M_{\AVar_g}^\heartsuit\to \M_{\BTtwog}^\heartsuit$ sending an abelian variety $\xx$ to its associated $p$-divisible group $\xx[p^\infty]$ (\cite[\textsection2]{tatepdiv}) is formally \'{e}tale \textbf{after} base change over $\Spf \Z_p$. This base change is crucial, as there only exists a map $\Spec R\to \Spf \Z_p$ is when $p$ is nilpotent inside $R$, as the continuous map of rings $\Z_p\to R$ must send $\{p^i\}_{i\geq 0}$ to a convergent sequence in $R$, where $R$ is equipped with the discrete topology. If we fail to make this base change, then (\ref{serretatediagram}) may not be Cartesian; see \cite[p.32]{mythesis} for a simple counter-example.
\end{example}

\begin{example}[Classical Lubin--Tate theorem]\label{lubintateexample}
Another classical example of a formally \'{e}tale map in $\P(\Aff^\heartsuit)$ comes from Lubin--Tate theory. The original source for this is \cite{lubintate} concerning formal group laws. Further exploration of the deformation theory of $p$-divisible groups has been done by Drinfeld \cite[Pr.4.5]{ellipticmodulesi} and Harris--Taylor \cite[II.1]{harristaylorshimuravar}; a summary aimed at homotopy theorists can be found in \cite[\textsection 7.1]{taf}. We will follow \cite[\textsection3]{ec2} as our intended application is for $p$-divisible groups; see \cite[Ex.3.0.5]{ec2} for a statement of the dictionary between deformations of formal and $p$-divisible groups. Let $\G_0$ be a $p$-divisible groups of height $0<n<\infty$ over a perfect field $\kappa$ of characteristic $p$. Then there exists a universal classical deformation $\G$ of $\G_0$ over the classical deformation ring $R^\LT_{\G_0}$; see \cite[Df.3.1.4]{ec2} or the proof of \Cref{lubintatefetexample}.
\end{example}

This implies that the map into $\M_\BTn^\heartsuit$ defining $\G$ is formally \'{e}tale. In fact, we slightly generalise the Lubin--Tate case above using \cite[\textsection3]{ec2}. To do this, we need the concept of a \emph{nonstationary} $p$-divisible group. The original definition is \cite[Df.3.0.8]{ec2}, but we will only need the alternative definition here for discrete Noetherian $\F_p$-algebras which combines \cite[Rmk.3.4.4 \& Th.3.5.1]{ec2}

\begin{mydef}\label{nonstationary}
	A $p$-divisible group $\G_0$ over a discrete Noetherian $\F_p$-algebra $R_0$ is \emph{nonstationary} if the Frobenius morphism on $R_0$ is finite and that the cotangent complex $L_{\Spec R/\M_\BT}$ induced by the defining morphism of $\G_0$ is $1$-connective. In particular, by \cite[Ex.3.0.10]{ec2}, all $p$-divisible groups over $\F_p$-algebras $R_0$ whose Frobenius is surjective are nonstationary.
\end{mydef}

Recall from \Cref{classicalcomesfromnonclassical} that if a spectral deformation ring $R_{\G_0}$ exists for a $p$-divisible group $\G_0$, then the classical deformation ring is simply $\pi_0 R_{\G_0}$. This will be the case below.

\begin{prop}\label{lubintatefetexample}
Fix $A=\Sph_p$ from \Cref{fixeda}. Let $R_0$ be a discrete $\F_p$-algebra such that $L_R$ is an almost perfect $R$-module (\Cref{almostperfectdefin})\footnote{See \cite[Pr.3.3.7 \& Th.3.5.1]{ec2} for many equivalent conditions to $L_R$ being almost perfect.} and $\G_0$ is a {nonstationary} $p$-divisible group over $R_0$ of height $n$. Then the map $\Spf R_{\G_0}\to {\M}_{\BTn}^\heartsuit$ induced by the universal classical deformation of $\G_0$ is formally \'{e}tale. Conversely, if $\G\colon \Spf R\to {\M}_{\BTn}^\heartsuit$ is formally \'{e}tale for a complete Noetherian discrete ring $R$, then for every maximal ideal $\m\subseteq R$ such that the residue field $R/\m=\kappa$ is perfect, then $\G_{B^\wedge_\m}$ is the universal classical deformation of $\G_\kappa$.
\end{prop}

The following proof is essentially an unwinding of definitions.

\begin{proof}
The existence of the spectral deformation ring follows from \cite[Th.3.4.1]{ec2} and the classical deformation ring is obtained by applying $\pi_0$. Let $R\to R/J$ be the quotient map where $R$ is discrete and $J$ is a square-zero ideal. First, we wish to show the commutative diagram of spaces
\begin{equation}\label{lubintateandformallyetatediagram}
\begin{tikzcd}
{(\Spf R_{\G_0})(R)}\ar[r]\ar[d, "l"]	&	{(\M_\BTn^\heartsuit)(R)}\ar[d, "r"]	\\
{(\Spf R_{\G_0})(R/J)}\ar[r, "{b}"]		&	{(\M_\BTn^\heartsuit)(R/J)}
\end{tikzcd}
\end{equation}
is Cartesian. Following \cite[Df.3.1.4]{ec2}, we define a functor $\Def_{\G_0}(-) \colon \CAlg^\cn_\ad \to \Cat_\infty$ by sending an adic $\E_\infty$-ring $A$ to the $\infty$-category
\[\Def_{\G_0}(A)=\underset{I}{\colim}\left(\BT (A)\underset{\BT(\pi_0 A/I)}{\times} \Hom_{\CAlg^\heartsuit}(R_0, \pi_0 A/I)\right)\]
where the colimit is indexed over the filtered system of finitely generated ideals of definition $I\subseteq \pi_0 A$. By \cite[Lm.3.1.10]{ec2}, if $A$ is complete, then $\Def_{\G_0}(A)$ is a space, which in particular holds if $A$ is a discrete ring equipped with the discrete topology. In fact, as $R_{\G_0}$ is the universal deformation of $\G_0$ one obtains for any such $A$ an equivalence of (discrete) spaces
\begin{equation}\label{fibreprocutincolimit}\Hom_{\CAlg^\heartsuit_\ad}(R_{\G_0}, A)\xrightarrow{\simeq} \Def_{\G_0}(A)=\underset{I\in \Nil_0(A)}{\colim}\left( \BT(A)\underset{\BT(A/I)}{\times} \Hom_{\CAlg^\heartsuit}(R_0, A/I) \right)\end{equation}
where the colimit is taken over all finitely generated nilpotent ideals $I$ inside $A$; see \cite[Th.3.1.15]{ec2}. By assumption, the cotangent complex $L_{R_0}$ is almost perfect in $\Mod_{R_0}$, and \cite[Pr.3.4.3]{ec2} then implies that the natural map
\[\Def_{\G_0}(A)\xrightarrow{\simeq} \underset{I\in \Nil(A)}{\colim} F_{A,I}\]
is an equivalence, where now the colimit is indexed over \emph{all} nilpotent ideals $I\subseteq A$ and $F_{A,I}$ is the fibre product of (\ref{fibreprocutincolimit}). Given a fixed nilpotent ideal $J\subseteq A$, denote by $\Nil_J(A)$ the poset of nilpotent ideals of $A$ which contain $J$. We obtain a natural inclusion functor $\Nil_J(A)\to\Nil(A)$, which is cofinal, as any nilpotent ideal $I$ lies within the nilpotent ideal $I+J$. Hence the natural map
\[\underset{I\in \Nil_J(A)}{\colim} F_{A,I} \xrightarrow{\simeq}\underset{I\in \Nil(A)}{\colim} F_{A,I}\]
is an equivalence. The map $l$ of (\ref{lubintateandformallyetatediagram}) is then equivalent to
\[\underset{I\in \Nil_J(R)}{\colim} F_{R,I} \xrightarrow{l} \underset{I\in \Nil_J(R)}{\colim} F_{R/J,I/J}\]
where we used the fact that ideals in $R/J$ correspond to ideals in $R$ containing $J$. If $(\Spf R_{\G_0})(R/J)$ is empty, then so is $(\Spf R_{\G_0})(R)$ and we are done. Otherwise, choose some $x$ in $(\Spf R_{\G_0})(R/J)$ and consider the fibre of $l$ over $x$. As filtered colimits of spaces commute with finite limits we calculate this fibre as
\[\fib_x(l)\simeq \underset{I\in \Nil_J(R)}{\colim} \fib_{x_I}\left(F_{R,I}\xrightarrow{g}F_{R/J,I/J}\right). \]
To simplify this further, we contemplate the commutative diagram in $\Cat_\infty$
\[\begin{tikzcd}
{\BT(R)\underset{\BT(R/I)}{\times} \Hom(R_0, R/I)}\ar[r, "{g}"]\ar[d]	&	{\BT(R/J)\underset{\BT(R/I)}{\times} \Hom(R_0, R/I)}\ar[r]\ar[d]	&	{\Hom(R_0, R/I)}\ar[d]	\\ 
{\BT(R)}\ar[r, "{f}"]	&	{\BT(R/J)}\ar[r]	&	{\BT(R/I).}
\end{tikzcd}\]
The right square and outer rectangle above are Cartesian by definition, so the left square is also Cartesian. This means the natural map $\fib(g)\to \fib(f)$ is an equivalence in $\Cat_\infty$, hence our fibre of $l$ becomes
\[\fib_x(l)\simeq\underset{I\in \Nil_J(R)}{\colim}\left(\fib_{b(x_I)}(\BT(R)\xrightarrow{f} \BT(R/J))\right)\simeq \fib_{b(x)}\left(\BT(R)\xrightarrow{f} \BT(R/J)\right).\]
This shows the fibre of $f$ lies in the essential image of $\Spc\rightarrow \Cat_\infty$ as $\fib_x(l)$ is. As $r$ is $f^\simeq$ we obtain a natural equivalence $\fib(l)\simeq \fib(r)$. As the fibres of $l$ and $r$ are naturally equivalent, we see that (\ref{lubintateandformallyetatediagram}) is Cartesian, so the composition
\[\Spf R_{\G_0}\to {\M}_{\BTn}^\heartsuit \times \Spf \Z_p \to \M_{\BTn}^\heartsuit\to \M_\BT^\heartsuit\]
is formally \'{e}tale---here we have written the usually implicit $-\times \Spf\Z_p$. To see the first map in the composition above is formally \'{e}tale, we use that the last map is open and hence formally \'{e}tale, the second last map is the base change of the formally unramified map $\Spf \Z_p\to\Spec \Z$, and the cancellation statement from \Cref{superbaseicproeprtoes}.\\

Let us omit a proof of the converse statement; the $\E_\infty$-version is \Cref{wehavesomeuniversaldeformations} and the proofs in both cases are analogous.
\end{proof}


\subsection{On presheaves of connective $\E_\infty$-rings}

We are now in the position to make a spectral definition. To that end, recall the definition of (trivial) square-zero extensions of $\E_\infty$-rings from \cite[Df.7.4.1.6]{ha}. 

\begin{mydef}\label{formallyetalemapsofeinftyrings}
Let $f\colon X\to Y$ be a natural transformation of functors in $\P(\Aff^\cn)$. For an integer $0\leq n\leq \infty$, we say $f$ is \emph{$n$-formally \'{e}tale} if for all square-zero extensions of connective $n$-truncated $\E_\infty$-rings $\widetilde{R}\to R$, meaning $\pi_d \widetilde{R}=\pi_d R=0$ for $d\not\in [0,n]$, the natural diagram of spaces
\[\begin{tikzcd}
{X(\widetilde{R})}\ar[r]\ar[d]	&	{X(R)}\ar[d]	\\
{Y(\widetilde{R})}\ar[r]		&	{Y(R)}
\end{tikzcd}\]
is Cartesian. We abbreviate $\infty$-formally \'{e}tale to \emph{formally \'{e}tale}.
\end{mydef}

It will be proven in \Cref{etalemapsofrings} that étale morphisms between connective $\E_\infty$-rings are in particular formally étale.

\begin{remark}[Spectral to classical]\label{biggerimplieslessthan}
If $f$ is $n$-formally \'{e}tale, then $f$ is also $m$-formally \'{e}tale for all $0\leq m\leq n\leq \infty$. In particular, for any $0\leq n\leq \infty$, if $f$ is $n$-formally \'{e}tale in $\P(\Aff^\cn)$, then $X^\heartsuit\to Y^\heartsuit$ is formally \'{e}tale inside $\P(\Aff^\heartsuit)$.
\end{remark}

A converse statement also holds.

\begin{remark}[Classical to spectral]\label{truncationsgiveetale}
Write $\tau_{\leq 0}^\ast\colon \P(\Aff^\heartsuit)\to \P(\Aff^\cn)$ for the functor induced by the truncation $\CAlg^\cn\to\CAlg^\heartsuit$. We claim that if $X\to Y$ is formally \'{e}tale in $\P(\Aff^\heartsuit)$, then it follows that $\tau_{\leq 0}^\ast X\to \tau_{\leq 0}^\ast Y$ is ($\infty$-) formally \'{e}tale inside $\P(\Aff^\cn)$. Indeed, for each square-zero extension of connective $\E_\infty$-rings $\widetilde{R}\to R$ we want to show the the diagram of spaces
\[\begin{tikzcd}
{X(\pi_0\widetilde{R})}\ar[r]\ar[d]	&	{X(\pi_0 R)}\ar[d]	\\
{Y(\pi_0\widetilde{R})}\ar[r]		&	{Y(\pi_0 R)}
\end{tikzcd}\]
is Cartesian. If the map $\rho\colon\pi_0 \widetilde{R}\to \pi_0 R$ is a square-zero extension of classical rings, then the above square would be Cartesian from our hypotheses. To see that $\rho$ is also a classical square-zero extension, notice that the fibre sequence $M\to \widetilde{R}\to R$ of connective $R$-modules shows that $\rho$ is surjective. Moreover, we see the kernel of $\rho$ is not $\pi_0 M$, but the image of the map $\pi_0 M\to \pi_0\widetilde{R}$. This does not worry us, as the multiplication map $M\otimes_{\widetilde{R}}M\to M$ is nullhomotopic by \cite[Pr.7.4.1.14]{ha}, hence the image of $\pi_0 M$ in $\pi_0\widetilde{R}$ squares to zero, and we see $\rho$ is a square-zero extension of classical rings.
\end{remark}

\begin{remark}[Classical stacks to spectral stacks]\label{classicaltonot}
If $\x\to \y$ is a formally \'{e}tale morphism of (locally Noetherian) classical formal Deligne--Mumford stacks inside $\P(\Aff^\heartsuit)$, then the corresponding morphism inside $\P(\Aff^\cn)$ is 0-formally \'{e}tale. This follows by the fully faithfulness of $\fDM\to \fSpDM$; see \Cref{spectralanddiscrete}.
\end{remark}

\begin{remark}[Formally étale vs étale]\label{differenttoluriedefinition}
Our definition of formally \'{e}tale deviates from Lurie's definition of \'{e}tale morphisms (\cite[Df.7.5.0.4]{ha}) as there is no flatness assumption. However, even in $\P(\Aff^\heartsuit)$ a formally \'{e}tale morphism of discrete rings need not be flat. This means there is no inherent descent theory for formally \'{e}tale morphisms. For more in this direction, see \Cref{followuphehe}.
\end{remark}

The basic properties of \Cref{superbaseicproeprtoes} also hold in $\P(\Aff^\cn)$.

\begin{prop}\label{bcandsuch}
Let $0\leq n\leq \infty$ and $X\xrightarrow{f} Y\xrightarrow{g} Z$ be composable morphisms in $\P(\Aff^\cn)$ where $g$ is $n$-formally \'{e}tale. Then $f$ is $n$-formally \'{e}tale if and only if $h$ is $n$-formally \'{e}tale. Moreover, $n$-formally \'{e}tale morphisms are closed under base change.
\end{prop}

There are a few alternative ways to test if a map $X\to Y$ is formally \'{e}tale in $\P(\Aff^\cn)$. Although Lurie does not directly discuss the adjective formally \'{e}tale in \cite[\textsection17]{sag}, many of the techniques below follow his ideas. Recall the definition of the cotangent complex $L_{X/Y}$ for a morphism $f\colon X\to Y$ between two functors in $\P(\Aff^\cn)$ of \cite[Df.17.2.4.2]{sag} as well as what it means for $f$ to be (infinitesimally) cohesive or nilcomplete; see \cite[Df.17.3.7.1]{sag}.

\begin{prop}\label{alternativefetexpressions}
Let $f\colon X\to Y$ be a map of functors in $\P(\Aff^\cn)$ and $0\leq n\leq \infty$. 
\begin{enumerate}
\item The map $f$ is $n$-formally \'{e}tale for \emph{finite} $n$ if and only if $X\to Y$ is $0$-formally \'{e}tale and for every connective $n$-truncated $\E_\infty$-ring $R$ the natural diagram of spaces
\[\begin{tikzcd}
{X(R)}\ar[r]\ar[d]	&	{X(\pi_0 R)}\ar[d]	\\	
{Y(R)}\ar[r]		&	{Y(\pi_0 R)}
\end{tikzcd}\]
is Cartesian. If $f$ is nilcomplete, then the $n=\infty$-case also holds.
\item If $f$ is infinitesimally cohesive, then $f$ is formally \'{e}tale if and only if for all \emph{trivial} square-zero extensions of connective truncated $\E_\infty$-rings $\widetilde{R}\to R$, so $n$-truncated for some finite $n$, the natural diagram of spaces
\[\begin{tikzcd}
{X(\widetilde{R})}\ar[r]\ar[d]	&	{X(R)}\ar[d]	\\
{Y(\widetilde{R})}\ar[r]		&	{Y(R)}
\end{tikzcd}\]
is Cartesian. 
\item If $f$ is infinitesimally cohesive and admits a connective cotangent complex $L_{X/Y}$, then $f$ is formally \'{e}tale if and only if $L_{X/Y}$ vanishes.
\end{enumerate}
\end{prop}

Thank you to an anonymous referee for correcting a previous version of (2) above.

\begin{remark}[Obstruction groups]
If $X\to Y$ is infinitesimally cohesive, nilcomplete, and $L_{X/Y}$ exists and is connective, then $X\to Y$ is $n$-formally \'{e}tale if certain Ext-groups $\Ext_R^m(\eta^\ast L_{X/Y}, M)$ vanish in a range, for certain discrete objects $(R,\eta,M)$ of $\Mod^X_\cn$, \`{a} la the deformation theory of \cite{illusie}. There is also a sharpening of part 3 above which deals with an $n$-connective cotangent complex $L_{X/Y}$, which is \textbf{not} equivalent to $X\to Y$ being $n$-formally \'{e}tale. These ideas will not be used here though.
\end{remark}

\begin{proof}
We start with part 1. Suppose $f$ is $n$-formally \'{e}tale for a finite $n\geq 0$, then $f$ is 0-formally \'{e}tale by \Cref{biggerimplieslessthan}. Given a connective $n$-truncated $\E_\infty$-ring $R$, then for any $0\leq m\leq n$ we can consider the commutative diagram of spaces
\begin{equation}\label{inductionwithsquares}\begin{tikzcd}
{X(\tau_{\leq m+1} R)}\ar[r]\ar[d]	&	{X(\tau_{\leq m} R)}\ar[r]\ar[d]	&	{X(\pi_0 R)}\ar[d]	\\
{Y(\tau_{\leq m+1} R)}\ar[r]		&	{Y(\tau_{\leq m} R)}\ar[r]		&	{Y(\pi_0 R).}
\end{tikzcd}\end{equation}
Above, the left square is always Cartesian by virtue of $f$ being $n$-formally \'{e}tale as $\tau_{\leq m+1} R\to \tau_{\leq m} R$ is a square-zero extension of $\E_\infty$-rings by \cite[Cor.7.1.4.28]{ha}. To show the outer rectangle Cartesian we use induction. The base case of $m=0$ is tautological. For $m\geq 1$, the right square is Cartesian by our inductive hypotheses, hence the whole rectangle is Cartesian. Conversely, if the second condition of part 1 holds, we consider a square-zero extension of $n$-truncated connective $\E_\infty$-rings $\widetilde{R}\to R$ and the natural commutative diagram of spaces
\[\begin{tikzcd}
X(\widetilde{R}) \arrow[dd] \arrow[rr] \arrow[rd] &                            & X(\pi_0\widetilde{R}) \arrow[dd] \arrow[rd] &                       \\
                                                  & X(R) \arrow[rr, crossing over] &                                             & X(\pi_0 R) \arrow[dd] \\
Y(\widetilde{R}) \arrow[rd] \arrow[rr]            &                            & Y(\pi_0\widetilde{R}) \arrow[rd]            &                       \\
                                                  & Y(R) \arrow[rr]\arrow[from=uu, crossing over]            &                                             & {Y(\pi_0 R).}            
\end{tikzcd}\]
The back and front faces are Cartesian by the second condition of part 1, and the right-most face is Cartesian as the second condition of part 1 also assumes $f$ is $0$-formally \'{e}tale. Hence by a base change argument, we see the left-most square is Cartesian, and we are done. For the $n=\infty$-case, suppose $X\to Y$ is nilcomplete, meaning that for every connective $\E_\infty$-ring $R$, the diagram of spaces
\[\begin{tikzcd}
X(R)\ar[r]\ar[d]	&	{\lim X(\tau_{\leq n} R)}\ar[d]	\\
Y(R)\ar[r]		&	{\lim Y(\tau_{\leq n} R)}
\end{tikzcd}\]
is Cartesian. Combining this diagram with the finite case above yields the desired conclusion.\\

For part 2, if $f$ is formally \'{e}tale, then we immediately see the second condition holds. Conversely, let $e\colon\widetilde{R}\to R$ be a square-zero extension of a connective $\E_\infty$-ring $R$ by a connective $R$-module $M$ and a derivation $d\colon L_R\to \Sigma M$. By definition (\cite[Df.7.4.1.6]{ha}) $\widetilde{R}$ is defined by the Cartesian diagram of connective $\E_\infty$-rings
\[\begin{tikzcd}
{\widetilde{R}}\ar[d]\ar[r, "e"]	&	{R}\ar[d, "\rho"]	\\
{R}\ar[r, "0"]			&	{R\oplus \Sigma M}
\end{tikzcd}\]
where the lower-horizontal map is induced by the zero map $L_R\to \Sigma M$ and the right-vertical map is induced by the derivation $d$. This Cartesian diagram of connective $\E_\infty$-rings then induces the commutative diagram of spaces
\begin{equation}\label{cubeforourplesure}\begin{tikzcd}
{X(\widetilde{R})}\ar[rr, "{X(e)}"]\ar[dd]\ar[rd]	&&	{X(R)}\ar[dd]\ar[rd, "{X(\rho)}"]			&&&	\\
	&	{X({R})}\ar[rr, crossing over]		&&	{X(R\oplus \Sigma M)}\ar[dd]\ar[rr]	&&	{X(R)}\ar[dd]	\\
{Y(\widetilde{R})}\ar[rr]\ar[rd]		&&	{Y(R)}\ar[rd, "{Y(\rho)}"]		&&&	\\
	&	{Y({R})}\ar[rr, "{Y(0)}"]\ar[from=uu, crossing over]	&&	{Y(R\oplus \Sigma M)}\ar[rr]			&& {Y(R).} 
\end{tikzcd}\end{equation}
The left cube is Cartesian from our assumption that $f$ is infinitesimally cohesive. By assumption the right-most square is Cartesian, and the only rectangle in the diagram is also Cartesian as the composition $R\to R\oplus \Sigma M\to R$ is equivalent to the identity, hence the left square in that same rectangle (the front face of the cube) is Cartesian. By a base change argument, we see that the desired back square of the cube (containing $X(e)$ and $Y(e)$) is also Cartesian, and we are done. Let us detail this base change argument. Write $I$ for the poset of nonempty subsets of $\{1,2,3\}$, ordered by inclusion, and use this poset to index the cube in (\ref{cubeforourplesure}) by setting $F_\varnothing=X(\widetilde{R})$, $F_1=X(R)$ (in the upper-right), $F_2=X(R)$ (in the centre), $F_3=Y(\widetilde{R})$, etc. As the whole cube is Cartesian we have $F_\varnothing\simeq \lim_{I_0\in I} F_{I_0}$ and as the front face is also Cartesian we have $F_2\simeq \lim \left(F_{12}\to F_{123}\gets F_{23}\right)$. These two facts, together with \cite[Ex.5.3.8]{cubicalhomotopytheory} give us the following natural chain of equivalences of spaces
\[F_\varnothing\simeq \underset{I_0\in I}{\lim} F_{I_0}\simeq \lim\left(F_2\to G_{123}\gets G_{13}\right)\simeq G_{13},\]
where $G_{123}=\lim\left(F_{12}\to F_{123}\gets F_{23}\right)\simeq F_2$ and $G_{13}= \lim\left(F_1\to F_{13}\gets F_3\right)$. This shows the back face of the cube (indexed by $\varnothing$, $\{1\}$, $\{3\}$, and $\{1,3\}$, is Cartesian.\\

Our proof of part 3 here is essentially that of \cite[Prs.17.3.9.3-4]{sag}. On the one hand, by \cite[Pr.6.2.5.2(1)]{sag} and \cite[Df.6.2.5.3]{sag}, we see that for some fixed integer $m$, an object $\FF$ of $\QCoh(X)$ is $m$-connective if and only if for all connective $\E_\infty$-rings $R$ and transformations $\eta\colon\Spec R\to X$, the object $\eta^\ast \FF$ is $m$-connective inside $\QCoh(\Spec R)\simeq \Mod_R$. Furthermore, if $\FF$ is connective and $m\geq 0$, the object $\eta^\ast \FF$ is $m$-connective if and only if the mapping space
\[\Map_{\Mod^\cn_R}(\eta^\ast\FF, N)\simeq \Map_{\Mod^\cn_R}(\tau_{\leq m}\eta^\ast \FF, N)\]
is contractible, for all connective $(m-1)$-truncated $R$-modules $N$, by the Yoneda lemma. On the other hand, the object $L_{X/Y}$ in $\QCoh(X)$ exists if and only if the functor $F\colon \Mod_\cn^X\to \Spc$, given on objects by
\begin{equation}\label{localrepctc}F(R,\eta,M)=\fib\left(X(R\oplus M)\to X(R)\underset{Y(R)}{\times} Y(R\oplus M)\right)\end{equation}
is locally almost representable. Here $\Mod_\cn^X$ denotes the $\infty$-category of triples $(A,\eta,M)$ where $A$ is a connective $\E_\infty$-ring, $\eta \in X(A)$, and $M$ is an $A$-module; see \cite[Not.17.2.4.1]{sag}. This means that we have a (locally almost; see \cite[Df.17.2.3.1]{sag}) natural equivalence for all triples $(R, \eta, M)$ in $\Mod_X^\cn$
\[F(R,\eta, M)\simeq \Map_{\Mod_R}(\eta^\ast L_{X/Y}, M)\]
where $R$ is a connective $\E_\infty$-ring, $\eta\colon\Spec R\to X$ a map in $\P(\Aff^\cn)$, and $M$ a connective $R$-module. If $L_{X/Y}$ vanishes, then we immediately see $F(R,\eta, M)$ is contractible for all triples $(R, \eta, M)$, which by part 3 implies $X\to Y$ is formally \'{e}tale, courtesy of the definition (\ref{localrepctc}) of $F$. Conversely, if $X\to Y$ is formally \'{e}tale, then $F(R,\eta, M)$ is contractible for all triples $(R,\eta,M)$, hence the mapping space
\[\Map_{\Mod_R}(\eta^\ast L_{X/Y}, M)\simeq F(R,\eta, M)\]
is contractible for all triples $(R,\eta, M)$ and $L_{X/Y}$ vanishes.
\end{proof}

There are many examples of formally \'{e}tale maps in spectral algebraic geometry.\\

Note that all formal spectral Deligne--Mumford stacks are cohesive, nilcomplete, and absolute cotangent complexes always exist, which follows by copying the proof of \cite[Cor.17.3.8.5]{sag} (the same statement for $\SpDM$), as all of the references made there also apply to $\fSpDM$.

\begin{example}[\'{E}tale morphisms of connective $\E_\infty$-rings]\label{etalemapsofrings}
Let $A\to B$ be an \'{e}tale morphism of connective $\E_\infty$-rings, then by \cite[Cor.7.5.4.5]{ha} we know $L_{B/A}$ vanishes, hence $A\to B$ is also a formally \'{e}tale morphism of $\E_\infty$-rings by \Cref{alternativefetexpressions}.
\end{example}

\begin{example}[Relatively perfect discrete $\F_p$-algebras]
Another classic example, which will not show up explicitly in this note but is at the heart of much of the work done in \cite{ec2}, is that a flat relatively perfect map of discrete commutative $\F_p$-algebras has a vanishing cotangent complex (\cite[Lm.5.2.8]{ec2}), and hence is formally \'{e}tale.
\end{example}

\begin{remark}[Noetherian, almost fp, and formally étale implies étale]\label{followuphehe}
In \Cref{differenttoluriedefinition}, we noted that formally \'{e}tale morphisms of connective $\E_\infty$-rings were not necessarily flat. However, \cite[Pr.3.5.5]{ec2} states that morphisms of (not necessarily connective) \emph{Noetherian} $\E_\infty$-rings with vanishing cotangent complex are flat. Combining this with \Cref{alternativefetexpressions}, we see formally \'{e}tale morphisms of connective Noetherian $\E_\infty$-rings are flat. It also follows (as in classical algebraic geometry, see \cite[\href{https://stacks.math.columbia.edu/tag/02HM}{02HM}]{stacks}) that formally \'{e}tale morphisms almost of finite presentation between connective Noetherian $\E_\infty$-rings are \'{e}tale. 
\end{remark}

The functor $\M_\BT$ is cohesive, nilcomplete, and admits a cotangent complex by \cite[Pr.3.2.2]{ec2}. It follows that $\M_\BTn$ (as well as the implicit base change ${\M}_{\BTn}\times \Spf A$) also satisfy these properties as $\M_\BTn\to \M_\BT$ is open---a consequence of the height of the $p$-divisible group being an open condition; see \cite[Rmk.2.1.2]{mythesis}.

\begin{example}[Spectral Serre--Tate theorem]\label{serretateexamplefet}
The spectral Serre--Tate theorem (\cite[Th.7.0.1]{ec1}) together with \Cref{alternativefetexpressions} shows that the map $[p^\infty]\colon{\M}_{\AVar_g}\to {\M}_{\BTtwog}$ is formally \'{e}tale after base change by $\Spf \Sph_p$.
\end{example}

Recall the definition of a {nonstationary} $p$-divisible group from \Cref{nonstationary}.

\begin{example}[Spectral Lubin--Tate theory]\label{lubintatespecteralformallyetaleexmaple}
For a nonstationary $p$-divisible group $\G_0$ over a discrete ring $R_0$ where $p$ is nilpotent and whose absolute cotangent complex $L_{R_0}$ is almost perfect, Lurie uses his de Rham space formalism to construct a map $\G\colon\Spf R\to {\M}_\BTn$ (\cite[Th.3.4.1]{ec2}) which is formally \'{e}tale by \cite[Cor.18.2.1.11(2)]{sag} and \Cref{alternativefetexpressions}. The $p$-divisible group $\G$ is the universal spectral deformation of $\G_0$ and $R$ its spectral deformation ring; see \Cref{deformationtheory}.
\end{example}

\begin{example}[Formal spectral completions]\label{vanishingofspfoverspec}
Let $\xx$ be a spectral Deligne--Mumford stack and $K\subseteq |\xx|$ be a cocompact closed subset, then the natural map from the formal completion of $\xx$ along $K$ (\cite[Df.8.1.6.1]{sag}) $\xx_K^\wedge\to \xx$ is formally \'{e}tale by \cite[Ex.17.1.2.10]{sag} and \Cref{alternativefetexpressions}.
\end{example}

\begin{example}[Spectral de Rham space]\label{bigderhambansingfromlurie}
Given a morphism $X\to Y$ of functors in $\P(\Aff^\cn)$, one can associate a \emph{de Rham space} $(X/Y)_\dR$ inside $\P(\Aff^\cn)$, whose value on a connective $\E_\infty$-ring is
\begin{equation}\label{derhamspaceequation} (X/Y)_\dR(R)=\underset{I}{\colim} \left( Y(R)\underset{Y(\pi_0 R/I)}{\times} X(\pi_0 R/I) \right)\end{equation}
where the colimit is taken over all nilpotent ideals $I\subseteq \pi_0 R$, which we note is a discrete filtered system; see \cite[\textsection18.2.1]{sag}. By \cite[Cor.18.2.1.11(2)]{sag}, the natural map $(X/Y)_{\dR}\to Y$ is nilcomplete, infinitesimally cohesive, and admits a vanishing cotangent complex, so by \Cref{alternativefetexpressions}, it is formally \'{e}tale.
\end{example}

This last example will help us study the moduli stack $\M_\BTn$.


\section{Deformation theory on the moduli stack of $p$-divisible groups}\label{defonmoduliofbdivsection}

Let us now apply the theory of formally \'{e}tale natural transformations to the functor ${\M}_{\BTn}$ and the categories $\C_{A_0}$ and $\C_A$ of \Cref{definitionofsites}.

\begin{notation}\label{nicefunctorthatareboring}
Write $\iota\colon \CAlg^\heartsuit\to \CAlg^\cn$ for the inclusion, a right adjoint, inducing a left adjoint $(-)^\heartsuit$ on presheaf categories. Write $\tau_{\leq 0}\colon \CAlg^\cn\to \CAlg^\heartsuit$ for the truncation functor, a left adjoint, inducing a right adjoint $\tau_{\leq 0}^\ast$ on presheaf categories. Also write $\tau_{\leq 0}$ for the composition $(-)^\heartsuit\circ \tau_{\leq 0}$ (this should seldom cause confusion). For each functor $\M$ in $\P(\Aff^\cn)$ there is a natural unit $\M\to \tau_{\leq 0}^\ast \M$ induced by the truncation $R\to \pi_0 R$ of a connective $\E_\infty$-ring $R$. The functor $\tau_{\leq 0}^\ast \M$ can also be described as the right Kan extension of $\M^\heartsuit$ along $\iota$.
\end{notation}

\begin{warn}
In \Cref{truncationssubsection}, we introduce the truncation of a locally Noetherian formal spectral Deligne--Mumford stack $\tau_{\leq 0}\x$ \`{a} la Lurie \cite[\textsection1.4.6]{sag} and we note that this is \textbf{not} equivalent to $\tau_{\leq 0}^\ast \x$. We have occasion to use both notations in this section.
\end{warn}

For mostly formal reasons, we obtain a functor $\C_A\to \C_{A_0}$.

\begin{prop}
The functor $(-)^\heartsuit\colon \C_{A}\to \P(\Aff^\heartsuit)_{/{\M}_{\BTn}^\heartsuit\times \Spf A_0}$ factors through $\C_{A_0}$.
\end{prop}

\begin{proof}
By definition, a object $\x$ of $\C_A$ is qcqs and so has an affine \'{e}tale hypercover $\u_\bullet\to \x$ ; see \Cref{etalehypercoversandsuch}. The formal spectral Deligne--Mumford stack $\tau_{\leq 0}\x=\x_0$ then lies in the essential image of $\fDM\to \fSpDM$ and hence can be considered as a classical spectral Deligne--Mumford stack. Moreover, $\x^\heartsuit$ and $\x_0^\heartsuit$ are naturally equivalent by \Cref{tunrcationsdotherightthingonFOP}. As each affine formal spectral Deligne--Mumford stack $\u_n$ is Noetherian, $\x^\heartsuit=\x_0$ has an affine \'{e}tale hypercover by $\u_\bullet^\heartsuit\to \x^\heartsuit=\x_0$ inside $\fDM$. By \Cref{etalehypercoversandsuch}, we see $\x_0$ is qcqs. From \Cref{tunrcationsdotherightthingonFOP} we also see that $\x_0$ is represented by $\tau_{\leq 0}\x$ inside $\P(\Aff^\heartsuit)$, so $\x_0$ and $\x$ have the same closed points. As $\u_0^\heartsuit$ is a Noetherian affine classical formal Deligne--Mumford stack, we also see $\x_0$ is locally Noetherian. It also follows from \Cref{biggerimplieslessthan} that $\x_0\to {\M}_{\BTn, A_0}^\heartsuit$ is formally \'{e}tale inside $\P(\Aff^\heartsuit)$. To see the cotangent complex $L$ of the map $\x_0\to {\M}_{\BTn}$ is almost perfect inside $\QCoh(\x_0)$, consider the following composition of maps in $\P(\Aff^\cn)$
\[\tau_{\leq 0} \x\simeq \x_0\to \x\to {\M}_{\BTn}\]
for which we obtain a fibre sequence in $\QCoh(\x_0)$ of the form
\[\left.L_{\x/{\M}_{\BTn}}\right|_{\x_0}\to L \to L_{\x_0/\x}.\]
By part 3 of \Cref{alternativefetexpressions}, the first term in the above fibre sequence vanishes, and our desired conclusion follows as $L_{\x_0/\x}$ is almost perfect by \Cref{truncationsofnoetherianfspdmaregood}.
\end{proof}

To see $(-)^\heartsuit$ is an equivalence (\Cref{maintheorembeforetopoi}), we will construct an explicit inverse.

\begin{mydef}\label{constructionofinverse}
Define a functor $\widetilde{\D}\colon \C_{A_0}\to \P(\Aff^\cn)_{/{\M}_{\BTn}}$ by sending an object $\G_0\colon\x_0\to {\M}_{\BTn}^\heartsuit\times \Spf A_0$ of $\C_{A_0}$ to Lurie's de Rham space of \cite[\textsection18.2.1]{sag} (and \Cref{bigderhambansingfromlurie})
\[\widetilde{\D}(\G_0)=\left(\x_0/{\M}_{\BTn}\right)_\dR.\]
\end{mydef}

\begin{theorem}\label{maintheorembeforetopoi}
The functor $\widetilde{\D}$ factors through the $\infty$-subcategory $\C_{A} \subseteq \P(\Aff^\cn)_{/{\M}_{\BTn}}$ as $\D\colon \C_{A_0} \to \C_A$. Moreover, this functor $\D$ preserves affine objects and \'{e}tale hypercovers, and is an inverse to $(-)^\heartsuit$.
\end{theorem}

The notation $\D$ is supposed to conjure the word ``deformation''.\\

This equivalence of $\infty$-categories fits into the general paradigm of spectral algebraic geometry---a well-behaved site over a classical moduli stack should be equivalent to the same site over the associated spectral moduli stack; see the example of the moduli stack of elliptic curves in \cite[Rmk.2.4.2]{ec1} and \cite[\textsection7]{ec2}, or the affine case in \cite[Th.7.5.0.6]{ha}.\\

To prove \Cref{maintheorembeforetopoi}, which will be done after \Cref{oneconnectiveisstandard}, we will use the interaction of the de Rham space technology of Lurie (\cite[\textsection18.2.1]{sag}) with formally \'{e}tale morphisms, and a representability theorem also due to Lurie.

\begin{prop}\label{formallyetaleontruncationsisnice}
Let $\x$ be a formal spectral Deligne--Mumford stack and $\x\to {\M}_{\BTn}$ be a 0-formally \'{e}tale map whose associated cotangent complex is almost perfect. Then the following natural diagram of functors in $\P(\Aff^\cn)$
\begin{equation}\label{diagramthatgotmegoing}\begin{tikzcd}
{(\x/{\M}_{\BTn})_\dR}\ar[d, "{\G_\dR}"]\ar[r, "{u}"]	&	{\tau_{\leq 0}^\ast \x}\ar[d]	\\
{{\M}_{\BTn}}\ar[r, "{d}"]			&	{\tau_{\leq 0}^\ast {\M}_{\BTn}}
\end{tikzcd}\end{equation}
is Cartesian, the natural map $\x\to(\x/{\M}_{\BTn})_\dR$ induces an equivalence when evaluated on discrete $\E_\infty$-rings, and $\G_\dR$ is formally \'{e}tale.
\end{prop}

The above proposition and its proof generalise to a wider class of functors in $\P(\Aff^\cn)$ of which we could not find a neater formulation than our leading example---we leave the reader to explore more general examples as they wish.

\begin{proof}
Recall the value of the de Rham space $(X/Y)_\dR$ on a connective $\E_\infty$-ring $R$ from (\ref{derhamspaceequation})
\[ (X/Y)_\dR(R)=\underset{I}{\colim} \left( Y(R)\underset{Y(\pi_0 R/I)}{\times} X(\pi_0 R/I) \right)\]
and define another functor $(X/Y)_\dR^0\colon \CAlg^\cn\to \Spc$ by the same formula except we index the colimit over finitely generated nilpotent ideals of $\pi_0 R$. One readily obtains a map of functors
\[(X/Y)_\dR^0\to (X/Y)_\dR\]
and we claim this map is an equivalence for $X=\x$ and $Y={\M}_{\BTn}$ in our hypotheses. Indeed, this is a key step in Lurie's proof of \cite[Pr.3.4.3]{ec2}, and the same proof holds \emph{mutatis mutandis} exchanging only $R_0$ for $\x$; the crucial step comes at the end and uses the almost perfect assumption on our cotangent complex. Writing $F_{R,I}$ for the fibre product within the colimit of (\ref{derhamspaceequation}) where $X=\x$ and $Y={\M}_{\BTn}$, we place $F_{R,I}$ into the commutative diagram of spaces
\[\begin{tikzcd}
{F_{R,I}}\ar[d]\ar[r]			&	{\x(\pi_0 R)}\ar[d, "f_{\pi_0 R}"]\ar[r]	&	{\x(\pi_0 R/I)}\ar[d, "f_{\pi_0 R/I}"]	\\
{{\M}_{\BTn}(R)}\ar[r]	&	{{\M}_{\BTn}(\pi_0 R)}\ar[r]	&	{{\M}_{\BTn}(\pi_0 R/I).}
\end{tikzcd}\]
The outer rectangle is Cartesian by definition and we claim that the right square is also Cartesian. Indeed, this follows as $I$ is finitely generated and hence is nilpotent of finite degree $n$ for some integer $n\geq 2$, and our $0$-formally \'{e}tale hypotheses can be sequentially applied to the composition of square-zero extensions
\[R\to R/I^n\to R/I^{n-1}\to \cdots \to R/I^2\to R/I.\]
This implies that the left square above is Cartesian, so the $R$-points of the de Rham space in question naturally take the form
\[\underset{I}{\colim} \left( {\M}_{\BTn}(R)\underset{{\M}_{\BTn}(\pi_0 R)}{\times} \x(\pi_0 R) \right) \simeq {\M}_{\BTn}(R)\underset{{\M}_{\BTn}(\pi_0 R)}{\times} \x(\pi_0 R)\]
as the diagram indexing our colimit is filtered. Hence (\ref{diagramthatgotmegoing}) is Cartesian. For the second statement, we use the facts that (\ref{diagramthatgotmegoing}) is Cartesian and $d$ induces equivalences when evaluated on discrete rings to see that $u$ induces an equivalence when evaluated on discrete rings. Noting that the maps
\[(\x/{\M}_{\BTn})_\dR\to \tau_{\leq 0}^\ast \x\gets \x\]
induce equivalences on discrete rings, we see the natural map $\x\to (\x/{\M}_{\BTn})_\dR$ induces an equivalence on discrete rings. Finally, to see $(\x/{\M}_{\BTn})_\dR\to {\M}_{\BTn}$ is formally \'{e}tale we refer to \Cref{bigderhambansingfromlurie}, or alternatively to the Cartesian diagram (\ref{diagramthatgotmegoing}), \Cref{truncationsgiveetale}, and \Cref{bcandsuch}.
\end{proof}

The proof of the above proposition uncovers something quite useful.

\begin{remark}[{$\D$ produces universal spectral deformations}]\label{specdeformationofD}
Recall that for each classical $p$-divisible group $\G_0\colon \Spf B_0\to {\M}_{\BTn}^\heartsuit$, such as those in $\C_{A_0}$, we can ask if there exists a \emph{universal spectral deformation} of $\G_0$ and its associated spectral deformation ring; see \Cref{deformationtheory}. It follows from the proof of \Cref{formallyetaleontruncationsisnice} above, that if $\G_0$ lies in $\C_{A_0}$, then the formal spectrum $\Spf$ of the spectral deformation ring of $\G_0$ is equivalent to the de Rham space $(\G_0\colon \Spf B_0\to {\M}_{\BTn})_\dR$. By \Cref{maintheorembeforetopoi}, we see that this de Rham space is represented by a formal spectral Deligne--Mumford stack $\Spf B$. This means that $\D(\G_0)$ is represented by the universal spectral deformation of $\G_0$.
\end{remark}

The following definition and representability theorem of Lurie is crucial. 

\begin{mydef}[{\cite[Df.18.2.2.1]{sag}}]
	A morphism $f\colon\x\to \y$ of formal spectral Deligne--Mumford stacks is called a \emph{formal thickening} if the induced map on reductions $\x^\red\to \y^\red$ is an equivalence (\cite[\textsection8.1.4]{sag}) and $f$ itself is representable by closed immersions which are locally almost of finite presentation.
\end{mydef}

\begin{theorem}[{\cite[Th.18.2.3.1]{sag}}]\label{impofrepofderham}
Let $f\colon \x\to Y$ be a map of functors in $\P(\Aff^\cn)$ such that $\x$ is a formal spectral Deligne--Mumford stack and $Y$ is nilcomplete, infinitesimally cohesive, admits a cotangent complex, and is an \'{e}tale sheaf. Suppose that $L_{\x/Y}$ is 1-connective and almost perfect. Then $(\x/Y)_\dR$ is represented by a {formal thickening} of $\x$.
\end{theorem}

Let us check that the hypotheses of this theorem are satisfied by objects in $\C_{A_0}$.

\begin{prop}\label{oneconnectiveisstandard}
	The cotangent complex $L_{\x_0/{\M}_{\BTn}\times \Spf A}$ corresponding to an object inside $\C_{A_0}$ is $1$-connective.
\end{prop}

\begin{proof}
By \cite[Cor.8.2.5.5]{sag} we may check this \'{e}tale locally on $\x_0$, so let us replace $\x_0$ with $\Spf B_0$ for some complete Noetherian discrete adic ring $B_0$. In particular, $L$ is now an almost perfect $J$-complete $B_0$-module, where $J$ is an ideal of definition for $B_0$. As $L$ is almost perfect, the fibrewise connectivity criterion of \cite[Cor.2.7.4.3]{sag} states that it suffices to check $L^\wedge_\m$ is 1-connective for every maximal ideal $\m\subseteq B_0$ which contains $J$. Moreover, considering the maps
\[\Spf (B_0)^\wedge_\m\to \Spf B_0\to \Spec B_0\]
the composition is formally \'{e}tale (as discussed for $\P(\Aff^\cn)$ by \Cref{vanishingofspfoverspec} and hence in $\P(\Aff^\heartsuit)$ by \Cref{biggerimplieslessthan}), and the latter map is unramified, so by \Cref{superbaseicproeprtoes} we see the first map is formally \'{e}tale. We may then assume $B_0$ is a complete local Noetherian ring. The morphism $\G\colon\Spf B_0\to {\M}_{\BTn}^\heartsuit \times \Spf A_0$ is formally \'{e}tale, so by the converse statement in \Cref{lubintatefetexample}, we see that $B_0$ is the classical deformation ring of $\G_\kappa$, where $\kappa$ is the residue field of $B_0$, which is necessarily perfect of characteristic $p$ by assumption. For such a pair $(\G_\kappa,\kappa)$, there exists a spectral deformation ring $B$ by \cite[Th.3.1.15]{ec2}, as $\kappa$ is perfect and $\G_\kappa$ is nonstationary by \cite[Ex.3.0.10]{ec2}, which implies $\pi_0 B\simeq B_0$ by \Cref{classicalcomesfromnonclassical}. This means the map $\Spf B_0\to {\M}_{\BTn} \times \Spf A$ in $\P(\Aff^\cn)$ factors as
\begin{equation}\label{compositionforderham}\Spf B_0\to \Spf B\to {\M}_{\BTn}\times \Spf A\end{equation}
where the first map is induced by the truncation. Associated with the above composition is the fibre sequence of complete $B_0$-modules
\[\left. L_{\Spf B/{\M}_{\BTn}\times \Spf A}\right|_{\Spf B_0}\to L_{\Spf B_0/{\M}_{\BTn}\times \Spf A}\to L_{\Spf B_0/\Spf B}.\]
The first object vanishes as $\Spf B$ is the de Rham space (see \Cref{specdeformationofD}) for the composite (\ref{compositionforderham}) and such objects always vanish; see \Cref{bigderhambansingfromlurie}. We then see the middle cotangent complex above is 1-connected and almost perfect as this holds for $L_{\Spf B_0/\Spf B}$ by \Cref{truncationsofnoetherianfspdmaregood}.
\end{proof}

\begin{proof}[Proof of \Cref{maintheorembeforetopoi}]
First, let us check $\D$ factors through $\C_A$. By \Cref{oneconnectiveisstandard}, we can apply \Cref{impofrepofderham} to see that $\D(\G_0)$ is represented by a formal thickening $\x$ of $\x_0$. To see $\x$ satisfies the conditions of \Cref{definitionofsites}, we note the following: 
\begin{itemize}
\item	${\x}$ is locally Noetherian, as it is a formal thickening of the locally Noetherian $\x_0$; see \cite[Cor.18.2.4.4]{sag}.
\item ${\x}$ is qcqs as a formal thickening of a qcqs formal spectral Deligne--Mumford stack is qcqs; see \Cref{qcqsformalthickenings}.
\item ${\x}$ has perfect residue fields at closed points as this is true for $\x_0$ and $\x_0=\tau_{\leq 0}\x$ has the same residue fields as $\x$.
\item $\G$ is formally \'{e}tale, as $L_{{\x}/{\M}_{\BTn}\times \Spf A}$ vanishes, either by \Cref{bigderhambansingfromlurie} or \Cref{formallyetaleontruncationsisnice}.
\end{itemize}

 Notice that by \cite[Cor.18.2.3.3]{sag}, if $\x_0\simeq\Spf B_0$ is affine, then the image of any $\G_0\colon \Spf B_0\to {\M}_{\BTn}^\heartsuit\times \Spf A_0$ in $\C_{A_0}$ under $\D$ is also affine. To see $\D$ is inverse to $(-)^\heartsuit$, notice the composite $(-)^\heartsuit \D$ is equivalent to the identity as $\G_0\to \D(\G_0)$ induces an equivalence on discrete rings by \Cref{formallyetaleontruncationsisnice}. For the other composition, part 1 of \Cref{alternativefetexpressions} states that the commutative diagram of spaces is Cartesian for every connective $\E_\infty$-ring $R$:
\[\begin{tikzcd}
{\x(R)}\ar[r]\ar[d, "{\G}"]				&	{\x(\pi_0R)}\ar[d]			\\
{{\M}_{\BTn} \times \Spf A(R)}\ar[r]	&	{{\M}_{\BTn}\times \Spf A(\pi_0R).}	
\end{tikzcd}\]
By \Cref{alternativefetexpressions}, we then see the natural map $\D((\G)^\heartsuit)\to \G$ is an equivalence in $\C_A$. Finally, to see $\D$ preserves \'{e}tale hypercovers, we first note this may be checked \'{e}tale locally, so take an \'{e}tale hypercover $\Spf C_0^\bullet\to \Spf B_0$ in $\C_{A_0}$ and write $\Spf C^\bullet\to \Spf B$ for its image under $\D$. From the above, we know that $\Spf C^\bullet\to \Spf B$ is an \'{e}tale hypercover on zeroth truncations, so it suffices to see each map $B\to C^n$ is \'{e}tale as morphism of $\E_\infty$-rings. Two applications of \Cref{formallyetaleontruncationsisnice} show the commutative diagram in $\P(\Aff^\cn)$
\[\begin{tikzcd}
{\y_\bullet}\ar[r]\ar[d]			&	{\x}\ar[r]\ar[d]			&	{{\M}_{\BTn}\times \Spf A}\ar[d]	\\
{\tau_{\leq 0}^\ast \y_\bullet^0}\ar[r]	&	{\tau_{\leq 0}^\ast \x_0}\ar[r]	&	{\tau_{\leq 0}^\ast {\M}_{\BTn}\times \Spf A}
\end{tikzcd}\]
consists of Cartesian squares, hence $\Spf C^n\to \Spf B$ is formally \'{e}tale by \Cref{truncationsgiveetale} and base change \Cref{bcandsuch}. It follows from \Cref{followuphehe} that $B\to C^n$ is \'{e}tale as a map of $\E_\infty$-rings, hence also as a map of adic $\E_\infty$-rings.
\end{proof}


\section{Orientations of $p$-divisible groups}\label{orientationparaphgra}

We are halfway to constructing our sheaf $\O^\top_\BTn$ of Lurie's theorem. The equivalence $\D\colon \C_{A_0} \simeq \C_A$ of \Cref{maintheorembeforetopoi} allows us to lift the classical algebraic geometry in $\C_{A_0}$ to spectral algebraic geometry. The next step, and the goal of this section, is to show that the assignment sending an affine object $\G\colon \Spf B \to \M_\BTn$ to an $\E_\infty$-ring classifying orientations of $\G^\circ$ is functorial. This should be intuitively clear, but spelling out the details also lets us set some notation.

\subsection{The sheaf of oriented $p$-divisible groups}\label{subsubsectiondefiningourguy}
We want to construct a \emph{moduli stack of oriented $p$-divisible groups of height $n$} $\M_\BTn^\o$. First, we need to recall Lurie's concept of an orientation of a {formal group} over an $\E_\infty$-ring from \cite[\textsection 4.3]{ec2}.

\begin{mydef}[{\cite[\textsection4.3.5]{ec2}}]\label{orientationsofformalgroupsetc}
Let $R$ be a complex periodic $\E_\infty$-ring and $\widehat{\G}$ be a formal group over $R$. An \emph{orientation} of $\widehat{\G}$ is an equivalence of formal groups $\widehat{\G}^\QQ_R\simeq \widehat{\G}$ over $R$, where $\widehat{\G}^\QQ_R$ is the Quillen formal group of $R$; see \Cref{quillenfg}. We write $\OrDat(\widehat{\G})$ the component of $\Map_{\FGroup(R)}(\widehat{\G}^\QQ_R, \widehat{\G})$ consisting of orientations. An \emph{orientation of a $p$-divisible group} $\G$ over a $p$-complete $\E_\infty$-ring is an orientation of $\G^\circ$, its identity component of \Cref{connectedcomponentremark}.
\end{mydef}

We will now define a moduli functor $\M^\o_\BT\colon\CAlg^p\to\Spc$ on the $\infty$-subcategory of $\CAlg$ spanned by $p$-complete $\E_\infty$-rings sending $R$ to a space of pairs $(\G,e)$ where $\G$ to is a $p$-divisible group over $R$ and $e$ is an orientation of $\G$. This space of pairs will be defined as the coCartesian unstraightening of a functor $F\colon\C\to \Cat_\infty$, so a higher categorical Grothendieck construction; see \cite[Th.3.2.0.1]{htt}.

\newcommand{\cp}{{\mathrm{cp}}}
\renewcommand{\co}{{\mathrm{cp}}}

\begin{construction}\label{howtodefineanhonestfunctor}
Let $\CAlg^\co$ be the full $\infty$-subcategory of $\CAlg^p$ spanned by $p$-complete $\E_\infty$-rings which are complex periodic. First, we define the functor $\M_\FGroup(-)$
\[\CAlg^\co\xrightarrow{\FGroup(-)} \Spc \qquad R\mapsto \FGroup(R)^\simeq\]
by sending a complex periodic $p$-complete $\E_\infty$-ring to the $\infty$-groupoid core of its associated $\infty$-category of formal groups; this assignment is a functor by \cite[Rmk.1.6.4]{ec2}. Write $F\colon\M_\FGroup\to \CAlg^\co$ for the associated left fibration. The functor $F$ has a section $\QQ$ which sends a $p$-complete complex periodic $\E_\infty$-ring $R$ to its Quillen formal group $\widehat{\G}^\QQ_R$ (\cite[Con.4.1.13]{ec2}). This assignment is functorial as $R$-homology and the cospectrum are both functors. Let $\M_\OrFGroup$ be the comma $\infty$-category $(\QQ F \downarrow \id_{\M_\FGroup})$, in other words, defined by the pullback square inside $\Cat_\infty$
\begin{equation}\label{arrowcategory}\begin{tikzcd}
{\M_\OrFGroup}\ar[rr]\ar[d]				&&	{(\M_\FGroup)^{\Delta^1}}\ar[d, "{(s,t)}"]	\\
{\M_\FGroup}\ar[rr, "{(\QQ F\times \id)\circ \Delta}"]	&&	{\M_\FGroup\times \M_\FGroup}
\end{tikzcd}\end{equation}
where $\Delta^1$ is the 1-simplex, $\Delta$ is the diagonal map, and $(s,t)$ sends an arrow in $\M_\FGroup$ to its source and target. More informally, an object of $\M_\OrFGroup$ is a complex periodic $p$-complete $\E_\infty$-ring $R$, a formal group $\widehat{\G}$ over $R$, and an orientation $e\colon \widehat{\G}^\QQ_R\simeq \widehat{\G}$ of $\widehat{\G}$. The functor
\begin{equation}\label{ihavetocheckthisdotdotdot}\M_\OrFGroup\to \M_\FGroup,\qquad (R,\widehat{\G}, e)\mapsto (R,\widehat{\G})\end{equation}
is a left fibration with associated functor
\[\M_\FGroup\to \Spc\qquad (R,\widehat{\G})\mapsto \OrDat(\widehat{\G}).\]
Consider the composition
\[F_\BT\colon \CAlg^\co \to \CAlg \xrightarrow{\tau_{\geq 0}} \CAlg^\cn \xrightarrow{\M_\BT} \Spc\]
of the natural inclusion, the connective cover functor, and the moduli functor $\M_\BT$ of \Cref{moduliofpdivdef}. Write the associated left fibration as $\M^\co_\BT \to \CAlg^\co$. The natural transformation $(-)^\circ \colon \BT(R) \to \FGroup(R)$ induces a functor $\M_\BT^\co \to \M_\FGroup$ of left fibrations over $\CAlg^\co$. We then define $\M^\co_\OrBT$ as the left fibration over $\CAlg^\co$ in the pullback of left fibrations
\[\begin{tikzcd}
{\M_\OrBT^\co}\ar[r]\ar[d]		&	{\M_\OrFGroup}\ar[d]	\\
{\M_\BT^\co}\ar[r, "{(-)^\circ}"]	&	{\M_\FGroup.}
\end{tikzcd}\]
As (\ref{ihavetocheckthisdotdotdot}) is a left fibration, then $\M_\OrBT^\co\to \M_\BT^\co$ is also a left fibration by base change. It follows that the left fibration $\M_\OrBT^\co\to \CAlg^\co$ has associated functor
\[F_\OrBT^\co\colon \CAlg^\co\to \Spc\qquad R\mapsto \OrBT(R),\]
as desired. By the same arguments, there is a version for $p$-divisible groups of a fixed height $\M_\OrBTn\to \CAlg^\co$ and its associated functor.
\end{construction}

\begin{remark}\label{remark4110}
Given a morphism $\varphi\colon A\to B$ in $\CAlg^p$, then if $A$ is complex periodic, $B$ is also complex periodic; this is the content of \cite[Rmk.4.1.10]{ec2}. Indeed, a complex orientation for $B$ is given by post-composing a complex orientation for $A$, so a morphism of spectra $\Sigma^{-2}\mathbf{CP}^\infty \to A$, with $\varphi$. The fact that $A$ is weakly $2$-periodic means that $\Sigma^2 A$ is a locally free $A$-module of rank $1$, so for some Zariski cover $A\to A'$, there is an isomorphism of $A'$ modules $A'\otimes_A \Sigma^2 A \simeq A'$. The $B$-module $\Sigma^2 B$ is then also locally free of rank $1$ as over the Zariski cover $B\to B\otimes_A A'=B'$ we have equivalences
\[\Sigma^2 B' \simeq B' \otimes_{A'} \Sigma^2 A' \simeq B' \otimes_{A'} (A'\otimes_A \Sigma^2 A) \simeq B' \otimes_{A'} (A') \simeq B'.\]
Hence $B$ is weakly $2$-periodic and therefore complex periodic.
\end{remark}

 We will use this remark in the following definition.

\begin{mydef}\label{orientationsandsuch}
Define a functor $\M_\BT^\o\colon \CAlg^p\to \Spc$ first on $\CAlg^\co$ as $F_\BT^\co$ and as the empty space on objects in $\CAlg^p$ which are not complex periodic. In other words, $\M_\BT^\o$ is the assignment
\[R\mapsto\begin{cases} \OrBT(R)^\simeq & \mbox{if $R$ is complex periodic}\\ \varnothing & \mbox{if $R$ is not complex periodic.}\end{cases}\]
Let us also write $\M_\BT\colon \CAlg^p \to \Spc$ for the functor which sends $R$ to $\BT(R)^\simeq$, a slight variant of \Cref{moduliofpdivdef}. From \Cref{howtodefineanhonestfunctor}, there is a natural transformation $\Omega \colon \M^\o_{\BT} \to \M_\BT$ of functors in $\P(\Aff^p)$ which on $R$-valued points sends a pair $(\G,e)$ to $\G$. Define the presheaf $\M_\BTn^\o$ in $\P(\Aff^p)$ by the Cartesian square
\[\begin{tikzcd}
{\M_\BTn^\o}\ar[r]\ar[d, "{\Omega}"]	&	{\M_\BT^\o}\ar[d, "{\Omega}"]	\\
{\M_\BTn}\ar[r]					&	{\M_\BT.}
\end{tikzcd}\]
\end{mydef}

The notation $\Omega$ is reminiscent of the word ``orientation''. 


\subsection{Orientation classifiers}
It is our goal now to try and understand \emph{universal} orientations and their relation to $\M_\BTn^\o$. We would like to construct a functor $\OO_\BTn^\o\colon \C_A^\op \to \CAlg$ which we think of as the structure sheaf for $\M^\o_\BTn$. In other words, it should send an affine object $\G\colon \Spf B \to \M_\BTn$ to the $p$-completion of the orientation classifier $\widehat{\OO}_{\G^\circ}$ of $\G^\circ$. Let us recall this notion.

\begin{mydef}[Orientation classifier {(\cite[Df.4.3.14]{ec2})}]
	Let $R$ be an $\E_\infty$-ring and $\widehat{\G}$ be a formal group over $R$. Then the \emph{orientation classifier} $\OO_{\widehat{\G}}$ of $\widehat{\G}$ is the corepresenting $R$-algebra for the functor $\CAlg_R\to \Spc$ sending $A$ to $\OrDat(\widehat{\G}_A)$.
\end{mydef}

Orientation classifiers help us to recognise that $\Omega\colon \M^\o_\BTn \to \M_\BTn$ is affine.

\begin{prop}\label{identifysomepushforward}
Given a $p$-complete $\E_\infty$-ring $R$ and a $p$-divisible group $\G$ of height $n$ over $R$, then the functor $P$ in the Cartesian square
\[\begin{tikzcd}
{P}\ar[r]\ar[d]		&	{\M^\o_\BTn}\ar[d, "\Omega"]	\\
{\Spec R}\ar[r, "\G"]	&	{\M_\BTn}
\end{tikzcd}\]
in $\P(\Aff^p)$ is represented by $\Spec\widehat{\OO}_{\G^\circ}$, the spectrum of the $p$-completion of the orientation classifier $\OO_{\G^\circ}$ for $\G^\circ$.
\end{prop}

\begin{proof}
It suffices to show that the commutative square of presheaves of $p$-complete $\E_\infty$-rings
\begin{equation}\label{cartesiansquareorientationthingy}\begin{tikzcd}
{\Spec \widehat{\OO}_{\G^\circ}}\ar[r]\ar[d]	&	{\M_\BTn^\o}\ar[d, "{\Omega}"]	\\
{\Spec R}\ar[r, "{\G}"]				&	{\M_\BTn}
\end{tikzcd}\end{equation}
is Cartesian. Fix a $p$-complete $\E_\infty$-ring $A$ and evaluate the above diagram at $A$. If there are no maps of $p$-complete $\E_\infty$-rings $R\to A$, then the two left-most spaces are empty and we are done, so let us then fix a map $\psi\colon R\to A$. We then note the chain of natural equivalences between the fibres of the vertical morphisms from left to right
\[\Spec \widehat{\OO}_{\G^\circ}(A)\simeq \Map_{\CAlg}(\OO_{\G^\circ},A)\simeq \OrDat(\G^\circ_A)\simeq \{\psi\}\underset{\M_\BTn^\un(A)}{\times} \M_\BTn^\o(A).\]
The first equivalence follows as $p$-completion is a left adjoint, the second from \cite[Pr.4.3.13]{ec2}, and the third from the construction of $\Omega$; see \Cref{howtodefineanhonestfunctor}. As these equivalences are natural in $A$, this shows (\ref{cartesiansquareorientationthingy}) is Cartesian. 
\end{proof}

\begin{remark}\label{orientationclassifiersaregood}
	Let $f\colon R\to R'$ be a morphism of connective $\E_\infty$-rings and $\widehat{\G}$ be a formal group over $R$. Then \cite[Rmk.6.4.2]{ec2} states that the orientation classifier for $\widehat{\G}_{R'}=f^\ast \widehat{\G}$ is the base change of $f$ along the canonical map $R \to \OO_{\widehat{\G}}$. In particular, this implies that the assignment of an $R$-algebra $A$ to $\OO_{\widehat{\G}_A}$ is an fpqc hypersheaf as if $A\to A^\bullet$ was a hypercover, then by base change $\OO_{\widehat{\G}_A} \to \OO_{\widehat{\G}_{A^\bullet}}$ is an fpqc hypercover too.
\end{remark}

Just as in the previous section, we construct the functor $\OO_\BTn^\o$ in a few formal steps. First, we need an algebraisation process on affines.

\begin{mydef}\label{baffnotation}
Write $\C_{A_0}^\aff$ (resp.\ $\C_{A}^\aff$) for the full $\infty$-subcategory of $\C_{A_0}$ (resp.\ $\C_{A}$) spanned by affine objects.
\end{mydef}

\begin{remark}\label{algebraisationremark}
	Let $A$ be an adic $\E_\infty$-ring and $\G$ be a $p$-divisible group over $\Spf A$. There is a natural map of spaces $\M_\BT(\Spec A)\to \M_\BT(\Spf A)$. We call any $\G^\alg$ over $\Spec A$ an \emph{algebraisation} of $\G$ if under the above map it is sent to $\G$. By \cite[Th.3.2.2(4)]{ec2}, the above map is an equivalence of spaces if $A$ is complete with respect to its ideal of definition. In particular, if $A$ is complete there always exists a unique algebraisation. In particular, this yields an equivalences of $\infty$-categories
\[\fSpDM^\aff_{/\M_\BTn} \simeq \SpDM^\aff_{/\M_\BTn}\]
which we think about as sending a pair $(\Spf B, \G)$ to $(\Spec B^\wedge_J, \G_{B^\wedge_J})$ where $J$ is a finitely generated ideal of definition for $B$. As mentioned in \Cref{completioninimplicitlol}, this completion is often left implicit.
\end{remark}

\begin{mydef}\label{sheafactuallyonb}
We define the functor $\OO^\aff_\BTn\colon (\C_{A}^\aff)^\op \to \CAlg$ as the composition
\[\C_A^\aff \subseteq  \fSpDM^{\aff,p}_{/\M_\BT} \simeq \SpDM^{\aff,p}_{/\M_\BT} \xrightarrow{\Omega^\ast} \SpDM^{\nc,\aff,p}_{/\M_\BT^\o} \to \SpDM^{\nc,\aff,p} =(\CAlg^p)^\op \to \CAlg^\op \]
given by the defining inclusion, the equivalence of \Cref{algebraisationremark}, base change along $\Omega$ (well-definedness comes from \Cref{identifysomepushforward}) into not necessarily connective spectral Deligne--Mumford stacks (\cite[Df.1.4.4.2]{sag}), the forgetful functor, and the inclusion. Moreover, all functors involved preserve étale hypercovers, the only questionable functor $\Omega^\ast$ being dealt with in \Cref{orientationclassifiersaregood}, so the composite $\OO^\aff_\BTn$ is an étale hypersheaf. We define $\OO_\BTn^\o\colon \C_A^\op\to \CAlg$ by right Kan extension along the inclusion $(\C_A^\aff)^\op\to \C_A^\op$. As a right Kan extensions preserve limits, we see $\OO_\BTn^\o$ is an \'{e}tale hypersheaf.
\end{mydef}

\begin{remark}\label{rightkanextensionevaluation}
The right Kan extension defining $\OO_\BTn^\o$ on $\C_A$ can be made more explicit. Indeed, by assumption, each object $\x$ in $\C_A$ is qcqs, so by \Cref{etalehypercoversandsuch}, we have an \'{e}tale hypercover $\y_\bullet\to \x$ such that each $\y_n=\Spf B_n$ is affine. The fact that $\OO_\BTn^\o$ is an \'{e}tale hypersheaf yields
\[\OO_\BTn^\o(\x)\simeq \lim\left(\OO_\BTn^\aff(\Spf B^0)\Rightarrow \OO_\BTn^\aff(\Spf B^1)\Rrightarrow \cdots\right).\]
By \Cref{identifysomepushforward}, the terms in the above limit take a known form.
\end{remark}

We have just defined $\OO_\BTn^\o$, so we are very close to our definition of $\O^\top_\BTn$ and proving \Cref{maintheorem}.


\section{The sheaf $\O^\top_\BTn$ and a proof of \Cref{maintheorem}}\label{finalsectioninproof}
The definition of $\O^\top_\BTn$ mirrors Lurie's definition of $\O^\top$ (\cite[\textsection7.3]{ec2}) and the proof that this definition satisfies \Cref{maintheorem} also follows Lurie's proof. 

\begin{mydef}\label{definitionofghmlsheaf}
Fix an adic $\E_\infty$-ring $A$ as in \Cref{fixeda}. Let $\O^\top_\BTn$ be the presheaf of $\E_\infty$-rings on $\C_{A_0}$ defined by the composition
\begin{equation}\label{compositiondefinitonofsheaf} \O^\top_\BTn\colon \C_{A_0}^\op\xrightarrow{\D^\op} \C_{A}^\op\xrightarrow{\OO^\o_\BTn} \CAlg.\end{equation}
It follows from \Cref{maintheorembeforetopoi} and \Cref{sheafactuallyonb} that $\O^\top_\BTn$ is an \'{e}tale hypersheaf.
\end{mydef}

\begin{remark}[Unpacking $\O^\top_n$]\label{inotherwordsremkar}
In other words, one first calculates the universal spectral deformation of $\G_0\colon\x_0\to {\M}_{\BTn}^\heartsuit$ giving $\D(\G_0)=\G$ (\Cref{specdeformationofD}), then the identity component $\G^\circ$ of $\G$, and $\O^\top_\BTn(\G_0)$ is then the $p$-completion orientation classifier of $\G^\circ$; we will see in the proof of \Cref{maintheorem} below that this final $p$-completion is unnecessary. This is the same recipe Lurie uses to define the Lubin--Tate theories $E_n$ in \cite[\textsection5]{ec2}. The fact that $\O^\top_n(\G_0)$ is an orientation classifier gives us the canonical isomorphisms of \Cref{maintheorem} when $\x_0=\Spf B_0$ is affine.
\end{remark}

Before we prove \Cref{maintheorem}, we will need the following connection between formally \'{e}tale morphisms and universal deformations. The following is analogous to \Cref{lubintatefetexample} and our proof follows that of \cite[Pr.7.4.2]{ec2} (we will even copy some of Lurie's notation).

\begin{prop}\label{wehavesomeuniversaldeformations}
Recall we are implicitly working over $\Spf A$ from \Cref{fixeda}. Let $\G\colon\Spf B\to {\M}_{\BTn}$ be formally \'{e}tale map where $B$ is a complete adic Noetherian $\E_\infty$-ring with ideal of definition $J$. Fix a maximal ideal $\m\subseteq\pi_0 B$ containing $J$ such that $\pi_0 B/\m$ is perfect of characteristic $p$. Then the $p$-divisible group $\G_{B^\wedge_\m}$ is the universal spectral deformation of $\G_\k$ (in the sense of \cite[Df.3.1.11]{ec2}), where $\k$ is the residue field of $B_\m$.
\end{prop}

\begin{proof}
First, the fact that $\k$ is perfect of characteristic $p$ means we can combine \cite[Ex.3.0.10]{ec2} with \cite[Th.3.1.15]{ec2} to obtain the spectral deformation ring $R^\un_{\G_\k}=B^\un$ with a universal $p$-divisible group $\G^\un$. By definition, $\G_{B^\wedge_\m}$ is a deformation over $\G_\k$ (\cite[Df.3.0.3]{ec2}) and the universality of the pair $(B^\un, \G^\un)$ yields a canonical continuous morphism of adic $\E_\infty$-rings $B^\un\xrightarrow{\al} B^\wedge_\m=\widehat{B}$ inducing the identity on the common residue field $\kappa$. Let us write $\C$ for the full $\infty$-subcategory of $(\CAlg^\cn_\ad)_{/\k}$ spanned by complete local Noetherian adic $\E_\infty$-rings whose augmentation to $\k$ exhibits $\k$ as its residue field. By \cite[Th.3.1.15]{ec2}, we see $B^\un$ belongs this $\infty$-subcategory $\C$. By assumption and \cite[Cor.7.3.8.3]{sag}, which states that the completion of a Noetherian $\E_\infty$-ring an ideal is also Noetherian, we see that $\widehat{B}$ also belongs to $\C$. We will then use a Yoneda-style argument to see that $\al$ is an equivalence in $\C$. Consider an arbitrary object $C$ of $\C$ and the induced map
\begin{equation}\label{desiredequivalence} \Map^\cont_{\CAlg_{/\k}}(\widehat{B}, C)\xrightarrow{\al^\ast} \Map^\cont_{\CAlg_{/\k}}(B^\un, C).\end{equation}
By writing $C$ as the limit of its truncations $C \simeq \lim \tau_{\leq n} C$ and using the universal property of a limit, we reduce ourselves to the case where $C$ is truncated. Similarly, by writing $\pi_0 C$ as a limit of Artinian subrings of $\pi_0 C$ we are further reduced to the case when $\pi_0 C$ is Artinian.\footnote{Our conventions demand that local adic $\E_\infty$-rings have their topology determined by the maximal ideal.} In this situation, when have a finite sequence of maps
\begin{equation}\label{chainofsqzeroextns} C=C_m\to C_{m-1}\to \cdots\to C_1\to C_0=\k\end{equation}
where each map is a square-zero extension by an almost perfect connective module. Hence, it would suffice to show that for every $C\to \k$ in $\C$, and every square-zero extension $\widetilde{C}\to C$ of $C$ by an almost perfect connective $C$-module, with $\widetilde{C}$ also in $\C$, the natural diagram of spaces
\begin{equation}\label{deofrmationdiagram}\begin{tikzcd}
{\Map^\cont_{\CAlg_{/\k}}(\widehat{B}, \widetilde{C})}\ar[r]\ar[d]	&	{\Map^\cont_{\CAlg_{/\k}}(B^\un, \widetilde{C})}\ar[d]	\\
{\Map^\cont_{\CAlg_{/\k}}(\widehat{B}, {C})}\ar[r]			&	{\Map^\cont_{\CAlg_{/\k}}(B^\un, {C})}
\end{tikzcd}\end{equation}
is Cartesian. Indeed, if the above square is Cartesian, then when $C=\k$ lower-horizontal map is the unique map between two contractible spaces, hence the upper map is an equivalence too. Continuing along the finite chain of square-zero extensions (\ref{chainofsqzeroextns}) and using the fact that (\ref{deofrmationdiagram}) is Cartesian, we inductively conclude that the desired map (\ref{desiredequivalence}) is an equivalence. \\

It suffices then to show that (\ref{deofrmationdiagram}) is Cartesian. Notice that for any $D$ in $\C$, this adic $\E_\infty$-ring is complete with respect to the kernel of its augmentation $D\to \kappa$. As $\widehat{B}$ is the $\m$-completion of $B_\m$, we see that for any $D$ in $\C$ the map
\[\Map^\cont_{\CAlg_{/\k}}(\widehat{B}, D)\xrightarrow{\simeq}  \Map^\cont_{\CAlg_{/\k}}(B_\m, D)\]
induced by $B_\m\to \widehat{B}$, is an equivalence. Moreover, for any $D$ inside $\C$ we have the following natural identifications:
\begin{equation*}
\begin{split}
\Map^\cont_{\CAlg_{/\k}}(B^\un, D) & \simeq \underset{B^\un\to \k}{\fib}\left( \Map_\CAlg^\cont(B^\un, D)\to \Map_\CAlg^\cont(B^\un, \k) \right) \\
 & \simeq \underset{B^\un\to \k}{\fib}\left( (\Spf B^\un)(D)\to (\Spf B^\un)(\k) \right) \\
 & \simeq \underset{\G^\un}{\fib}\left( \Def_{\G_\k}(D)\to \Def_{\G_\k}(\k) \right) \\
 & \simeq \Def_{\G_\k}(D, (D\to \k)) \\
 & \simeq{\mathrm{BT}^p_n}(D)\underset{{\mathrm{BT}^p_n}(\k)}{\times}\{\G_\k\}.
\end{split}
\end{equation*}

The first equivalence is a categorical fact about over/under categories, the second is the identification of the $R$-valued points of $\Spf B^\un$ (\cite[Lm.8.1.2.2]{sag}), the third is from universal property of spectral deformation rings (\cite[Th.3.1.15]{ec2}), and the fourth and fifth can be taken as two alternative definitions of $\Def_{\G_\k}(D,(D\to \k))$ (\cite[Df.3.0.3 \& Rmk.3.1.6]{ec2}). The natural equivalences show that (\ref{deofrmationdiagram}) is equivalent to the upper-left square in the commutative diagram of spaces
\begin{equation}\label{Istolethisfromjacob}\begin{tikzcd}
{(\Spf B_\m)(\widetilde{C}) \simeq \Map^\cont_{\CAlg_{/\k}}(B_\m, \widetilde{C})}\ar[r]\ar[d]	&	{{\mathrm{BT}^p_n}(\widetilde{C})\underset{{\mathrm{BT}^p_n}(\k)}{\times} \{\G_\k\}}\ar[d]\ar[r]	&	{{\mathrm{BT}^p_n}(\widetilde{C})^\simeq}\ar[d]	\\
{(\Spf B_\m)({C}) \simeq \Map^\cont_{\CAlg_{/\k}}(B_\m, {C})}\ar[r]				&	{{\mathrm{BT}^p_n}({C})\underset{{\mathrm{BT}^p_n}(\k)}{\times} \{\G_\k\}}\ar[r]\ar[d]			&	{{\mathrm{BT}^p_n}({C})^\simeq}\ar[d]			\\
											&	{\{\G_\k\}}\ar[r]											&	{{\mathrm{BT}^p_n}({\k})^\simeq.}
\end{tikzcd}\end{equation}
The lower-right square and outer-right rectangle are both Cartesian by definition, so the upper-right square is Cartesian. Moreover, the upper-outer rectangle is Cartesian as $\Spf B \to \M_\BTn$, hence also $\Spf B_\m \to \M_\BTn$, is formally étale. This implies the upper-left square is also Cartesian, and we are done.
\end{proof}

We can now prove \Cref{maintheorem}; we follow Lurie's proof of \cite[Th.7.0.1]{ec2} closely.

\begin{proof}[Proof of \Cref{maintheorem}]
In \Cref{definitionofghmlsheaf} we defined an \'{e}tale hypersheaf $\O^\top_\BTn$ on $\C_{A_0}$ taking values in $\E_\infty$-rings. It remains to show that when restricted to objects $\G_0\colon\Spf B_0\to {\M}_{\BTn}^\heartsuit$ in $\C_{A_0}^\aff$, the $\E_\infty$-ring $\EE=\O_\BTn^\top(\G_0)$ has the expected properties 1-4 of \Cref{maintheorem}. Let us use \Cref{inotherwordsremkar} to break a part $\O^\top_\BTn$. Under $\D$, the object $\G_0$ is sent to the affine object $\G_\un\colon \Spf B\to {\M}_{\BTn}$ of $\C_{A}^\aff$ such that $\pi_0 B\simeq B_0$ and $\G_\un$ is equivalent to $\G_0$ when restricted to $\Spf B_0$; see \Cref{constructionofinverse}. The $\E_\infty$-ring $\EE$ is the $p$-completion of the orientation classifier of the identity component $\G^\circ$ of $\G$, denoted by $\OO_{\G^\circ}$. First, we will argue that the $\E_\infty$-ring $\OO_{\G^\circ}$ satisfies the desired properties 1-3, and then for $\EE$.\\

Firstly, note that as $\OO_{\G^\circ}$ is an orientation classifier, \cite[Pr.4.3.23]{ec2} states that $\OO_{\G^\circ}$ is complex periodic (we will discuss Landweber exactness at the very end). By \Cref{remark4110}, it follows that $\EE$ is complex periodic as there is a map of $\E_\infty$-rings $\OO_{\G^\circ}\to\EE$.\\

Recall from \cite[\textsection6.4.1]{ec2}, that a formal group $\widehat{\G}$ over a connective $\E_\infty$-ring $R$ is \emph{balanced} if the unit map $R\to \OO_{\widehat{\G}}$ induces an equivalence on $\pi_0$ and the homotopy groups of $\OO_{\widehat{\G}}$ are concentrated in even degree. We claim that $\G^\circ$ is balanced over $B$. To see this, we use \cite[Rmk.6.4.2]{ec2} (twice) to reduce ourselves to showing that $\G^\circ_{B^\wedge_\m}$ is balanced over $B^\wedge_\m$ for every maximal ideal $\m\subseteq\pi_0B\simeq B_0$; these ideals contain $J$ as $B_0$ is $J$-complete. By \Cref{wehavesomeuniversaldeformations}, we see $\G_{B^\wedge_\m}$ is the universal spectral deformation of $\G_\k$, where $\k$ is the residue field of $B^\wedge_\m$, and a powerful statement of Lurie \cite[Th.6.4.6]{ec2} then implies the identity component $\G^\circ_{B^\wedge_\m}$ of $\G_{B^\wedge_\m}$ is balanced.\\

As $\G^\circ$ is balanced over $B$ and we have already proven condition 1 of \Cref{maintheorem}, the $p$-completion map of $\E_\infty$-rings $\OO_{\G^\circ} \to \EE$ is an equivalence. Indeed, as $\G^\circ$ is balanced over $B$, we have an equivalence $B_0 \simeq \pi_0 \OO_{\G^\circ}$. We know that $B_0$ is $J$-complete by assumption, so we want to now observe that $\OO_{\G^\circ}$ is also $J$-complete. This would follow if all of its homotopy groups are $J$-complete. As $\G^\circ$ is balanced over $B$ and we have already proven condition 1 of \Cref{maintheorem}, we see that the homotopy groups of $ \OO_{\G^\circ}$ are either zero or line bundles over $B_0$. However, line bundles over complete discrete rings are complete in a classical sense, and also in a spectral sense as $B_0$ is Noetherian by \cite[Pr.7.3.6.1]{sag}. Therefore, $\OO_{\G^\circ}$ is also $J$-complete, hence also $\m_A$-complete and also $p$-complete. Therefore the $\E_\infty$-rings $\OO_{\G^\circ} \simeq \EE$ are naturally equivalent. The fact that $\G^\circ$ is balanced over $B$ and this equivalence then proves conditions 2 and 3, except for the identification of $\pi_{2k}\EE$.\\

For condition 4, \cite[Pr.4.3.23]{ec2} states that the canonical orientation of the $p$-divisible group $\G$ over $\EE$ supplies us with an equivalence $\widehat{\G}^\QQ_\EE\xrightarrow{\simeq} \G^\circ$ between the Quillen formal group of $\EE$ and the identity component of $\G$. In particular, this implies the classical Quillen formal group $\widehat{\G}^{\QQ_0}_\EE$ is isomorphic to the formal group $\G_0^\circ$ after an extension of scalars along the unit map $B_0\simeq \pi_0 B\to \pi_0 \EE$. As $\G^\circ$ is a balanced formal group over $B$, this unit map is an isomorphism, giving us property 4.\\

Let us now round off condition 3 and calculate $\pi_{2k}\EE$. We know that $\EE$ is weakly 2-periodic, that the $p$-divisible group $\G$ over $\EE$ comes equipped with a canonical orientation and hence a chosen equivalence of locally free $\EE$-modules of rank 1 $\be\colon \omega_{\G}\to \Sigma^{-2}\EE$, and that there is an equivalence of $B_0$-modules $\pi_0 \omega_{\G}\simeq\omega_{\G_0}$. Combining these facts, we obtain the desired computation:
\[\pi_{2k}\EE\simeq (\pi_2 \EE)^{\otimes k}\simeq (\pi_0 \omega_{\G})^{\otimes k}\simeq \omega_{\G_0}^{\otimes k}.\]
Finally, to finish condition 1 and the Landweber exactness of $\EE$, we appeal directly to Behrens--Lawson's arguments in \cite[Lm.8.1.6 \& Cor.8.1.7]{taf}, as they are checking the same conditions on a sheaf with the same properties as ours above.
\end{proof}

\begin{remark}\label{comparison}
Let us close this section by stating that there have, of course, been other iterations of Lurie's theorem, see \cite[Th.8.1.4]{taf} and \cite[\textsection6.7]{handbooktmf}, and the statements made there are not identical to our \Cref{maintheorem}. Lurie's theorem is an existence theorem though, so potentially two different proofs of Lurie's theorem might conjure up two different sheaves. The universality of the above construction of $\O^\top_\BTn$ and the statement of Lurie's theorem shows that at least each section of two potentially different sheaves are equivalent as $\E_\infty$-rings, but we do not go as far here to claim some uniqueness of $\O^\top_\BTn$ as done for $\O^\top$ in \cite{uniqueotop}. In the section to follow, detailing applications of Lurie's theorem, we will see that \Cref{maintheorem} applies to all of the expected examples. We can construct Lubin--Tate theories, $\TMF$, and $\TAF$, all using \Cref{maintheorem}. Practically speaking, there is no difference between any forms of Lurie's theorem we can find.
\end{remark}

\section{Applications of Lurie's theorem}\label{applicationsandactions}
To advertise Lurie's theorem to a wider audience and lay some (known) groundwork for future applications, let us now discuss how the titular theorem of this article can be used. A vast majority of the applications below can be found in either \cite{taf}, \cite{ec2}, or \cite[\textsection6.7]{handbooktmf}, in some form or another. Further applications to $\TMF$ can be found in \cite{heckeontmf}.


\subsection{Lubin--Tate and Barsotti--Tate theories}\label{ltbttheoryeiessection}
First, we want to apply Lurie's theorem to the Lubin--Tate deformation theory of \Cref{lubintateexample}. See \cite[\textsection5]{ec2} for Lurie's approach to spectral Lubin--Tate theory, which plays a key role here.\\

For now, we can set $A=\Sph_p$ and $A_0 =\Z_p$ from \Cref{fixeda}.

\begin{prop}\label{lubintateexampleunderway}
Let $\widehat{\G}_0$ be a formal group of exact height $n$ over a perfect field $\kappa$ and $\G_0$ for the associated formal $p$-divisible group over $\kappa$ whose identity component is equivalent to $\widehat{\G}_0$; the existence of such a $\G_0$ is guaranteed by \cite[Pr.4.4.22]{ec2}. Write $\G$ for the classical universal deformation of $\G_0$, which is a $p$-divisible group over the discrete ring $R^\LT$. The object $\G\colon \Spf R^\LT\to {\M}_{\BTn}^\heartsuit\times \Spf \Z_p$ lies in $\C_{\Z_p}$. Moreover, there is an equivalence of $\E_\infty$-rings $\O_\BTn^\top(\G)\simeq E_n$ where $E_n=E(\widehat{\G}_0)$ is the Lubin--Tate $\E_\infty$-ring of $\widehat{\G}_0$; see \cite[\textsection5]{ec2}.
\end{prop}

This will follow from a more general family of $p$-divisible groups in $\C_{\Z_p}$.

\begin{prop}\label{generalbtspectra}
Let $R_0$ be a discrete Noetherian $\F_p$-algebra such that the Frobenius endomorphism on $R_0$ is finite and $\G_0$ be a nonstationary $p$-divisible group of height $n$ over $R_0$. Write $R$ for the universal spectral deformation adic $\E_\infty$-ring of $\G_0$ from \cite[Th.3.4.1]{ec2} and assume the residue fields of $\pi_0 R$ are perfect of characteristic $p$. Then the morphism $\G\colon\Spf \pi_0 R\to {\M}_{\BTn}^\heartsuit$ lies in $\C_{A_0}$. 
\end{prop}

Key to Goerss--Hopkins--Miller construction of Lubin--Tate spectra is the vanishing of the algebraic cotangent complex $L^\alg_{\Spf \pi_0 E_n/\M_\BTn}$ as this implies the vanishing of all of the necessary obstruction groups; see \cite[7.3-7.6]{gh04}. In our proof of \Cref{generalbtspectra}, this condition will be consolidated into the fact that $\Spf \pi_0 R \to \M_\BTn^\heartsuit$ is formally étale which is proven in \Cref{lubintatefetexample}. The fact that the universal spectral deformation $(R,\G)$ has vanishing relative cotangent complex is a consequence of \Cref{maintheorembeforetopoi}, hence is hidden in our construction of $\O^\top_\BTn$.\\

The $\E_\infty$-rings produced by applying $\O^\top_\BTn$ to the $p$-divisible groups $\G$ occurring in \Cref{generalbtspectra} seem interesting enough to name.

\begin{mydef}
Let $R_0$, $\G_0$, and $\G$ be as in \Cref{generalbtspectra}. We call $\O^\top_\BTn(\G) = E(R_0,\G_0)$ the \emph{Barsotti--Tate $\E_\infty$-ring} associated with $(R_0,\G_0)$.
\end{mydef}

\begin{proof}[Proof of \Cref{generalbtspectra}]
Let us first see $\G$ lies in $\C_{A_0}$ by checking the conditions of \Cref{definitionofsites}. It is shown in \Cref{lubintatefetexample} that the morphism $\G$ is formally \'{e}tale. As $R_0$ is Noetherian, then \cite[Th.3.4.1(6)]{ec2} tells us that $R$ and hence also $\pi_0 R$ are Noetherian as well. We are left to check a particular cotangent complex is almost perfect. Consider the maps in $\P(\Aff^\cn)$
\[\Spf \pi_0 R\to \Spf R\to {\M}_{\BTn}\]
and the associated fibre sequence of complete $\pi_0 R$-modules
\[\left. L_{\Spf R/{\M}_{\BTn}}\right|_{\Spf \pi_0 R}\to L_{\Spf \pi_0 R/{\M}_{\BTn}}\to L_{\Spf \pi_0 R/\Spf R}.\]
By construction (\cite[Pr.3.4.3]{ec2}), $R$ corepresents the de Rham space of the map $\Spec R_0\to \M_\BT$, or equivalently, the de Rham space of $\Spec R_0\to {\M}_{\BTn}$, as $R_0$ is an $\F_p$-algebra and $\G_0$ is of height $n$. As $R$ represents this de Rham space and using \Cref{bigderhambansingfromlurie}, we see that $L_{\Spf R/{\M}_{\BTn}}$ vanishes. Hence $L_{\Spf \pi_0 R/{\M}_{\BTn}}$ is almost perfect as $L_{\Spf \pi_0 R/\Spf R}$ is almost perfect; see \Cref{truncationsofnoetherianfspdmaregood}. It follows that $\G$ lies in $\C_{A_0}$.
\end{proof}

\begin{proof}[Proof of \Cref{lubintateexampleunderway}]
The fact that $\G$ lies in $\C_{\Z_p}$ follows from \Cref{generalbtspectra}. The fact that $\O_\BT^\top(\G)$ is equivalent to $E_n$ follows as the universal spectral deformation of $\G_0$ is given by $\D(\G)$ (\Cref{generalbtspectra}) and the orientation classifier of $\D(\G)$ is $E_n$ (\cite[Cor.6.0.6]{ec2}).
\end{proof}

From the functorality of $\O^\top_\BTn$ we obtain an action of the automorphism group $\GG(\kappa,\widehat{\G}_0)$ of the pair $(\kappa,\widehat{\G}_0)$ on the $\E_\infty$-ring $E_n$. In other words, $E_n$ obtains an action of the \emph{extended Morava stabiliser group}; see \cite[\textsection5]{ec2} and \cite[\textsection7]{gh04}. It is not clear from these techniques that these account for all $\E_\infty$-endomorphisms of $E_n$; to prove this requires a dash of chromatic homotopy theory as done in \cite[\textsection5]{ec2}.\\

Regardless, for any finite subgroup $F$ of the automorphism group of the pair $(R_0, {\G}_0)$, we can also model $E_n^{hF}$ using the fact that $\O^\top_\BTn$ is an étale sheaf.

\begin{prop}
	Consider the situation of \Cref{generalbtspectra} and let $F$ be a finite subgroup of the automorphism group of the pair $(R_0,\G_0)$. Then the defining map $\G\colon \Spf \pi_0 R \to \M_\BTn$ factors through the stacky quotient $\Spf \pi_0 R/F$, this factorisation lies in $\C_{A_0}$, and there is a natural equivalence of $\E_\infty$-rings
	\[\O^\top_\BTn(\Spf \pi_0 R/F) \simeq E(R_0, \G_0)^{h F}.\]
\end{prop}

As $\O^\top_\BTn$ is an étale hypersheaf, there is also potential for extensions to profinite subgroups $F$. These ideas are carried out in \cite{rokdevinatzhopkins} to reprove the famous Devinatz--Hopkins equivalence $E_n^{h\GG_n} \simeq \Sph_{\mathrm{K}(n)}$ of \cite{devinatzhopkins}.

\begin{proof}
	The homotopy fixed points $E(R_0, \G_0)^{h F}$ can be modelled by a limit
	\[E(R_0, \G_0)^{h F} = \lim \left( E(R_0, \G_0) \Rightarrow \prod_K E(R_0, \G_0)\Rrightarrow \cdots \right).\]
	Using \Cref{generalbtspectra} and the fact that $\O^\top_\BTn$ is an étale sheaf, we obtain natural equivalences
	\[\lim \left( E(R_0, \G_0) \Rightarrow \prod_K E(R_0, \G_0)\Rrightarrow \cdots \right) \simeq
	\lim \left( \O^\top_\BTn(\Spf \pi_0 R) \Rightarrow \O^\top_\BTn\left(\coprod_K \Spf \pi_0 R\right) \Rrightarrow \cdots \right).\]
	The map of stacks $\Spf \pi_0 R \to \Spf \pi_0 R /F$ is an $F$-torsor by construction, hence we can make the identification
	\[\Spf \pi_0 R \times_{\Spf \pi_0 R/F} \Spf \pi_0 R \simeq \coprod_K \Spf \pi_0 R\]
	and the same goes for higher fibre products. In particular, the $F$-torsor above gives rise the the \v{C}ech nerve
	\[\cdots \Rrightarrow \coprod_K \Spf \pi_0 R \Rightarrow \Spf \pi_0 R \to \Spf \pi_0 R/F\]
	occurring in the limit above. The fact that $\O^\top_\BTn$ is an étale sheaf gives us the desired equivalence
	\[E(R_0, \G_0)^{h F} \simeq \lim \left( \O^\top_\BTn(\Spf \pi_0 R) \Rightarrow \O^\top_\BTn\left(\coprod_K \Spf \pi_0 R\right) \Rrightarrow \cdots \right) \simeq \O^\top_\BTn(\Spf \pi_0 R/F) .\qedhere\]
\end{proof}

Let us mention one explicit example of such a Lubin--Tate $\E_\infty$-ring. 

\begin{mydef}\label{themultplicatiovepdivisiblegroup}
Let $\mu_{p^\infty}^\heartsuit$ denote the \emph{multiplicative $p$-divisible group} over $\Spec \Z$, whose $R$-valued points (for a discrete ring $R$) are defined as
\[\mu_{p^n}^\heartsuit(R)=\left\{ x\in R \,|\, x^{p^n}=1 \right\}.\]
This lifts to a $p$-divisible group $\mu_{p^\infty}$ over $\Spec \Sph$ by \cite[Pr.2.2.11]{ec2}.
\end{mydef}

 The following is a consequence of the above or one can explicitly prove this using Lurie's proof of Snaith's theorem in \cite[\textsection 6.5]{ec2}.

\begin{prop}\label{ktheorycomesfromlurie}
The $p$-divisible group $\mu_{p^\infty}^\heartsuit$ over $\Spf \Z_p$ defines an object of $\C_{\Z_p}$ (for $n=1$) and there is a natural equivalence of $\E_\infty$-rings $\O^\top_\BTone(\mu_{p^\infty}^\heartsuit)\simeq \KU_p$. Moreover, the associated map $\Spf \Z_p \to \M_\BTone$ factors through $\Spf \Z_p/C_2$, this factorisation lies in $\C_{\Z_p}$, and there is an equivalence of $\E_\infty$-rings $\O^\top_\BTone(\Spf \Z_p /C_2) \simeq \KO_p$.
\end{prop}


\subsection{Topological modular forms}\label{tmfsection}

Another exciting application of \Cref{maintheorem} is to construct the $\E_\infty$-ring $\TMF$ of \emph{periodic topological modular forms}. Of course, this also uses the ideas of Lurie from \cite{ec2} and \cite{lurieecsurveyname}, but reinterpreting $\TMF_p$ as a section of $\O^\top_\BTtwo$ yields additional endomorphisms to those previously known; the case of stable Adams operations is outlined in \Cref{adamsoperationsection}, and further constructions for $\TMF$ and $\Tmf$ appear in \cite{heckeontmf} and \cite{adamsontmf}, respectively.

\begin{prop}\label{pinftyliesinA}
The map $[p^\infty]\colon {\M}_{\Ell}^\heartsuit\to {\M}_{\BTtwo}^\heartsuit$ base changed over $\Spf \Z_p$ lies inside $\C_{\Z_p}$.
\end{prop}

\begin{proof}
Using \Cref{simplercrierion}, we only need to show that the map $[p^\infty]$ above is formally \'{e}tale inside $\P(\Aff^\heartsuit)$ and that ${\M}_{\Ell}^\heartsuit$ is finitely presented over $\Spf \Z_p$. The former follows from \Cref{classicalserretatetheorem}; a consequence of the classical Serre--Tate theorem. The latter follows from \cite[Th.13.1.2]{olsson}, which states that $\M^\heartsuit_\Ell$ is locally of finite presentation over $\Spec \Z$, the fact that this adjective is stable under base change, and the fact that ${\M}_{\Ell}^\heartsuit\times \Spf\Z_p$ is qcqs.
\end{proof}

As promised in the introduction, we should relate $\O^\top_\BTtwo$ to a more classical object:

\begin{mydef}
Let $\O^\top$ denote the Goerss--Hopkins--Miller sheaf of $\E_\infty$-rings on the small \'{e}tale site $\DM^\et_{/\M_\Ell^\heartsuit}$ of $\M_\Ell^\heartsuit$; see \cite[Th.7.0.1]{ec2} or \cite[Th.1.2]{bourbakigoerss} for a version over the compactification of $\M_\Ell^\heartsuit$. We define $\TMF = \O^\top(\M_\Ell)$.
\end{mydef}

\begin{theorem}\label{tmfcomesfromBTptwo}
The sheaf $\O^\top_\BTtwo$, pulled back to $\M_{\Ell}^\heartsuit\times \Spf \Z_p$, is equivalent to the $p$-completion of $\O^\top$. In particular, there is an equivalence of $\E_\infty$-rings
\[\O^\top_\BTtwo\left([p^\infty]\colon {\M}_{\Ell}^\heartsuit \times \Spf \Z_p\to {\M}_{\BTtwo, \Z_p}^\heartsuit\right)\simeq \TMF_p.\]
\end{theorem}

Most of the details for this proof can be found in \cite[Th.7.0.1]{ec2} which proves a similar integral statement. A full proof of this theorem can be found in \cite[Cor.1.6]{heckeontmf}, and an integral version of this appears as \cite[Th.1.9]{heckeontmf}. The important step is to realise that $\O^\top$ is uniquely determined up to homotopy as a sheaf on the small étale site of $\M_\Ell$ by four conditions which mirror those conditions in \Cref{maintheorem}; see \cite[Th.A]{uniqueotop}.\\

Let us also mention a few variations on $\TMF$ that one can obtain from $\O^\top$.

\begin{mydef}\label{levelstructuresonellipticcurbves}
There exist moduli functors $\CAlg^\heartsuit\to \Spc$ denoted as $\M_\Ga$ for each congruence subgroup $\Ga\leq \SL_2(\Z)$. Of particular interest are $\Ga=\Ga(n)$, $\Ga_1(n)$, and $\Ga_0(n)$, which yield moduli stacks $\M(n)$, $\M_1(n)$, and $\M_0(n)$ for each $n\geq 1$. These are defined in \cite[\textsection 3]{km}, and they sit in a commutative diagram in $\P(\Aff^\cn)$
\[\begin{tikzcd}
{\M(n)}\ar[r]\ar[rd]	&	{\M_1(n)}\ar[d]\ar[r]	&	{\M_0(n)}\ar[dl]	\\
				&	{\M_\Ell}			&
\end{tikzcd}\]
where all the transformations above are some kind of forgetful functor. Moreover, by \cite[Th.3.7.1]{km}, we see that when working over $\Spec \Z[\frac{1}{n}]$, all of the morphisms above are finite \'{e}tale. Using these maps one then defines the $\E_\infty$-rings $\TMF(\Ga)=\O^\top(\M_\Ga)$ called \emph{periodic topological modular forms with level structure}. These $\E_\infty$-rings are naturally $\E_\infty$-$\TMF[\frac{1}{n}]$-algebras for $\Ga=\Ga(n)$, $\Ga_1(n)$, or $\Ga_0(n)$.
\end{mydef}

One can also use these moduli stacks with level structure together with \Cref{tmfcomesfromBTptwo} to define stable Hecke operators on $\TMF$. This is already explored in \cite[\textsection11.2]{taf} and when working with $6$ inverted, such operations were originally defined and studied by Baker; see \cite{bakerhecketwo} and \cite{bakerhecke}. Higher categorical and integral refinements appear in \cite{heckeontmf}.\\

Another prominent use of the extra functoriality on $\TMF_p$ afforded by $\O^\top_\BTtwo$ appears in the work of Behrens constructing his $Q(\ell)$ spectra, which form resolutions of the $K(2)$-local sphere; see \cite{ktwospheremark} for the original definition and \cite[Con.1.12]{heckeontmf} for an (almost) integral reformation. Stable Adams operations on $\TMF$ also make an appearance in \cite[\textsection7]{boss}.

\begin{remark}\label{nocompactifcation}
	There are difficulties applying Lurie's theorem to the compactification $\overline{\M}_\Ell$ of $\M_\Ell$. Firstly, as this compactification classifies generalised elliptic curves, and such curves do not necessarily have a group structure, it is not immediately clear how to define a map $\overline{\M}_\Ell \to \M_\BT$. If we took the $p$-divisible group associated with an elliptic curve to be the $p$-divisible group of its smooth locus, then the $p$-divisible group associated with a cuspidal elliptic curve would be $\mu_{p^\infty}$ and have height one. In other words, this map $\overline{\M}_\Ell$ would not be flat and would not factor through $\M_\BTtwo$. There are methods to extending the ideas of this article to $\Tmf$, the global sections of the Goerss--Hopkins--Miller sheaf $\O^\top$ over $\overline{\M}_\Ell$, using obstruction theory, see \cite{adamsontmf} and \cite{realspectra}, and using spectral algebraic geometry, see \cite[\textsection4.3]{lurieecsurveyname} and \cite{globaltate}.
\end{remark}


\subsection{Topological automorphic forms}\label{tafsection}
The first examples of new cohomology theories constructed with \Cref{maintheorem} come from Behrens--Lawson \cite{taf}. The main idea is that the Serre--Tate theorem, which was vital in our construction of $\TMF_p$ from $\O^\top_\BTtwo$, actually applies to the moduli stack of dimension $g$ abelian varieties for any $g\geq 1$; the $g=1$ case recovers the moduli stack of elliptic curves. A new problem now arises: we need our $p$-divisible groups to be of dimension $1$, and then and only then can they have an orientation. To obtain a 1-dimensional $p$-divisible group from an abelian variety $A$ of dimension $g\geq 2$, one needs more structure on $A$ to split its associated $p$-divisible group into one of dimension $1$ and another of dimension $g-1$ (which we forget about). This comes in the form of\textbf{p}olarisations, \textbf{e}ndomorphisms, and \textbf{l}evel structure, leading us \emph{PEL-Shimura varieties}; for a full introduction to the subject and the intended application to stable homotopy theory, see \cite{taf}. What appears below is simply a restatement of \cite{taf} and \cite{handbooktmf}.

\begin{notation}\label{piecesoftaf}
Fix an integer $n\geq 1$. Let $F$ be a quadratic imaginary extension of $\Q$, such that $p$ splits as $u\overline{u}$. Write $\O_F$ be the ring of integers of $F$. Let $V$ be an $F$-vector space of dimension $n$ equipped with a $\Q$-valued nondegenerate Hermitian alternating form of signature $(1,n-1)$. Finally, choose a $\O_F$-lattice $L$ in $V$ such that the alternating form above takes integer values on $L$ and makes $L_{(p)}$ self-dual.
\end{notation}

\begin{mydef}
Write $\X_{V,L}$ for the formal Deligne--Mumford stack over $\Spf \Z_p$ (of \cite[Th.6.6.2]{taf} with $K^p=K^p_0$) where a point in $\X_{V,L}(S)$ for a locally Noetherian formal scheme $S$ over $\Spf \Z_p$, is a triple $(A,i,\lambda)$ where $A$is an abelian scheme over $S$ of dimension $n$, $\lambda\colon A\to A^\vee$ is a polarisation (principle at $p$), with Rosati involuation $\dagger$ on $\Endo(A)_{(p)}$, and $i\colon \O_{F, (p)}\to \Endo(A)_{(p)}$ is an inclusion of rings satisfying $i(\overline{z})=i(z)^\dagger$. These triples have to satisfy two conditions assuring they are locally modelled by $V$ and $L$; see \cite[\textsection6.7]{handbooktmf}.
\end{mydef}

In the situation above, the splitting $p=u\overline{u}$ induces a splitting of $p$-divisible groups
\[A[p^\infty]\simeq A[u^\infty]\oplus A[\overline{u}^\infty]\]
and our assumptions on $(A,i,\lambda)$ ensure that $A[u^\infty]$ is a 1-dimensional $p$-divisible group of height $n$. This yields a morphism of stacks $[u^\infty]\colon\X_{V,L}\to {\M}_{\BTn}^\heartsuit$ which sends $(A,\lambda, i)$ to $A[u^\infty]$.

\begin{prop}\label{PELShimuracanbefoundinA}
Given $V$ and $L$ from \Cref{piecesoftaf}, then the morphism $[u^\infty]\colon \X_{V,L}\to {\M}_{\BTn}^\heartsuit$ base changed over $\Spf \Z_p$ is an object of $\C_{\Z_p}$.
\end{prop}

\begin{proof}
\Cref{simplercrierion} reduces us to show that $[u^\infty]$ is formally \'{e}tale inside $\P(\Aff^\cn)$ and that $\X_{V,L}$ is of finite presentation over $\Spf \Z_p$. The first statement, that $[u^\infty]$ is formally \'{e}tale, follows straight from the definitions of a formally \'{e}tale morphism and \cite[Th.7.3.1]{taf}, which itself is a consequence the classical Serre--Tate theorem and an analysis of $\X_{V,L}$. We now use \cite[Cor.7.3.3]{taf} to see $\X_{V,L}$ is of locally finite presentation over $\Spf \Z_p$, so it suffices to show now that $\X_{V,L}$ is qcqs. To do this, we first use \cite[Th.6.6.2]{taf}, which states that $\X_{V,L}$ has an \'{e}tale cover by a quasi-projective scheme. As a quasi-projective formal scheme $X$ is separated and qc, we see $X$ itself has a Zariski cover by an affine formal scheme $\Spf B$, meaning $\X_{V,L}$ has an \'{e}tale cover by $\Spf B$. By \Cref{etalehypercoversandsuch}, this implies $\X_{V,L}$ is qcqs.
\end{proof}

We can now define the spectra of topological automorphic forms following Behrens--Lawson \cite[\textsection8.3]{taf}.

\begin{mydef}\label{tafworks}
Let $V$ and $L$ be as in \Cref{piecesoftaf}. Define the $\E_\infty$-ring of \emph{topological automorphic forms}
\[\TAF_{V,L}=\O^\top_\BTn\left(\X_{V,L}\xrightarrow{[u^\infty]} {\M}_{\BTn}^\heartsuit\times \Spf\Z_p \right).\]
\end{mydef}

As mentioned in \Cref{comparison}, as there is no construction of the sheaf $\O^\top_\BTn$ of Lurie's theorem in \cite{taf}, so we will not explicitly compare our two definitions. In \Cref{PELShimuracanbefoundinA} above, we have only checked that \Cref{maintheorem} applies to the stacks that Behrens--Lawson use to define $\TAF_{V,L}$ and then applied Lurie's theorem as done in \cite{taf}.\\

As with topological \emph{modular} forms (\Cref{levelstructuresonellipticcurbves}), we can also define variants of $\TAF_{V,L}$ which incorporate level structures. Such extra structure can then be used to define restriction maps, transfers, and Hecke operators on $\TAF_{V,L}$; see \cite[\textsection11]{taf}.


\subsection{Stable Adams operations}\label{adamsoperationsection}
The next example exploits the intrinsic functorality of the sheaf $\O^\top_\BTn$.

\begin{mydef}\label{stableadamsopsdef}
Let $k=(k_1,k_2,\ldots)$ be a $p$-adic integer and $\G$ be a $p$-divisible group over an arbitrary scheme (or stack) $S$. Write $[k]\colon \G\to \G$ for the endomorphism of $\G$ given on $p^n$-torsion by the $k_n$-fold multiplication $[k_n]\colon \G_n\to \G_n$. These assemble to an endomorphism of $\G$ as the sequence $(k_1,k_2,\ldots)$ represents a $p$-adic integer and the closed immersions $\G_n\to \G_{n+1}$ witness the equality $\G_n=\G_{n+1}[p^n]$. If $k$ is a unit inside $\Z_p$ then each $[k_n]$ is an isomorphism of finite flat group schemes on $S$, hence $[k]$ is an automorphism of $\G$. If $\G$ defines a morphism $S\to {\M}_{\BTn}^\heartsuit$ inside $\C_{A_0}$ and $k\in \Z_p^\times$, then write
\[\psi^k = [k]^\ast\colon \O^\top_\BTn(\G)\to \O^\top_\BTn(\G)\]
for the induced endomorphism of $\E_\infty$-rings. These are the ($p$-adic) \emph{stable Adams operations} $\O^\top_\BTn(\G)$; we will justify this name shortly.
\end{mydef}

Many properties expected of Adams operations are formal.

\begin{prop}
Let $l,k$ be two units in $\Z_p$, $\G$ be an object of $\C_{A_0}$, and write $\EE=\O^\top_\BTn(\G)$. Then $\psi^1$ is homotopic to the identity map on the $\E_\infty$-ring $\EE$, and the maps of $\E_\infty$-rings $\psi^l\psi^k$ and $\psi^{l k}$ on $\EE$ are homotopic.
\end{prop}

The homotopy $H$ between $\psi^{l}\psi^k$ and $\psi^{l k}$ above are coherent in the following sense: if $j$ is another $p$-adic unit, then the homotopy between $\psi^{j}\psi^l\psi^k$ and $\psi^{jl k}$ factors through $H$. This follows straight from the fact that $\O^\top_\BTn\colon \C_{A_0}^\op\to \CAlg$ is first and foremost a functor of $\infty$-categories, and the calculations $[l][k]=[l k]$ hold up to equality in $\C_{A_0}$. This is made precise for $\TMF$ in \cite[Th.D]{heckeontmf} and for $\tmf$ in \cite[Ths.C \& D]{realspectra}.

\begin{proof}
As these facts hold for $[k]$ in $\C_{A_0}$ and $\O^\top_\BTn$ is a functor, we obtain the result.
\end{proof}

Using the information we already have at hand, we can calculate $[k]^\ast$ on the homotopy groups of the $\E_\infty$-rings $\O^\top_\BTn(\G)$ over affine objects of $\C_{A_0}$.

\begin{prop}\label{admscals}
Let $k$ be a unit in $\Z_p$ and $\G$ be a $p$-divisible group defining an affine object in $\C_{A_0}$. Then for each integer $n$, we have the equality of morphisms of $\Z_p$-modules
\[[k]^\ast=k^n\colon \pi_{2n} \O^\top_\BTn(\G)\to \pi_{2n} \O^\top_\BTn(\G).\]
\end{prop}

\begin{proof}
Using \Cref{maintheorem}, we see that $\pi_{2n}\O^\top_\BTn(\G)$ is naturally isomorphic to the line bundle $\omega^{\otimes n}_{\G}$ over $\pi_0 \O^\top_\BTn(\G)=B$. It then suffices to calculate the $n=1$ case. As $\omega_{\G}$ is the dualising line for the identity component $\G^\circ$ of $\G$, the $B$-module $\omega_{\G}$ is naturally equivalent to the dual of the Lie algebra $\Lie(\G^\circ)$. We will now calculate $[k]^\ast$ on this Lie algebra. This is quite elementary, but let us recall some details. The question can be answered by localising at each minimal ideal $\m$ of $B$ containing its ideal of definition $J$, and over $B_\m$ the $1$-dimensional formal group $\G^\circ$ has coordinate $t$ and an associated formal group law $G$---the choice of coordinates forms a line bundle over $B_\m$ and line bundles over local rings are trivial; see \cite[\textsection 2]{goerssquasicoherent}. Assume $B$ is local then. If $k$ is an integer, can write $[k]$ on $B\llbracket t\rrbracket$, the global sections of $\G^\circ$ using the coordinate $t$, as the composite
\begin{equation}\label{multiplicationusingfgl}[k]\colon B\llbracket t\rrbracket \xrightarrow{c_k} B\llbracket t_1,\ldots, t_k\rrbracket \xrightarrow{\mu} B\llbracket t\rrbracket\end{equation}
where the first map is the comultiplication on $B\llbracket t\rrbracket$ induced by $G$ and the second is the completed multiplication map. As $c_k(t)\equiv t_1+\cdots + t_k$ modulo terms of higher degree, then $[k](t)\equiv kt$ modulo $t^2$. Finally, the Lie algebra $\Lie(\G^\circ)$ can be written as a Zariski tangent space
\[\Lie(\G^\circ)\simeq \Hom_{\Mod_{B}}(tB\llbracket t\rrbracket/(tB\llbracket t\rrbracket)^2, B).\]
It is now clear that $[k]^\ast\colon \Lie(\G^\circ)\to \Lie(\G^\circ)$ is simply multiplication by $k$ if $k$ is an integer. For a general $p$-adic unit $k$, we approximate $k$ by integers using its $p$-adic expansion, then take the limit.
\end{proof}

\begin{remark}
	In \Cref{admscals} we identified the action of $[k]^\ast$ on the homotopy groups of the affine sections of $\O^\top_\BTn$. Given a general object $\G_0\colon \x_0 \to \M_n^\heartsuit$ of $\C_{A_0}$, there is a descent spectral sequence
	\[E_2 = H^s(\x_0, \omega_{\G_0}^{\otimes t}) \implies \pi_{2t-s} \O^\top_n(\G_0).\]
	One can reinterpret \Cref{admscals} as giving an action of $[k]^\ast$ on this $E_2$-page; this is detailed for various operations on $\TMF$ in \cite[Cor.2.12]{heckeontmf}. This idea is used in \cite[\textsection2]{adamsontmf} to compute the Adams operations on $\pi_\ast \Tmf$. For more on this spectral sequence, see \cite{smfcomputation}.
\end{remark}

The following justifies why we call the operations $[k]^\ast$ stable Adams operations.

\begin{prop}\label{adamsoperationsareeasy}
For integers $k$ not divisible by $p$, the map of $\E_\infty$-rings $[k]^\ast\colon \KU_p\to \KU_p$ is homotopic to classical stable Adams operation $\psi^k$ of \cite[\textsection3.2]{atiyahpoweroperations}.
\end{prop}

\begin{proof}
By restricting ourselves to the case of an integer $k$ not divisible by $p$, we have assured that $[k]\colon \mu_{p^\infty}^\heartsuit\to \mu_{p^\infty}^\heartsuit$ is an automorphism of $p$-divisible groups. Let us write $\EE=\O^\top_\BTone(\mu_{p^\infty}^\heartsuit)$. We claim that $[k]^\ast$ can be calculated on the universal line bundle over $\complexproj^\infty$ using just the algebraic geometry of $\widehat{\G}_m$. By \Cref{ktheorycomesfromlurie}, or more appropriately \cite[Th.6.5.1]{ec2}, the equivalence of $\E_\infty$-rings $\rho\colon \EE\simeq \KU_p$ sends the canonical complex orientation $x_\EE$ of $\EE$ to the usual complex orientation $x_\KU$ of $\KU_p$. We obtain orientations (now in the sense of \Cref{orientationsofformalgroupsetc}) $e_\EE$ and $e_\KU$ of the formal multiplicative group $\widehat{\G}_m$ over $\EE$ and $\KU_p$, respectively, (\cite[Ex.4.3.22]{ec2}) such that $\rho(e_\EE)=e_\KU$. As these orientations of $\widehat{\G}_m$ determine morphisms from the associated Quillen formal group to $\widehat{\G}_m$ (\cite[Pr.4.3.23]{ec2}) and $\rho(e_\EE)=e_\KU$, we obtain the commutative diagram of equivalences of formal groups over $\Z_p$ courtesy of \cite[Pr.4.3.23]{ec2}
\[\begin{tikzcd}
{\widehat{\G}_{\KU_p}^{\QQ_0}}\ar[rr, "{\rho^\ast}"]\ar[rd]	&&	{\widehat{\G}_{\EE}^{\QQ_0}}\ar[ld]	\\
	&	{\widehat{\G}_{m, \Z_p}.}	&
\end{tikzcd}\]
Focusing on $\KU_p$ now, we can rewrite the above diagram of equivalences of formal groups over $\Spf \Z_p$
\begin{equation}\label{iamstilloffbyalittle} \Spf \KU_p^0(\complexproj^\infty)\xrightarrow{\simeq} \Spf \Z_p\llbracket t\rrbracket=\widehat{\G}_m \times \Spf \Z_p. \end{equation}
We know that $[k]$ acts by taking $k$-fold multiplication, which on the multiplicative formal group is an operation represented by the map of rings
\[[k]^\ast\colon\Z_p\llbracket t\rrbracket\xrightarrow{c^k} \Z_p\llbracket t_1,\ldots, t_k\rrbracket\xrightarrow{\mu} \Z_p\llbracket t\rrbracket,\qquad t\mapsto (t+1)^k-1.\]
Recall from (\ref{multiplicationusingfgl}) that the first map is the $k$-fold iteration of the comultiplication, itself given by the formula
\[\Z_p\llbracket t\rrbracket\to \Z_p\llbracket x,y\rrbracket,\qquad t+1\mapsto xy+x+y+1=(x+1)(y+1).\]
The second map above is the completed multiplication map. As the map (\ref{iamstilloffbyalittle}) induces a map of adic rings sending $t$ to $\be x_\KU$, we then obtain the same formulae for $[k]^\ast$ in $\KU_p^0(\complexproj^\infty)$
\[[k]^\ast(\be x_\KU)=(\be x_{\KU}+1)^k-1.\]
As $\be x_\KU\in{\KU}_p^0(\complexproj^\infty)$ is represented by $[\xi_1]-1$ one obtains the equalities
\[[k]^\ast([\xi_1])=[k]^\ast(\be x_\KU+1)=[k]^\ast(\be x_\KU)+1=(\be x_\KU+1)^k-1+1=[\xi_1^{\otimes k}].\]
It follows that for any finite space $X$ and any complex line bundle $\L$ over $X$ with corresponding map $g\colon X\to \complexproj^\infty$, the inherent naturality of $[k]^\ast$ gives us the formula
\[[k]^\ast([\L])=[k]^\ast([g^\ast \xi_1])=g^\ast([k]^\ast(\xi_1))=g^\ast [\xi_1]^k=[\L^{\otimes k}].\]
It follows from \cite[Pr.3.2.1(3)]{atiyahpoweroperations} that the operations $[k]^\ast$ on $\KU_p^0(X)$ are the Adams operations $\psi^k$ as maps of cohomology theories. To lift this statement from one about cohomology theories to one about the spectra that represent them, we now show there are no phantom maps of spectra $\KU_p\to \KU_p$, as this is the only obstacle to the fully faithfulness of the functor
\[\h\Sp\to \CohomTh\qquad E\mapsto E^\ast(-)\]
where $\CohomTh$ denotes the 1-category cohomology theories on finite spaces; see \cite[Cor.2.15]{hoveysticklandMoravalocal} or \cite[Lec.17]{lurielecturenotes}. As $\KU_p$ represents an even periodic Landweber exact cohomology theory, it follows that there exist no phantom endomorphisms of $\KU_p$.
\end{proof}

\begin{remark}
There is a $C_2$-action on the sections of $\O^\top_\BTn$ coming from the inversion action on $p$-divisible groups, ie, coming from $\psi^{-1}$. Any $C_2$-action on an $\E_\infty$-ring $\EE$ can be used to upgrade $\EE$ to a genuinely commutative $C_2$-ring spectrum (the kind with norms); see \cite[Th.2.4]{lennartandhill}. When $p=2$, this has interesting results, for example, the $C_2$-structure on sections of $\O^\top_\BTn$ can be used to obtain a $C_2$-equivariant refinement of part 1 of \Cref{maintheorem}: the complex orientability and Landweber exactness of affine sections of $\O^\top_\BTn$ can be upgraded to \emph{Real orientability} and \emph{Real Landweber exactness} \`{a} la \cite[\textsection3]{lennartandhill}. This essentially follows from the \emph{regular homotopy fixed point spectral sequences} of \cite{tmfwls}, the descent theory developed by Lurie in \cite[\textsection6]{ec2}, and the analogous result of Hahn--Shi \cite{hahnshi} for Lubin--Tate spectra.
\end{remark}

\begin{remark}
Let $p$ be an odd prime. Using the multiplicative lift $\F_p^\times\to\Z_p^\times$, which sends $d$ to the limit of the Cauchy sequence $\{d^{p^n}\}_{n\geq 0}$, one obtains an action of $\F_p^\times\simeq C_{p-1}$ on any sections of $\O^\top_\BTn$. In particular, for any $\G$ in $\C_{A_0}$ (it need not be just an affine object), then the $\E_\infty$-ring $\EE=\O^\top_\BTn(\G)$ has an $\E_\infty$-$\F_p^\times$-action, and the homotopy fixed points $\EE^{h\F_p^\times}$ split off a summand of $\EE$ using the idempotent map
\[\frac{1}{p-1}\sum_{d\in \F_p^\times\subseteq \Z_p^\times}\psi^d\colon \EE\to \EE.\]
In particular, if $\EE=\KU_p$, this summand is known as the \emph{periodic Adams summand} $\mathrm{L}$. The case $\EE=\TMF_p$ is expanded upon in \cite[\textsection3.1]{adamsontmf}.
\end{remark}

\appendix
\section{Appendix on formal spectral Deligne--Mumford stacks}\label{appendix}
Throughout this article we have used basic properties of formal spectral Deligne--Mumford stacks that are not explicitly contained in \cite{sag} (at least not obviously to the author), so we have arranged this appendix to prove these statements. Every single statement below is an extension of a proof of Lurie's in \cite{sag} and the author claims no originality for the ideas below.

\subsection*{Truncations}\label{truncationssubsection}
The following is a generalisation of \cite[Pr.1.4.6.3]{sag} to formal spectral Deligne--Mumford stacks; we will even use the same proof and notation.

\begin{prop}\label{truncatationsoffspdmstacksaregood}
Let $\x=(\X, \O_\x)$ be a locally Noetherian formal spectral Deligne--Mumford stack. For each $n\geq 0$, the object $\tau_{\leq n}\x=(\X, \tau_{\leq n} \O_\x)$ is a locally Noetherian formal spectral Deligne--Mumford stack. Moreover, for every $(\Y, \O_\Y)$ inside $\infty\Top_{\CAlg}^{\sHen}$, if $\O_\Y$ is connective and $n$-truncated, then the canonical map $\tau_{\leq n}\x\to \x$ induces an equivalence
\[\Map_{\infty\Top_{\CAlg}^{\sHen}}((\Y,\O_\Y), \tau_{\leq n}\x)\to \Map_{\infty\Top_{\CAlg}^{\sHen}}((\Y,\O_\Y),\x).\]
\end{prop}

\begin{proof}
The first half of the proof of \cite[Pr.1.4.6.3]{sag} applies \emph{mutatis mutandis}. That is, by copying that proof we see that for every strictly Henselian spectrally ringed $\infty$-topos $(\Y, \O_\Y)$ which is connective and $n$-truncated, the canonical map

\begin{equation}\label{tobereferto}\Map_{\infty\Top_{\CAlg}^{\sHen}}((\Y,\O_\Y), \tau_{\leq n}\x)\to \Map_{\infty\Top_{\CAlg}^{\sHen}}((\Y,\O_\Y),\x)\end{equation}
is an equivalence of spaces. Hence, we are left to show that $\tau_{\leq n} \x=(\X, \tau_{\leq n}\O_\x)$ is a locally Noetherian formal spectral Deligne--Mumford stack. By \cite[Prs.8.1.3.3 \& 8.4.2.7]{sag}, being a formal spectral Deligne--Mumford stack and being locally Noetherian are local conditions, hence we may assume $\x=\Spf A$ for a complete Noetherian adic $\E_\infty$-ring $A$. Set $B=\tau_{\leq n} A$, equipped with the same topology as $A$ induced by $I\subseteq \pi_0 A$ using the isomorphism $\pi_0 A\simeq\pi_0 B$. We need to show $\Spf B$ is connective, $n$-truncated, and equivalent to $\tau_{\leq n}\x$ as a spectrally ringed $\infty$-topos.\\

By \cite[Pr.8.1.1.13]{sag}, we see $\Spf B=(\X_{\Spf B}, \O_{\Spf B})$ is connective. For $n$-truncatedness, one can argue as follows: for affine objects $U$ of $\X_{\Spf B}$ we have $\O_{\Spf B}(U)\simeq C^\wedge_I$ for some \'{e}tale $B$-algebra $C$. As $C$ is an \'{e}tale $\E_\infty$-$B$-algebra, then it is almost of finite presentation, and as $B$ is Noetherian (as a truncation of the Noetherian $\E_\infty$-ring $A$), then the spectral Hilbert basis theorem (\cite[Pr.7.2.4.31]{ha}) implies that $C$ is also Noetherian. It then follows from \cite[Cor.7.3.6.9]{sag} that the natural map of $\E_\infty$-$A$-algebras $C\to  C^\wedge_I$ is flat. As the composition
\[B\to C\to C^\wedge_I\simeq \O_{\Spf B}(U)\]
is flat, we see $\O_{\Spf B}(U)$ is $n$-truncated as $B$ is so. The $\infty$-topos $\X_{\Spf B}$ is generated by affine objects under small colimits (\cite[Pr.8.1.3.7]{sag}) and the structure sheaf $\O_{\Spf B}\colon \X_{\Spf B}^\op\to \CAlg$ preserves limits, so it follows that $\O_{\Spf B}(X)$ is $n$-truncated for all $X\in \X_{\Spf B}$, hence $\Spf B$ is $n$-truncated; see \cite[Rmk.1.3.2.6]{sag}. By (\ref{tobereferto}), the natural map $\Spf B\to \Spf A=\x$ factors as
\[\Spf B\xrightarrow{\phi}\tau_{\leq n}\x= (\X, \tau_{\leq n}\O_\x)\to (\X, \O_\x)=\x.\]
Using \cite[Rmk.8.1.1.9]{sag}, we see the map of underlying $\infty$-topoi induced by $\phi\colon A\to \tau_{\leq n}A=B$ is an equivalence
\[\Shv^\et_{\pi_0 B/I}\simeq \Shv_B^\ad\xrightarrow{\phi_\ast} \Shv_A^\ad\simeq \Shv^\et_{\pi_0 A/I},\]
where we used the notation of \cite[Nt.8.1.1.8]{sag}. Under this map, the structure sheaf of $\Spf B$ is sent to the functor
\begin{equation}\label{letsseeifspfbandtruncationsaregood}\phi_\ast \O_{\Spf B}\colon\CAlg_A^\et\to \CAlg^\cn,\qquad D\mapsto (D\otimes_A B)^\wedge_I\simeq(\tau_{\leq n} D)_I^\wedge.\end{equation}
The equivalence above comes from the facts that $A\to D$ is \'{e}tale and a degenerate Tor-spectral sequence calculation; see \cite[Pr.7.2.1.19]{ha}. To see $\phi$ is an equivalence, it therefore suffices to see that (\ref{letsseeifspfbandtruncationsaregood}) is equivalent to $\tau_{\leq n}\O_{\Spf A}$. This is a slight variation on an argument made above. As $D$ is \'{e}tale over the Noetherian $\E_\infty$-ring $A$, then the spectral Hilbert basis theorem implies that $D$ is also Noetherian. It follows straight from the definition that the $\E_\infty$-ring $\tau_{\leq n} D$ is Noetherian, so the natural completion map of $\E_\infty$-$A$-algebras $\tau_{\leq n} D\to \left(\tau_{\leq n} D\right)^\wedge_I$ is flat. This implies that $(\tau_{\leq n} D)^\wedge_I$ is $n$-truncated. As $\tau_{\leq n}(D^\wedge_I)$ is $I$-complete by \cite[Cor.7.3.4.3]{sag}, there is a natural equivalence of $\E_\infty$-$A$-algebras $(\tau_{\leq n} D)^\wedge_I\simeq \tau_{\leq n}(D^\wedge_I)$. Hence $\phi$ is an equivalence of spectrally ringed $\infty$-topoi.
\end{proof}

The following is a formal generalisation of \cite[Cor.1.4.6.4]{sag}:

\begin{cor}\label{adjointworks}
For each integer $n\geq 0$, write $\fSpDM_\lN^{\leq n}$ for the full $\infty$-subcategory of $\fSpDM_\lN$ spanned by those $n$-truncated locally Noetherian formal spectral Deligne--Mumford stacks. The inclusion $\fSpDM^{\leq n}_\lN\hookrightarrow \fSpDM_\lN$ has a right adjoint, given on objects by
\[\x=(\X, \O_\x)\mapsto \tau_{\leq n} \x=(\X, \tau_{\leq n} \O_\x).\]
\end{cor}

\begin{proof}
This follows straight from the universal property of \Cref{truncatationsoffspdmstacksaregood} and the observation and truncations of locally Noetherian formal spectral Deligne--Mumford stacks remain locally Noetherian.
\end{proof}

\begin{cor}\label{tunrcationsdotherightthingonFOP}
Let $\x$ be a locally Noetherian formal spectral Deligne--Mumford stack. Then for any integer $n\geq0$ the truncation $\tau_{\leq n}\x$ and $\x$ represent the same functor on $n$-truncated $\E_\infty$-rings.
\end{cor}

\begin{proof}
Follows straight from \Cref{truncatationsoffspdmstacksaregood}, as $\Spec R$ is a connective $n$-truncated spectrally ringed $\infty$-topos when $R$ is a connective $n$-truncated $\E_\infty$-ring; see \cite[Ex.1.4.6.2]{sag}.
\end{proof}

\subsection*{The fully faithful embedding $\fDM\to \fSpDM$}
Next, we formalise the relationship between the classical and spectral worlds of formal algebraic geometry. Let us begin by defining these categories.

\begin{mydef}\label{classicalformalspectra}
Let $A$ be a discrete adic Noetherian ring with finitely generated ideal of definition $I\subseteq A$, cutting out a closed subset $V\subseteq |\Spec A|$.
\begin{enumerate}
\item Define the topos $\Shv_\Set^\ad(\CAlg_A^\et)$ is the full $\infty$-subcategory of $\Shv_\Set^\et(\CAlg_A^\et)$ spanned by those \'{e}tale sheaves $\FF$ such that if the space $V\times_{|\Spec A|}|\Spec B|$ is empty, then $\FF(B)$ is a point. 
\item One has a sheaf of discrete rings $\O_{\Spec A}$ on $\Shv^\et_\Set(\CAlg_A^\et)$ as in \cite[Df.1.2.3.1]{sag}, which we complete at $I$ to obtain a sheaf $\widehat{\O}$. This sheaf factors through $\Shv_\Set^\ad(\CAlg_A^\et)$ as $\widehat{\O}(B)\simeq B^\wedge_I$ vanishes if whenever the image of $I$ generates the unit ideal of $B$.
\end{enumerate}
Define the ringed topos $\Spf A=(\Shv_\Set^\ad(\CAlg^\et_A), \widehat{\O})$, the \emph{formal spectrum of $A$}, leaving the dependency on the specific topology on $A$ implicit. A \emph{locally Noetherian formal Deligne--Mumford stack} is a ringed topos $\x=(\X, \O_\x)$ such that $\X$ has a cover $U_\al$ such that each ringed topos $\x_{/\U_\al}$ is equivalent (in the 2-category of ringed topoi of \cite[Df.1.2.1.1]{sag}) to $\Spf A_\al$ for some discrete adic Noetherian ring $A_\al$. Write $\fDM$ for the full $2$-category of $1\Top^\sHen_{\CAlg^\heartsuit}$ spanned by locally Noetherian formal Deligne--Mumford stacks.
\end{mydef}

The $\infty$-category of formal spectral Deligne--Mumford stacks $\fSpDM$ is defined similarly; see \Cref{formalspdm} or \cite[Df.8.1.3.1]{sag}.\\

As in \cite[\textsection8]{sag}, when dealing with classical formal Deligne--Mumford stacks, we restrict ourselves to the locally Noetherian case \emph{by definition}, as opposed to the spectral case, when we only add this assumption when we need it. As mentioned in \cite[Warn.8.1.0.4]{sag}, this is due to the incompatibility between completions in the classical and derived worlds.

\begin{remark}
If an adic discrete ring $A$ has a nilpotent ideal of definition, then $\Spf B$ is naturally equivalent to $\Spec B$ by definition. In this way, we can see (Noetherian) affine Deligne--Mumford stacks as affine formal Deligne--Mumford stacks. It then also immediately follows from the definitions that $\DM_\lN$ is a full $2$-subcategory of $\fDM$.
\end{remark}

The following is a clarification of Lurie's \cite[Rmk.1.4.1.5]{sag}.

\begin{construction}\label{discreteintospectralconstruction}
There is a fully faithful embedding of $\infty$-categories from classical ringed topoi to spectrally ringed $\infty$-topoi
\[1\Top_{\CAlg^\heartsuit}\hookrightarrow \infty\Top_\CAlg,\qquad (\X,\O_\X)\mapsto (\Shv(\X), \O).\]
In other words, it associates to a classical Grothedieck topos $\X$ the associated $\infty$-topos $\Shv(\X)$ (using \cite[Pr.6.4.5.7]{htt}) and by \cite[Rmk.1.3.5.6]{sag} we obtain a connective 0-truncated structure sheaf on $\Shv(\X)$, denoted as $\O$. In fact, the essential image of the above embedding is spanned by the spectrally ringed $\infty$-topoi $(\X,\O_\X)$ where $\X$ is 1-localic and $\O_\X$ is connective and 0-truncated.
\end{construction}

By definition \cite[Df.6.4.5.8]{htt}, we see the $\infty$-topoi $\Shv(\X)$ produced by \Cref{discreteintospectralconstruction} are \emph{1-localic}. By \cite[Rmk.1.4.8.3]{sag}, the fully faithful embedding of \Cref{discreteintospectralconstruction} restricts to the full-faithful embedding $\DM\to \SpDM$. Let us show that the same holds for \emph{formal} Deligne--Mumford stacks.

\begin{prop}\label{spectralanddiscrete}
The functor of \Cref{discreteintospectralconstruction}, when restricted to $\fDM$ factors through $\fSpDM$. Moreover, the essential image of this fully faithful functor $\fDM\to \fSpDM$ consists of those locally Noetherian formal spectral Deligne--Mumford stacks $\x=(\X, \O_\x)$ for which the $\infty$-topos $\X$ is 1-localic (\cite[Df.6.4.5.8]{htt}) and the structure sheaf $\O_\x$ is 0-truncated.
\end{prop}

\begin{proof}
The fully faithful functor of \Cref{discreteintospectralconstruction} descends to a fully faithful functor between (not full) $\infty$-subcategories of local topoi
\[1\Top_{\CAlg^\heartsuit}^\sHen\hookrightarrow \infty\Top_\CAlg^\sHen.\]
Let $\x_0=(\X_0, \O_0)$ be a classical formal Deligne--Mumford stack, and write $\x=(\X, \O)$ for the image of $\x_0$ under \Cref{discreteintospectralconstruction}, so $\X=\Shv(\X_0)$. By \cite[Pr.8.1.3.3]{sag}, the property of being a formal spectral Deligne--Mumford stack is a local one, so it suffices to show that there exists a cover $U_\al$ of $\X$ such that each $\x_{/U_\al}$ is in $\fSpDM$. Consider a formal affine cover of $\x_0$ in $1\Top_{\CAlg^\heartsuit}$, so a collection of $U_\al$ inside $\X_0$ such that $\coprod U_\al\to \1_{\X_0}$ is an effective epimorphism and $(\x_0)_{/U_\al}$ is equivalent in $1\Top_{\CAlg^\heartsuit}$ to $\Spf A_\al$. Considering $U_\al$ as a discrete object $V$ of $\X$ (as in \cite[Pr.6.4.5.7]{htt}), then \cite[Lm.1.4.7.7(2)]{sag} states that $\X_{/V}$ is $1$-localic, as $\X$ is 1-localic and $V$ is 0-truncated in $\X$. One then notes the natural equivalences
\[\X_{/V}\xrightarrow{\simeq} \Shv((\X_{/V})^\heartsuit)\simeq \Shv((\X_0)_{/U_\al})\simeq \Shv(\Shv_\Set^\ad(\CAlg^\et_{A_\al}))\xleftarrow{\simeq} \Shv^\ad(\CAlg_{A_\al}^\et).\]
The first equivalence holds as $\X_{/V}$ is 1-localic, the second by identifying $\X_0$ as the underlying discrete objects of $\X$ (and then \cite[Rmk.7.2.2.17]{htt}), the third from the choice of $U_\al$ as an affine object of $\X_0$, and the forth from the fact that affine formal spectral Deligne--Mumford stacks are 1-localic; see \cite[Rmk.8.1.1.9]{sag}. Furthermore, as $\O$ was defined as the sheaf of connective 0-truncated $\E_\infty$-rings on $\X$ associated to the commutative ring object $\O_0$ on $\X_0$, we claim that by \cite[Rmk.1.3.5.6]{sag} the spectrally ringed $\infty$-topos $\x_{/U_\al}$ is equivalent to $\Spf A_\al$. To see this, one notes that $\O(\Spf B)=B^\wedge_I$ for some \'{e}tale morphism $\Spf B\to \Spf A_\al$ in $\X_0\subseteq \X$, and one also has a natural equivalence $\O_{\Spf A_\al}(\Spf B) \simeq B^\wedge_I$ by \cite[Con.8.1.1.10]{sag}. The ``moreover'' statement follows by \cite[Rmk.1.4.1.5]{sag}.
\end{proof}

Combining the functor of points approach with the above, we obtain the following:

\begin{cor}\label{fffffstacks}
The diagram of $\infty$-categories and fully faithful functors commutes
\[\begin{tikzcd}
{\Aff^\heartsuit_\lN}\ar[r, "a"]\ar[d, "c"]	&	{\Aff^\heartsuit_{\ad, \lN}}\ar[r, "b"]\ar[d, "d"]	&	{\fDM}\ar[d, "e"]	&		\\
{\Aff^\cn}\ar[r, "f"]					&	{\Aff^\cn_\ad}\ar[r, "g"]					&	{\fSpDM}\ar[r, "h"]		&	{\P(\Aff^\cn).}
\end{tikzcd}\]
\end{cor}

\begin{warn}
One might want to place $\P(\Aff^\heartsuit)$ in the top-right corner of the diagram above, however, we do not see a functor $\P(\Aff^\heartsuit)\to \P(\Aff^\cn)$ such that the diagram above commutes. Indeed, the right Kan extension mentioned in \Cref{nicefunctorthatareboring} does not commute with the other constructions above by inspection and a left Kan extension would not necessarily preserve sheaves. The existence of the functors $c$, $d$, and $e$ above, are all due to nontrivial theorems of Lurie, and the lack of a similar functor $\P(\Aff^\heartsuit)\to \P(\Aff^\cn)$ indicates one reason why we restrict our attention to (formal) Deligne--Mumford stacks.
\end{warn}

\begin{proof}[Proof of \Cref{fffffstacks}]
The funtors $a$, $b$, $f$, and $g$ are all the inclusions of full $\infty$-subcategories, $c$ and $d$ are the inclusions of $\infty$-subcategories as shown by Lurie (\cite[Pr.7.1.3.18]{ha}), $e$ is \Cref{discreteintospectralconstruction}, and $h$ is the functor of points functor. The diagram commutes as $c$ and $d$ are restrictions of $e$. We then notice the following: by definition, we see that $a$, $b$, $f$, and $g$ are fully faithful; by \cite[Pr.7.1.3.18]{ha}, we see $c$ and hence $d$ are fully faithful; \Cref{spectralanddiscrete} shows $e$ is fully faithful; and he fact that $h$ is fully faithful is the content of \cite[Th.8.1.5.1]{sag}.
\end{proof}

\subsection*{Finiteness and compactness in $\fSpDM$}
Finally, let us discuss finiteness and compactness conditions in $\fSpDM$.

\begin{prop}\label{truncationsofnoetherianfspdmaregood}
Let $\x$ be a locally Noetherian formal spectral Deligne--Mumford stack. Then for any $n\geq 0$ the natural map $\tau_{\leq n}\x\to \x$ admits an $(n+1)$-connective and almost perfect cotangent complex.
\end{prop}

\begin{proof}
These are local conditions, so we may take $\x=\Spf A$ for a complete Noetherian adic $\E_\infty$-ring $A$ with finitely generated ideal of definition $I\subseteq \pi_0 A$. By the Hilbert basis theorem for connective $\E_\infty$-rings (\cite[Pr.7.2.4.31]{ha}) we see $\tau_{\leq n} A$ is almost finitely presented as an $\E_\infty$-$A$-algebra and the cofibre of map $A\to \tau_{\leq n} A$ is $(n+1)$-connective. By \cite[Cor.7.4.3.2]{ha} and \cite[Th.7.4.3.18]{ha}, we then see $L=L_{\tau_{\leq n} A/A}$ is $(n+1)$-connective and almost perfect inside $\Mod_{\tau_{\leq n}A}$. It follows from \cite[Pr.7.3.5.7]{sag} that $L$ is in fact $I$-complete, hence we have a natural equivalence $L_1\simeq L_{\Spf \tau_{\leq n}A/\Spf A}$ by \cite[Df.17.1.2.8]{sag}, and we are done.
\end{proof}

\begin{mydef}\label{qcqsdefinition}
A formal spectral Deligne--Mumford stack $\x=(\X, \O_\x)$ is \emph{quasi-compact} (qc) if the underlying $\infty$-topos $\X$ is quasi-compact, ie, every cover of $\X$ has a finite subcover; see \cite[Df.A.2.0.12]{sag}. A morphism of formal spectral Deligne--Mumford stacks
\[f\colon \x=(\X, \O_\x)\to \y=(\Y, \O_\y)\]
is \emph{qc} if for any qc object $U$ of $\Y$, the pullback $f^\ast(U)$ is qc in $\X$, meaning $\X_{f^\ast U}$ is qc. A morphism of formal spectral Deligne--Mumford stacks is called \emph{quasi-separated} (qs) if the diagonal map $\Delta\colon\y\to \y\times_\x\y$ is qc. We say $\x$ is qs if $\x\to\Spec\Sph$ is qs. 
\end{mydef}

By formal nonsense, we see that qc (and qs) maps are stable under base change; a fact we will use without further reference.

\begin{prop}\label{spfisqc}
Let $A$ be an adic $\E_\infty$-ring. Then $\Spf A$ is qc.
\end{prop}

\begin{proof}
By \cite[Rmk.8.1.1.9]{sag}, we see the underlying $\infty$-topos of $\Spf A$ is equivalent to $\Shv^\et_{\pi_0 A/I}$ where $I$ is a finitely generated ideal of definition for the topology on $\pi_0 A$. As this is the same underlying $\infty$-topos of $\Spec (\pi_0 A/I)$, it follows from \cite[Pr.2.3.1.2]{sag} that $\Spf A$ is qc. 
\end{proof}

The following is a formal generalisation of a special case of \cite[Pr.2.3.2.1]{sag}.

\begin{prop}\label{whatisqs}
Let $\x=(\X, \O_\x)$ be a formal spectral Deligne--Mumford stack. Then the following are equivalent.
\begin{enumerate}
\item $\x$ is qs.
\item For all qc objects $U,V$ of $\X$, the product $U\times V$ in $\X$ is qc.
\item For all affine objects $U,V$ of $\X$, the product $U\times V$ is qc.
\end{enumerate}
\end{prop}

\begin{proof}
It is clear that 1 implies 2 as $U\times V=\Delta^\ast(U,V)$ inside $\X\times\X$, and 2 also implies 1 as the quasi-compact objects of $\X\times\X$ are all of the form $(U,V)$ for $U$ and $V$ quasi-compact in $\X$. \Cref{spfisqc} shows that 2 implies 3. Conversely, for two arbitrary qc objects $U$ and $V$ of $\X$, using the fact they are qc, there exist two effective epimorphisms $U'\to U$ and $V'\to V$ where $U'$ and $V'$ are affine. It then follows that $U\times V$ is qc as there is an effective epimorphism $U'\times V'\to U\times V$ from a qc object of $\X$.
\end{proof}

\begin{cor}\label{affinesareqcqs}
Let $A$ be an adic $\E_\infty$-ring. Then $\Spf A$ is qcqs.
\end{cor}

\begin{proof}
By \Cref{spfisqc}, we see $\Spf A$ is qc, and by \Cref{whatisqs} it suffices to see that for all affine objects $U=\Spf B$ and $V=\Spf C$ inside $\X_{\Spf A}$, that the product $U\times V$ in $\X_{\Spf A}$ is qc. This product can be recognised as the fibre product (\cite[Lm.8.1.7.3]{sag})
\[\Spf B\underset{\Spf A}{\times} \Spf C\simeq \Spf\left(B\underset{A}{\otimes} C\right)^\wedge_I\]
where $I$ is an ideal of definition for the topology on $\pi_0 A$, which is qc by \Cref{whatisqs}.
\end{proof}

\begin{prop}\label{etalehypercoversandsuch}
Let $\x$ be a formal spectral Deligne--Mumford stack. Then $\x$ is qcqs if and only if there exists an \'{e}tale hypercover $\u_\bullet$ of $\x$ such that each $\u_n$ is an affine formal spectral Deligne--Mumford stack for every $n\geq 0$. In particular, the same holds for classical Deligne--Mumford stacks.
\end{prop}

\begin{proof}
First, let us assume $\x$ is qcqs and write $\x=(\X, \O_\x)$ and set $\u_{-1}=\x$. As a formal spectral Deligne--Mumford stack, there exists a collection of affine objects $U_\al$ in $\X$ such that $\coprod_\al U_\al$ cover $\X$, and as $\x$ is qc, this collection can be taken to be finite. As $\X_{/U_\al}\simeq \Spf A_\al$ for some adic $\E_\infty$-ring $A_\al$, we see the fact that $\coprod U_\al$ covers $\X$ is equivalent to the statement that
\[\Spf A_0=\Spf \left(\prod A_\al\right)\simeq \coprod \Spf A_\al\to \x\]
is an \'{e}tale surjection, where we have used the finiteness of the above co/product. Set $\u_0=\Spf A_0$ and $\u_0\to M_0(\u_\bullet^{\leq -1})\simeq \u_{-1}=\x$ to be the \'{e}tale surjection above. One can now inductively construct $\u_\bullet$ by defining $\u_{n+1}$ as something which covers $M_{n+1}(\u_\bullet^{\leq n})$; see the proof of \cite[Pr.A.3.6]{mythesis} for more details. Conversely, assume that $\x$ has an \'{e}tale hypercover $\u_\bullet\to \x$ where each $\u_n$ is affine, which we write as $U_\bullet\to \1$ when considered as objects in $\X$. Given an arbitrary cover $\{V_\al\}_{\al\in I}$ of $\x$, so an effective epimorphism $\coprod V_\al \to \1$, then we can consider the Cartesian square inside $\X$ of the form
\[\begin{tikzcd}
{W}\ar[r]\ar[d]		&	{U_0}\ar[d]	\\
{\coprod_I V_\al}\ar[r]	&	{\1.}
\end{tikzcd}\]
All of the maps above are effective epimorphisms either by assumption or as the class of such maps is stable under pullback; see \cite[Pr.6.2.3.15]{htt}. Products commute with colimits in an $\infty$-topos as colimits in $\infty$-topoi are universal,\footnote{We say that colimits in a presentable $\infty$-category $\C$ are \emph{universal} if pullbacks commute with all small colimits; see \cite[Df.6.1.1.2]{htt}. This holds in an $\infty$-topos due to the $\infty$-categorical version of Giraud's axioms; see \cite[Th.6.1.0.6]{htt}.} hence we have a natural equivalence in $W\simeq \coprod_I W_\al$ in $\X$, where $W_\al=V_\al\times U_0$. As $U_0$ is quasi-compact (as an affine object of $\X$; see \Cref{spfisqc}), we can choose a finite subset of $I$, say $I_0$, such that $\coprod_{I_0} W_\al\to U_0$ is an effective epimorphism. We then consider the commutative diagram inside the $\infty$-topos $\X$
\[\begin{tikzcd}
{\coprod_{I_0}W_\al}\ar[r]\ar[d]	&	{U_0}\ar[d]	\\
{\coprod_{I_0} V_\al}\ar[r]		&	{\1.}
\end{tikzcd}\]
The top and right maps are effective epimorphisms by assumption, and the bottom map is an effective epimorphism by \cite[Cor.6.2.3.12(2)]{htt}, hence $\X$ is qc. To see $\X$ is qs, we look at the Cartesian diagram of formal spectral Deligne--Mumford stacks
\[\begin{tikzcd}
{\u_0}\ar[r, "{\Delta_{\u_0}}"]\ar[d]	&	{\u_0\times \u_0}\ar[d]	\\
{\X}\ar[r, "{\Delta_\X}"]			&	{\X\times\X.}
\end{tikzcd}\]
As $\u_\bullet\to \x$ is an \'{e}tale hypercover, the map $\U_0\times \U_0\to \X\times \X$ is an effective epimorphism. As $\u_0$ is the $\infty$-topos of an affine formal Deligne--Mumford stack, then by \Cref{affinesareqcqs} we see $\u_0$ is qs and the map $\Delta_{\u_0}$ is qc. It follows from \cite[Cor.A.2.1.5]{sag} that $\Delta_{\X}$ is qc; in \emph{ibid}, a qc morphism is called \emph{relatively $0$-coherent}. Hence, $\X$, and therefore $\x$, is qs.
\end{proof}

Let us now show the \emph{formal thickenings} of \cite[\textsection18.2.2]{sag} preserve the adjective qcqs.

\begin{prop}\label{qcqsformalthickenings}
Let $\x_0$ be a qcqs formal spectral Deligne--Mumford stack and $\x_0\to \x$ a formal thickening. Then $\x$ is qcqs.
\end{prop}

\begin{proof}
The adjective qcqs depends only on the underlying $\infty$-topoi, so it suffices to show if that $\x_0\to \x$ is an equivalence of $\infty$-topoi. To see this, consider the \emph{reduction} of a formal spectral Deligne--Mumford stack of \cite[Pr.8.1.4.4]{sag}. From this one obtains the commutative diagram of formal spectral Deligne--Mumford stacks
\[\begin{tikzcd}
{\x_0^\red}\ar[r]\ar[d]	&	{\x_0}\ar[d]	\\
{\x^\red}\ar[r]		&	{\x.}
\end{tikzcd}\]
We know the natural map from the reduction of a formal spectral Deligne--Mumford stack $\x$ back into $\x$ is an equivalence of underlying $\infty$-topoi (by \cite[Pr.8.1.4.4]{sag}), and the underlying $\infty$-topoi of the reduction of a formal thickening is also an equivalence (by \cite[Pr.18.2.2.6]{sag}). Hence the horizontal and the left vertical maps are equivalences of underlying $\infty$-topoi, hence the right vertical map is as well. 
\end{proof}

\newcommand{\Bibkeyhacktwo}[2]{}
\newcommand{\Bibkeyhackthree}[3]{}
\newcommand{\Bibkeyhackfour}[4]{}
\newcommand{\Bibkeyhackfive}[5]{}

\addcontentsline{toc}{section}{References}
\scriptsize{

\bibliography{/Users/jackdavies/Dropbox/Work/references} 

\begin{thebibliography}{{Dav}24b}

\bibitem[{Ada}74]{bluebook}
J.~F. {Adams}.
\newblock {Stable homotopy and generalised homology.}
\newblock {Chicago Lectures in Mathematics. Chicago - London: The University of
  Chicago Press. X, 373 p. \textsterling 3.00}, 1974.

\bibitem[And71]{andrehomology}
Michel Andr{\'e}.
\newblock Homologie des alg{\'e}bres commutatives. ({Homology} of commutative
  algebras).
\newblock Actes {Congr}. internat. {Math}. 1970, 1, 301-308, 1971.

\bibitem[{Ati}67]{atiyahpoweroperations}
Michael~F. {Atiyah}.
\newblock {\(K\)-theory. Lecture notes by D. W. Anderson. Fall 1964. With
  reprints of M. F. Atiyah: Power operations in \(K\)-theory; \(K\)-theory and
  reality.}
\newblock {New York-Amsterdam: W.A. Benjamin, Inc. 166 p.}, 1967.

\bibitem[{Bak}90]{bakerhecketwo}
Andrew {Baker}.
\newblock {Hecke operators as operations in elliptic cohomology.}
\newblock {\em {J. Pure Appl. Algebra}}, 63(1):1--11, 1990.

\bibitem[{Bak}98]{bakerhecke}
Andrew {Baker}.
\newblock {Hecke algebras acting on elliptic cohomology}.
\newblock In {\em {Homotopy theory via algebraic geometry and group
  representations. Proceedings of a conference on homotopy theory, Evanston,
  IL, USA, March 23--27, 1997}}, pages 17--26. Providence, RI: American
  Mathematical Society, 1998.

\bibitem[Bas99]{taqofbasterra}
M.~Basterra.
\newblock Andr{\'e}-{Quillen} cohomology of commutative {{\(S\)}}-algebras.
\newblock {\em J. Pure Appl. Algebra}, 144(2):111--143, 1999.

\bibitem[BB20]{agnestobi}
Agn\`es {Beaudry} and Tobias {Barthel}.
\newblock {Chromatic structures in stable homotopy theory}.
\newblock In {\em Handbook of homotopy theory}, pages 163--220. Boca Raton, FL:
  CRC Press, 2020.

\bibitem[{Beh}06]{ktwospheremark}
Mark {Behrens}.
\newblock {A modular description of the \(K(2)\)-local sphere at the prime 3}.
\newblock {\em {Topology}}, 45(2):343--402, 2006.

\bibitem[{Beh}20]{handbooktmf}
Mark {Behrens}.
\newblock {Topological modular and automorphic forms}.
\newblock In {\em {Handbook of homotopy theory}}, pages 221--261. Boca Raton,
  FL: CRC Press, 2020.

\bibitem[BL10]{taf}
Mark {Behrens} and Tyler {Lawson}.
\newblock {\em {Topological automorphic forms.}}, volume 958.
\newblock Providence, RI: American Mathematical Society (AMS), 2010.

\bibitem[BOSS19]{boss}
M.~{Behrens}, K.~{Ormsby}, N.~{Stapleton}, and V.~{Stojanoska}.
\newblock {On the ring of cooperations for 2-primary connective topological
  modular forms}.
\newblock {\em {J. Topol.}}, 12(2):577--657, 2019.

\bibitem[CD24]{heighttwojatthree}
Christian Carrick and Jack~Morgan Davies.
\newblock Nonvanishing of products in $v_2$-periodic families at the prime $3$.
\newblock \href{arXiv : 2410.02564}{https://arxiv.org/abs/2410.02564}, 2024.

\bibitem[CDvN24]{smfcomputation}
Christian Carrick, Jack~Morgan Davies, and Sven van Nigtevecht.
\newblock The descent spectral sequence for topological modular forms.
\newblock Preprint, {arXiv}:2412.01640 [math.{AT}] (2024), 2024.

\bibitem[CM21]{clausenmathew1}
Dustin {Clausen} and Akhil {Mathew}.
\newblock {Hyperdescent and \'etale \(K\)-theory}.
\newblock {\em {Invent. Math.}}, 225(3):981--1076, 2021.

\bibitem[CS15]{serretatevolII}
Pierre {Colmez} and Jean-Pierre {Serre}, editors.
\newblock {\em {Correspondance Serre -- Tate. Volume II.}}, volume~14.
\newblock Paris: Soci\'et\'e Math\'ematique de France (SMF), 2015.

\bibitem[Dav21]{adamsontmf}
Jack~Morgan Davies.
\newblock {Constructing and calculating Adams operations on dualisable
  topological modular forms}, 2021.
\newblock {arXiv preprint,
  \href{https://arxiv.org/abs/2104.13407}{arXiv:2104.13407}. Accepted in
  Documenta Mathematicae}.

\bibitem[{Dav}22]{mythesis}
Jack~Morgan {Davies}.
\newblock {\em {Stable operations and topological modular forms}}.
\newblock Utrecht University, 2022.
\newblock available at \url{https://dspace.library.uu.nl/handle/1874/422774}.

\bibitem[Dav23]{uniqueotop}
Jack~Morgan Davies.
\newblock Elliptic cohomology is unique up to homotopy.
\newblock {\em J. Aust. Math. Soc.}, 115(1):99--118, 2023.

\bibitem[Dav24a]{heckeontmf}
Jack~Morgan Davies.
\newblock Hecke operators on topological modular forms.
\newblock {\em Advances in Mathematics}, 452:109828, 2024.

\bibitem[{Dav}24b]{realspectra}
Jack~Morgan {Davies}.
\newblock {Uniqueness of real ring spectra up to higher homotopy}.
\newblock {\em Annals of {$K$}-theory}, 9(3):447--473, August 2024.

\bibitem[DFHH14]{tmfbook}
Christopher~L. {Douglas}, John {Francis}, Andr\'e~G. {Henriques}, and
  Michael~A. {Hill}, editors.
\newblock {\em {Topological modular forms. Based on the Talbot workshop, North
  Conway, NH, USA, March 25--31, 2007}}, volume 201.
\newblock Providence, RI: American Mathematical Society (AMS), 2014.

\bibitem[DH04]{devinatzhopkins}
Ethan~S. Devinatz and Michael~J. Hopkins.
\newblock Homotopy fixed point spectra for closed subgroups of the {Morava}
  stabilizer groups.
\newblock {\em Topology}, 43(1):1--47, 2004.

\bibitem[DL25]{globaltate}
Jack~Morgan {Davies} and Sil {Linskens}.
\newblock {On the derived {Tate} curve and global smooth {Tate} {$K$}-theory}.
\newblock In preparation, 2025.

\bibitem[Dri76]{ellipticmodulesi}
V.~G. Drinfel'd.
\newblock Elliptic modules.
\newblock {\em Math. USSR, Sb.}, 23:561--592, 1976.

\bibitem[GH04]{gh04}
P.~G. {Goerss} and M.~J. {Hopkins}.
\newblock {Moduli spaces of commutative ring spectra.}
\newblock In {\em {Structured ring spectra}}, pages 151--200. Cambridge:
  Cambridge University Press, 2004.

\bibitem[Goe08]{goerssquasicoherent}
Paul~G. Goerss.
\newblock Quasi-coherent sheaves on the moduli stack of formal groups, 2008.

\bibitem[{Goe}10]{bourbakigoerss}
Paul~G. {Goerss}.
\newblock {Topological modular forms [after Hopkins, Miller and Lurie]}.
\newblock In {\em {S\'eminaire Bourbaki. Volume 2008/2009. Expos\'es
  997--1011}}, pages 221--255, ex. Paris: Soci\'et\'e Math\'ematique de France
  (SMF), 2010.

\bibitem[Gre23]{rokdevinatzhopkins}
Rok Gregoric.
\newblock The {Devinatz}-{Hopkins} theorem via algebraic geometry.
\newblock {\em Algebr. Geom. Topol.}, 23(7):3015--3042, 2023.

\bibitem[GW10]{hortzwedhorn}
Ulrich {G{\"{o}}rtz} and Torsten {Wedhorn}.
\newblock {\em {Algebraic geometry I. Schemes. With examples and exercises.}}
\newblock Wiesbaden: Vieweg+Teubner, 2010.

\bibitem[HM17]{lennartandhill}
Michael~A. {Hill} and Lennart {Meier}.
\newblock {The \(C_2\)-spectrum \(\mathrm{Tmf}_1(3)\) and its invertible
  modules.}
\newblock {\em {Algebr. Geom. Topol.}}, 17(4):1953--2011, 2017.

\bibitem[{Hop}95]{hopkinsfirsttmficm}
Michael~J. {Hopkins}.
\newblock {Topological modular forms, the Witten genus, and the theorem of the
  cube}.
\newblock In {\em Proceedings of the international congress of mathematicians,
  ICM '94, August 3-11, 1994, Z\"urich, Switzerland. Vol. I}, pages 554--565.
  Basel: Birkh\"auser, 1995.

\bibitem[HS99]{hoveysticklandMoravalocal}
Mark {Hovey} and Neil~P. {Strickland}.
\newblock {\em {Morava \(K\)-theories and localisation.}}, volume 666.
\newblock Providence, RI: American Mathematical Society (AMS), 1999.

\bibitem[HS20]{hahnshi}
Jeremy {Hahn} and XiaoLin~Danny {Shi}.
\newblock {Real orientations of Lubin-Tate spectra}.
\newblock {\em {Invent. Math.}}, 221(3):731--776, 2020.

\bibitem[HT01]{harristaylorshimuravar}
Michael Harris and Richard Taylor.
\newblock {\em The geometry and cohomology of some simple {Shimura} varieties.
  {With} an appendix by {Vladimir} {G}. {Berkovich}}, volume 151 of {\em Ann.
  Math. Stud.}
\newblock Princeton, NJ: Princeton University Press, 2001.

\bibitem[{Ill}71]{illusie}
Luc {Illusie}.
\newblock {\em {Complexe cotangent et d\'eformations. I. (The cotangent complex
  and deformations. I.).}}, volume 239.
\newblock Springer, Cham, 1971.

\bibitem[KM85]{km}
Nicholas~M. {Katz} and Barry {Mazur}.
\newblock {\em {Arithmetic moduli of elliptic curves.}}, volume 108.
\newblock Princeton University Press, Princeton, NJ, 1985.

\bibitem[{Lan}88]{landweberconferencebook}
Peter~S. {Landweber}, editor.
\newblock {\em {Elliptic curves and modular forms in algebraic topology.
  Proceedings of a conference held at the Institute for Advanced Study,
  Princeton, NJ, Sept. 15-17, 1986.}}, volume 1326.
\newblock Berlin etc.: Springer-Verlag, 1988.

\bibitem[LT66]{lubintate}
J.~{Lubin} and J.~{Tate}.
\newblock {Formal moduli for one-parameter formal Lie groups.}
\newblock {\em {Bull. Soc. Math. Fr.}}, 94:49--59, 1966.

\bibitem[{Lur}09a]{lurieecsurveyname}
Jacob {Lurie}.
\newblock {A survey of elliptic cohomology}.
\newblock In {\em {Algebraic topology. The Abel symposium 2007. Proceedings of
  the fourth Abel symposium, Oslo, Norway, August 5--10, 2007}}, pages
  219--277. Berlin: Springer, 2009.

\bibitem[Lur09b]{htt}
Jacob Lurie.
\newblock {\em Higher topos theory}, volume 170 of {\em Ann. Math. Stud.}
\newblock Princeton, NJ: Princeton University Press, 2009.

\bibitem[{Lur}10]{lurielecturenotes}
Jacob {Lurie}.
\newblock {Chromatic homotopy theory (for Math 252x (offered Spring 2010 at
  Harvard))}.
\newblock Available at \url{https://www.math.ias.edu/~lurie/}, 2010.

\bibitem[{Lur}17]{ha}
Jacob {Lurie}.
\newblock Higher algebra.
\newblock Available at \url{https://www.math.ias.edu/~lurie/}, September 2017.

\bibitem[{Lur}18a]{ec1}
Jacob {Lurie}.
\newblock {Elliptic Cohomology I: Spectral Abelian Varieties}.
\newblock Available at \url{https://www.math.ias.edu/~lurie/}, February 2018.

\bibitem[{Lur}18b]{ec2}
Jacob {Lurie}.
\newblock {Elliptic Cohomology II: Orientations}.
\newblock Available at \url{https://www.math.ias.edu/~lurie/}, April 2018.

\bibitem[{Lur}18c]{sag}
Jacob {Lurie}.
\newblock {Spectral Algebraic Geometry}.
\newblock Available at \url{https://www.math.ias.edu/~lurie/}, February 2018.

\bibitem[{Lur}19]{ec3}
Jacob {Lurie}.
\newblock {Elliptic Cohomology III: Tempered cohomology}.
\newblock Available at \url{https://www.math.ias.edu/~lurie/}, 2019.

\bibitem[{Mei}22]{tmfwls}
Lennart {Meier}.
\newblock {Topological modular forms with level structure: decompositions and
  duality}.
\newblock {\em {Trans. Am. Math. Soc.}}, 375(2):1305--1355, 2022.

\bibitem[Mor85]{moravacobordism}
Jack Morava.
\newblock Noetherian localisations of categories of cobordism comodules.
\newblock {\em Ann. Math. (2)}, 121:1--39, 1985.

\bibitem[MRW77]{MRW77}
Haynes~R. {Miller}, Douglas~C. {Ravenel}, and W.~Stephen {Wilson}.
\newblock {Periodic phenomena in the Adams-Novikov spectral sequence.}
\newblock {\em {Ann. Math. (2)}}, 106:469--516, 1977.

\bibitem[MV15]{cubicalhomotopytheory}
Brian~A. {Munson} and Ismar {Voli\'c}.
\newblock {\em {Cubical homotopy theory.}}, volume~25.
\newblock Cambridge: Cambridge University Press, 2015.

\bibitem[{Nau}07]{naumann07}
Niko {Naumann}.
\newblock {The stack of formal groups in stable homotopy theory.}
\newblock {\em {Adv. Math.}}, 215(2):569--600, 2007.

\bibitem[{Ols}16]{olsson}
Martin {Olsson}.
\newblock {\em {Algebraic spaces and stacks.}}, volume~62.
\newblock Providence, RI: American Mathematical Society (AMS), 2016.

\bibitem[Qui69]{quillenstheorem}
D.~Quillen.
\newblock On the formal group laws of unoriented and complex cobordism theory.
\newblock {\em Bull. Am. Math. Soc.}, 75:1293--1298, 1969.

\bibitem[{Qui}70]{quillencotangentcomplex}
Daniel {Quillen}.
\newblock {On the (co)-homology of commutative rings.}
\newblock {Appl. categorical Algebra, Proc. Sympos. Pure Math. 17, 65-87},
  1970.

\bibitem[{Rav}92]{orangebook}
Douglas~C. {Ravenel}.
\newblock {\em {Nilpotence and periodicity in stable homotopy theory.}}, volume
  128.
\newblock Princeton, NJ: Princeton University Press, 1992.

\bibitem[{Rav}04]{greenbook}
Douglas~C. {Ravenel}.
\newblock {\em {Complex cobordism and stable homotopy groups of spheres. 2nd
  ed.}}
\newblock Providence, RI: AMS Chelsea Publishing, 2nd ed. edition, 2004.

\bibitem[{Sil}86]{silvermaneasy}
Joseph~H. {Silverman}.
\newblock {\em {The arithmetic of elliptic curves.}}, volume 106.
\newblock Springer, New York, NY, 1986.

\bibitem[{Sta}]{stacks}
The {Stacks Project Authors}.
\newblock \textit{Stacks Project}.
\newblock \url{https://stacks.math.columbia.edu}.

\bibitem[{Tat}67]{tatepdiv}
J.~T. {Tate}.
\newblock {\(p\)-divisible groups.}
\newblock {Proc. Conf. local Fields, NUFFIC Summer School Driebergen 1966,
  158-183}, 1967.

\bibitem[{Zin}84]{zinkcartiertheorie}
Thomas {Zink}.
\newblock {\em {Cartiertheorie kommutativer formaler Gruppen. Unter Mitarb. von
  Harry Reimann}}, volume~68.
\newblock Teubner, Leipzig, 1984.

\end{thebibliography}
\bibliographystyle{alpha}
}

\end{document}